\documentclass[12pt]{article} 
\usepackage{amsfonts,amsmath}
\usepackage{amssymb,latexsym}
\usepackage[dvips]{epsfig}
\usepackage{amscd} 

\title{Matrix factorizations and link homology} 
\author{ Mikhail Khovanov and Lev Rozansky} 
\date{}

\newtheorem{prop}{Proposition}
\newtheorem{theorem}{Theorem}
\newtheorem{lemma}{Lemma}
\newtheorem{definition}{Definition} 
\newtheorem{corollary}{Corollary}
 
\newcommand{\oplusop}[1]{{\mathop{\oplus}\limits_{#1}}}

\begin{document}
\maketitle
\baselineskip 14pt 

\def\C{\mathbb C}
\def\R{\mathbb R}
\def\Q{\mathbb Q}
\def\Z{\mathbb Z}
\def\l{\lbrace}
\def\r{\rbrace}
\def\o{\otimes}
\def\lra{\longrightarrow}
\def\Hom{\mathrm{Hom}}
\def\hom{\mathrm{hom}}          
\def\Ext{\mathrm{Ext}}
\def\ext{\mathrm{ext}}          
\def\RHom{\mathrm{RHom}}
\def\hmf{\mathrm{hmf}}
\def\HMF{\mathrm{HMF}}
\def\MF{\mathrm{MF}}
\def\Id{\mathrm{Id}}
\def\Tr{\mathrm{Tr}}
\def\mc{\mathcal} 
\def\mf{\mathfrak}
\def\gdim{\mathrm{gdim}}
\def\cF{\mc{F}}
\def\cC{\mc{C}}
\def\sbinom#1#2{\left( \hspace{-0.06in}\begin{array}{c} #1 \\ #2 
\end{array}\hspace{-0.06in} \right)}
\newcommand{\bba}{\mathbf{a}}
\newcommand{\bbb}{\mathbf{b}}
\newcommand{\bracket}[1]{\langle #1 \rangle}
\newcommand{\otimesop}[1]{{\mathop{\otimes}\limits_{#1}}}
\newcommand{\mtwobyone}[2]{\left( \begin{array}{c} #1 \\ #2
   \end{array} \right)} 
\newcommand{\define}{\stackrel{\mbox{\scriptsize{def}}}{=}}
\def\drawing#1{\begin{center}\epsfig{file=#1}\end{center}}

 \begin{abstract} 
For each positive integer $n$ the HOMFLY polynomial of links 
specializes to a one-variable polynomial that can be recovered 
from the representation theory of quantum $sl(n).$ For 
each such $n$ we build a doubly-graded homology theory of links 
with this polynomial as the Euler characteristic. The core of 
our construction utilizes the theory of matrix factorizations, 
which provide a linear algebra description of maximal 
Cohen-Macaulay modules on isolated hypersurface singularities. 
 \end{abstract} 

 \tableofcontents

\section{Introduction} 

\begin{figure} [htb] \drawing{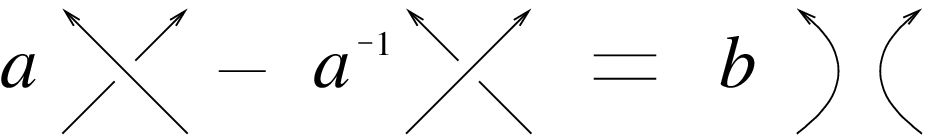}\caption{The  
 HOMFLY skein relation} \label{homfly}
 \end{figure}

The HOMFLY polynomial of oriented links in $\mathbb{R}^3$ is 
uniquely determined by the skein relation in figure~\ref{homfly} 
and its value on the unknot, see [HOMFLY]. 

The specialization $a=q^n$ and $b=q-q^{-1},$ for integer $n,$ 
produces a one-variable link polynomial which 
can be interpreted via representation theory of the quantum 
group $U_q(sl(n))$ if $n$ is positive, see [RT], $U_q(sl(-n))$ 
if $n$ is negative, and $U_q(gl(1|1))$ if $n=0,$ see [KS].
This polynomial is invariant under changing $q$ to $q^{-1}$ 
simultaneously with passing to the 
mirror image of the link; therefore, we don't lose any 
information by restricting to non-negative $n.$ 
We denote this one-variable polynomial by $P_n(L),$ where 
 $L$ is an oriented link; the normalization is 
 $$P_n(\mathrm{unknot})=[n]\define\frac{q^{n}-q^{-n}}{q-q^{-1}}$$ 
if $n>0,$ and $P_0(\mathrm{unknot})=1.$  

For $n=0,1,2,3$ there exists a doubly-graded homology 
theory of links whose Euler characteristic is 
$P_n(L).$ Let us denote this theory by $H_n'(L).$ 

\begin{itemize} 
\item $P_0(L)$ is the Alexander polynomial, and 
$H_0'(L)$ was constructed by Peter Ozsv\'ath and Zolt\'an Szab\'o 
[OS], and, independently, by Jacob Rasmussen [Ra].  
Their theory exists in greater generality and, in particular, 
 encompasses knots in homology spheres. 
\item $P_1(L)=1$ and $H_1'(L)\cong \Z$ for any oriented link 
 $L,$ with $\Z$ in bidegree $(0,0).$ It will be clear 
 subsequently that this is a natural choice for $H_1'(L).$ 
\item $P_2(L)$ is the Jones polynomial; $H_2'(L)$ was defined in 
 [Kh1], and denoted by $\mc{H}(L)$ there. $H_2'(\mathrm{unknot})$ 
 is isomorphic to the integral cohomology ring of the 2-sphere. 
\item  $H_3'(L)$ was constructed in [Kh3]. $H_3'(\mathrm{unknot})$ 
 is isomorphic to the integral cohomology ring of 
 $\mathbb{CP}^2.$  
\end{itemize} 

\begin{figure} [htb] \drawing{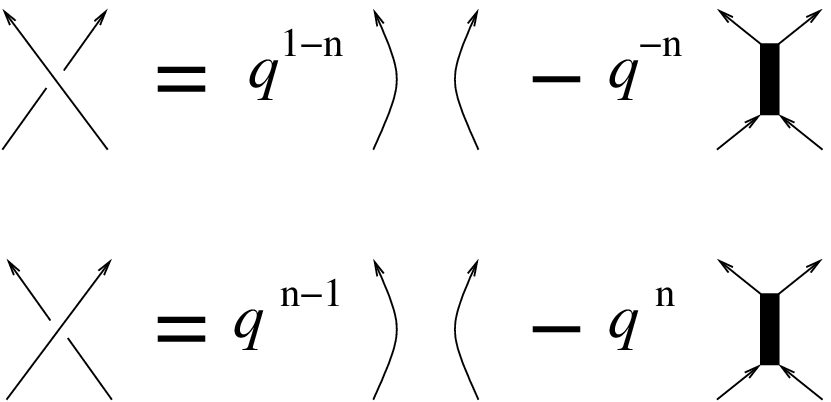}\caption{Reducing to 
 planar graphs} \label{diff}
 \end{figure}

The goal of the present paper is to construct, for each $n>0,$
a doubly-graded homology theory $H_n(L)$ with the Euler 
characteristic $P_n(L).$
The polynomial $P_n(L)$ can be computed by
breaking up each crossing into a linear combination of 
diagrams of flat trivalent graphs, as in figure~\ref{diff}.  
 
In these planar graphs some edges are oriented so that the 
neighbourhood of each unoriented edge (depicted by a thick 
line, and referred to from now on as  a "wide edge") looks 
as on the rightmost pictures in figure~\ref{diff}. Two 
oriented edges "enter" the wide edge at one vertex and two 
oriented edges "leave" it the other. Of course, this arrangement 
could be used to provide each wide edge with a canonical 
orientation, but we won't need it. In addition, oriented 
loops are allowed (an oriented loop is a crossingless plane 
projection of the oriented unknot). 

\begin{figure} [htb] \drawing{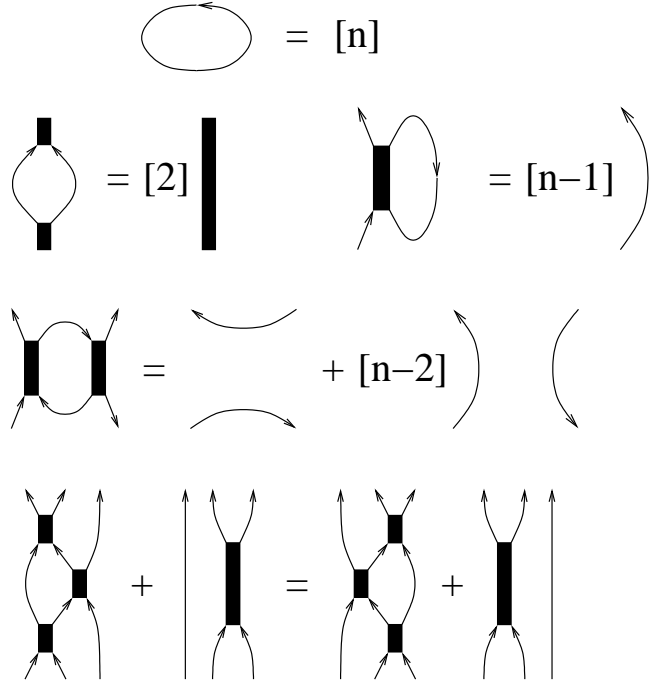}\caption{Graph skein 
relations, $[i]=\frac{\textstyle q^i -q^{-i}}{\textstyle q-q^{-1}}$} 
\label{graphrel}
 \end{figure}

There is a unique way to assign a Laurent polynomial 
$P_n(\Gamma)\in\Z[q,q^{-1}]$ to each such graph 
$\Gamma$ so as to satisfy all skein relations in 
figure~\ref{graphrel}. 

In the representation theory language, an oriented edge 
stands for the vector representation $V$ of quantum 
$\mf{sl}(n),$ and a wide edge for its (quantum) exterior 
power $\Lambda^2 V.$ The trivalent vertex is the unique 
(up to scaling) intertwiner between $V^{\otimes 2}$ 
and $\Lambda^2 V.$ The polynomial $P_n(\Gamma)$ has 
\emph{non-negative} coefficients, see [MOY]. 

\begin{figure} [htb] \drawing{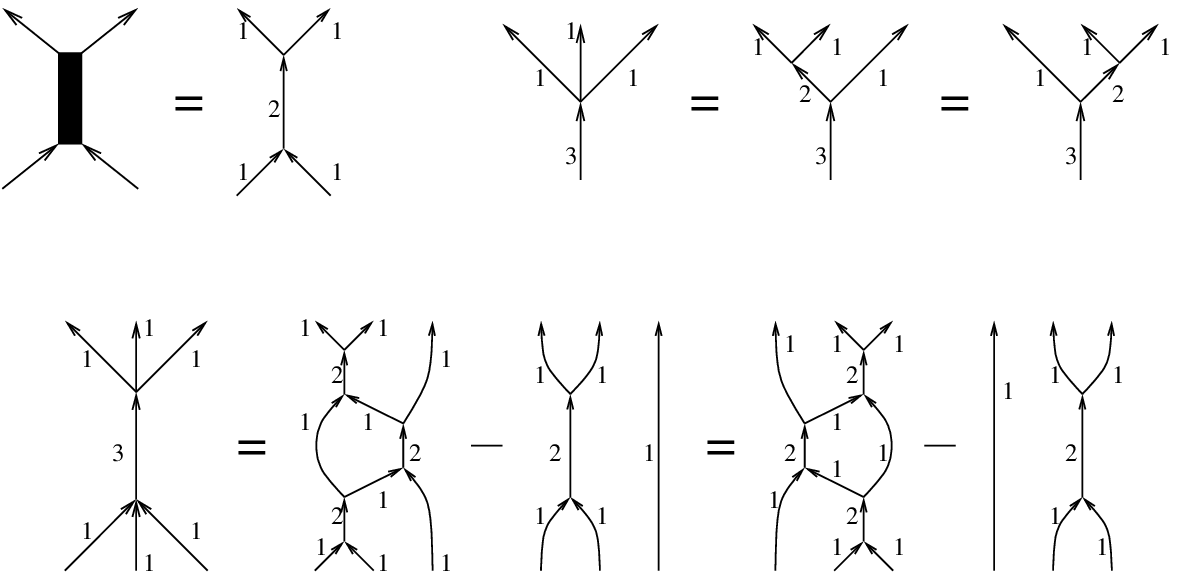}
\caption{Murakami-Ohtsuki-Yamada terminology and the appearance 
of edges labelled by $3.$}  \label{triple}
 \end{figure}

This calculus of planar graphs and its generalization 
to arbitrary exterior powers of $V$ was developed by 
Murakami, Ohtsuki and Yamada [MOY]. Oriented edges of 
graphs in their calculus carry labels from $1$ to $n-1$ 
that denote fundamental weights of $sl(n).$  
The part of their calculus that we use deals only with 
edges labelled by $1$ and $2.$ In our notation we omit these 
labels; instead, we indicate edges labelled by $2$ by wide 
lines, see figure~\ref{triple} top left. Consistency of 
figure~\ref{graphrel} relations is shown in [MOY].  
 
The last relation in figure~\ref{graphrel} can be rewritten to 
introduce edges labelled by $3$ (see figure~\ref{triple}), and 
to say that the difference of certain two endomorphisms of 
$V^{\otimes 3}$ is a multiple of the projection onto the 
irreducible summand $\Lambda^3 V.$  

The polynomial invariant $P_n(L)$ of an oriented link $L$
can be computed by choosing a plane diagram $D$ of $L,$ 
resolving each crossing as shown in figure~\ref{diff} and 
summing $P_n(\Gamma)$ weighted by powers of $q$ over all resolutions 
 $\Gamma$ of $D:$ 
 \begin{equation*}
  P_n(L) = P_n(D) \define \sum_{\mathrm{resolutions }
 \hspace{0.05in}\Gamma} q^{\alpha(\Gamma)} P_n(\Gamma),
 \end{equation*} 
with $\alpha(\Gamma)$ determined by figure~\ref{diff} rules. 
Independence from the choice of $D$ follows from the equations  
in figure~\ref{graphrel}. They imply that 
$P_n(D_1)=P_n(D_2)$ whenever $D_1$ and $D_2$ are 
related by a Reidemeister move. 

To construct homology theory $H_n,$ we first categorify 
 $P_n(\Gamma),$ by defining in a rather roundabout way  
a graded $\Q$-vector space 
$H(\Gamma)= \oplusop{j\in \Z} H^j(\Gamma)$ such that
 \begin{equation*} 
  P_n(\Gamma)=\sum_{j\in \Z}\mathrm{dim}_{\Q}H^j(\Gamma)
   \hspace{0.05in} q^j.
 \end{equation*}

With $\Gamma$ we associate a 2-periodic complex $C(\Gamma)$ 
of graded $\Q$-vector spaces
 $$ C^0(\Gamma) \lra C^1(\Gamma) \lra C^0(\Gamma) $$ 
and define $H(\Gamma)$ as the degree $i$ cohomology of this 
complex, where $i$ is the parity of the number of components of 
link $L.$ 
The construction of $C(\Gamma)$ is based on the notion of a matrix 
factorization. An $(R,w)$-factorization $M$ over a 
commutative ring $R$ (where $w\in R$) consists of two free 
$R$-modules and two $R$-module maps 
 \begin{equation*} 
  M^0 \stackrel{d}{\lra} M^1 \stackrel{d}{\lra} M^0
 \end{equation*} 
such that $d^2(m)=wm$ for any $m\in M.$ The most commonly considered 
case in the literature is when $R$ is the ring of power series, 
and $w$ satisfies a certain 
non-degeneracy assumption (that the quotient $R/(w)$ is an 
isolated singularity). Such $M$ is called a \emph{matrix factorization.} 
When $w$ is homogeneous, one can switch from power series to 
polynomials. 

Matrix factorizations appeared in commutative algebra in 
early and mid-eighties [E1], [B], [Kn], [S], [BEH] in the 
study of isolated hypersurface singularities, and 
much more recently in string theory, as boundary conditions 
for strings in Landau-Ginzburg models [KL1-3]. 

The tensor product 
$M\otimes_R N$ of an $(R,w)$-factorization $M$ and an 
$(R,-w)$-factorization $N$ is a 2-periodic complex of $R$-modules, 
with well-defined cohomology. More generally, given a finite set 
 $\{ w_1, \dots, w_m\}$ of elements of $R$ that sum to zero, and 
an $(R,w_i)$-factorization $M_i,$ for each $i,$ the tensor product 
 $$M_1\otimes_R M_2 \otimes \dots \otimes_R M_m$$ 
is a 2-periodic complex and its cohomology $H(\otimesop{i} M_i)$ 
is a $\Z_2$-graded $R$-module. 

\begin{figure} [htb] \drawing{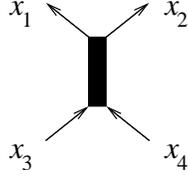}\caption{Near a 
wide edge} \label{fg-xonefour}
 \end{figure}

Starting with a resolution $\Gamma$ of a link diagram, we denote by 
$E$ the set of oriented edges of $\Gamma$ and by $R$ the ring of 
polynomials in $x_j, j\in E.$ We give each $x_j$ degree $2,$ making 
$R$ graded. Assume for simplicity that $\Gamma$ has no oriented 
loops and each wide edge of $\Gamma$ borders exactly four oriented 
edges (no oriented edge shares both endpoints with the same wide 
edge). Let $T$ be the set of wide edges. Choosing a  $t\in T,$ 
denote the oriented edges at $t$ by $1,2,3,4$ (we think of $1,2,3,4$ 
as elements of $E$) so that the four corresponding variables 
are $x_1, x_2, x_3, x_4,$ see figure~\ref{fg-xonefour}.  

Assign the polynomial  
    $$w_t = x_1^{n+1}+x_2^{n+1}-x_3^{n+1}-x_4^{n+1} $$ 
to the edge $t.$ This polynomial lies in the ideal generated 
by $x_1+ x_2 - x_3 - x_4$ and $x_1x_2 - x_3 x_4$ (since 
$x^{n+1} + y^{n+1}$ is a polynomial $g(x+y,xy)$ in $x+y$ and 
 $xy$). Therefore, we can write 
    $$ w_t = (x_1 + x_2 - x_3 - x_4) u_1 + 
    (x_1 x_2 - x_3 x_4) u_2 $$
for some polynomials $u_1, u_2.$ The latter are not uniquely 
defined, but the indeterminacy is easy to describe. We choose 
 \begin{eqnarray*} 
 u_1 & = & \frac{x_1^{n+1}+x_2^{n+1}-g(x_3+x_4,x_1x_2)}{x_1+
 x_2-x_3-x_4} , \\
 u_2 & = & \frac{g(x_3+x_4,x_1x_2)- x_3^{n+1}-x_4^{n+1}}{x_1x_2-
 x_3 x_4} .
 \end{eqnarray*}  
Let $C_t$ be the tensor product of factorizations 
 \begin{equation*} 
 R\xrightarrow{x_1+x_2-x_3-x_4} R \stackrel{u_1}{\lra} R
 \end{equation*} 
 and 
 \begin{equation*} 
 R \xrightarrow{x_1x_2-x_3x_4} R \stackrel{u_2}{\lra} R.  
 \end{equation*}
 $C_t$ is an $(R,w_t)$-factorization. Define $C(\Gamma)$ 
to be the tensor product of $C_t,$ over all wide edges $t,$
 \begin{equation*} 
   C(\Gamma) \define \otimesop{t\in T} C_t. 
 \end{equation*}
The square of the differential in $C(\Gamma)$ is 
the sum of $w_t,$ over all wide edges $t,$
 $$ d^2 = \sum_t w_t = 0.$$
The sum is $0,$ since for each oriented edge $i$ the term 
$x_i^{n+1}$ appears twice in the sum, once with positive 
and once with the negative sign. We see that $C(\Gamma)$ 
is a 2-periodic complex of $\Q$-vector spaces. 
 
\begin{figure} [htb] \drawing{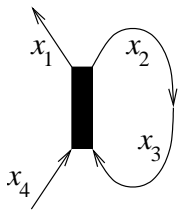}\caption{$x_2$ and $x_3$ 
 are on the same edge} \label{loop}
 \end{figure}

If, in addition, $\Gamma$ contains $k$ oriented loops, 
to define $C(\Gamma)$ we  
tensor the product of $C_t$'s with $k$ copies of the vector space 
 $\mathrm{H}^{\ast}(\mathbb{CP}^{n-1},\Q)$
and, if $k$ is odd, shift the complex. The shift $M\bracket{1}$ 
of a factorization $M^0 \lra M^1 \lra M^0$ is 
 $$ M^1 \lra M^0 \lra M^1.$$
 If, for some $t,$ some 
variables (say, $x_2$ and $x_3$) belong to the same oriented 
edge (see figure~\ref{loop}), we quotient the ring $R$ and the 
complex by the corresponding relation ($x_2=x_3$).  
  
We let $H(\Gamma)$ be the cohomology of $C(\Gamma).$ 
The latter has cohomology only in degree equal to the parity 
of the number of components of link $L.$  
Additional internal grading in $C(\Gamma)$ coming from the 
grading in the polynomial ring $R$ induces a $\Z$-grading 
on $H(\Gamma).$ This, in the nutshell, is how we define 
 $H(\Gamma),$ which serve as building blocks for complexes 
 $C(D)$ assigned to plane diagrams $D.$  

\vspace{0.1in} 

\begin{figure} [htb] \drawing{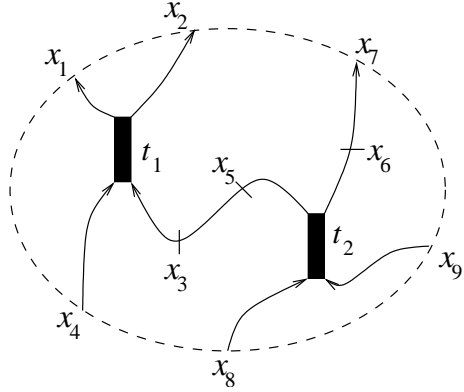}\caption{A graph 
with boundary points and marks} 
\label{diagbound}
 \end{figure}

In this paper we work in greater generality, first by allowing 
graphs that lie inside a disc and have points on its 
boundary. To such graph we assign a factorization, rather 
than a 2-periodic complex. Second, we place one or more marks 
on each oriented edge, with variables $x_i$ assigned to the 
marks, rather than just to the edges (see figure~\ref{diagbound}). 
Boundary points also count 
as marks. Let $E$ be the set of marks and $R$ the ring of 
polynomials in variables $x_i,$ for $i\in E.$  

To an arc bounded by marks $x_i,x_j$ and oriented from $x_j$ 
to $x_i$ we assign the factorization 
$L^i_j:$ 
\begin{equation*} 
  R \xrightarrow{\pi_{ij}} R \xrightarrow{x_i - x_j} R,
\end{equation*} 
where
\begin{equation*} 
  \pi_{ij}= \frac{x_i^{n+1} - x_j^{n+1}}{x_i - x_j}.
\end{equation*}  
We then form the tensor product $C(\Gamma)$ of 
factorizations $C_t,$ over all wide edges $t,$ and 
factorizations $L^i_j,$ over all arcs. Now we ignore variables 
$x_i$ assigned to internal marks, and view $C(\Gamma)$ 
as an $(R',w)$-factorization, with $R'$ being the ring 
of polynomials in $x_i,$ over all boundary points $i,$ and 
$$w = \sum_i \pm x_i^{n+1},$$ 
with signs determined by 
the orientation of $\Gamma$ near boundary points. 
We prove that, in the homotopy category of factorizations, 
$C(\Gamma)$ does not depend on how we place internal marks. 

For example, if $\Gamma$ is given in figure~\ref{diagbound}, 
to the wide edges $t_1$ and $t_2$  we assign polynomials  
$w_{t_1}$ and $w_{t_2}$ in variables $x_1,x_2,x_3,x_4$ and 
$x_5,x_6,x_8,x_9,$ respectively, and factorizations 
$C_{t_1}, C_{t_2}.$  
Form the tensor product factorization 
 $$C(\Gamma)=C_{t_1}\otimes_R C_{t_2} \otimes_R L^3_5 
 \otimes_R L^7_6,$$ 
where $R$ is the ring of polynomials in $x_1, \dots, x_9.$ 
We then ignore the variables $x_3,x_5,x_6$ that came from 
internal marks, and treat  $C(\Gamma)$ as an 
$(R',w)$-factorization with 
$$w=x_1^{n+1}+ x_2^{n+1}-x_4^{n+1} + x_7^{n+1}-x_8^{n+1}-x_9^{n+1},$$ 
and $R'$ being the ring of polynomials in "boundary" variables 
$x_1,x_2,$  $x_4,x_7,$  $x_8,x_9.$ As an $R'$-module, $C(\Gamma)$ is 
free but has infinite rank. However, by stripping off contractible 
summands, $C(\Gamma)$ can be reduced 
to a finite rank factorization. 

\begin{figure} [htb] \drawing{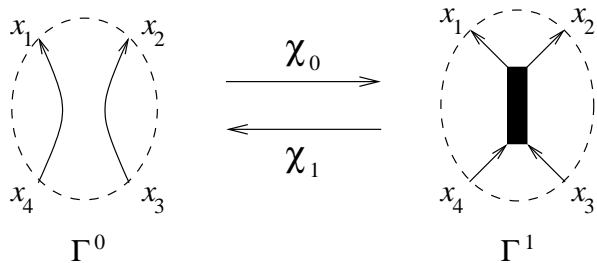}\caption{Graphs $\Gamma^0$
and $\Gamma^1$} 
\label{pair2}
 \end{figure}

We prove that factorizations $C(\Gamma)$ have direct sum 
decompositions that mimic skein relations in 
figure~\ref{graphrel}, and define homomorphisms $\chi_0$ and 
$\chi_1$ between factorizations $C(\Gamma^0)$ and $C(\Gamma^1)$
in figure~\ref{pair2}.  

\begin{figure} [htb] \drawing{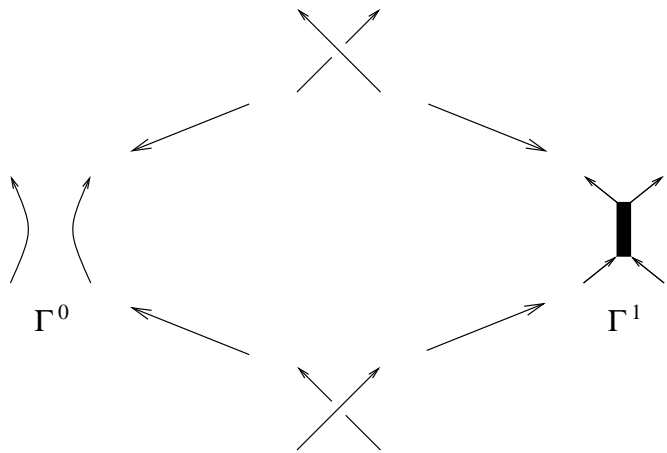}\caption{Resolving 
  crossings} \label{flats}
 \end{figure}

Then we consider oriented tangles in a 3-ball $B^3,$ such that all boundary 
points of a tangle lie on a fixed great circle of the boundary 
sphere. Let $D$ be a generic projection of tangle $L$ onto the 
plane of this great circle. We separate crossings of $D$ 
into positive and negative, following the rule in figure~\ref{posneg}. 
 To each crossing we assign two planar graphs (resolutions of this crossing), 
see figure~\ref{flats}. 

To a diagram of a single crossing we assign the 
complex of factorizations 
$$ 0 \lra C(\Gamma^0) \stackrel{\chi_0}{\lra} C(\Gamma^1)\lra 0 $$
if the crossing is positive, and the complex 
$$0 \lra  C(\Gamma^1) \stackrel{\chi_1}{\lra} C(\Gamma^0)\lra 0 $$
if the crossing is negative (shifts in the internal grading should be added to match 
powers of $q$ that appear in figure~\ref{diff} formulas). 
In both cases we place $C(\Gamma^0)$ in cohomological degree $0.$  

In general, we place marks on each internal edge of $D,$ 
form a commutative cube of factorizations $C(\Gamma),$ over 
all resolutions $\Gamma$ of $D,$ and take the total complex $C(D)$  
of the cube. Each factorization $C(\Gamma)$ 
has additional $\Z$-grading, induced by the grading of the polynomial 
ring $R,$ and the differential is grading-preserving. If we ignore 
the differentials, $C(\Gamma)$ is a $\Z\oplus \Z\oplus \Z_2$-graded 
$R$-module. 
 
We prove that if $D_1$ and $D_2$ are related by a Reidemeister 
move, $C(D_1)$ and $C(D_2)$ are isomorphic as objects in a suitable 
homotopy category. Therefore, the isomorphism class of $C(D)$ 
in this category is an invariant of tangle $L.$

A cobordism $S$ between tangles $T_0$ and $T_1$ is an oriented 
surface properly embedded in $B^3\times [0,1]$ with $S\cap (B^3\times \{i\}) = T_i,$ 
and some additional standard assumptions.  $S$ can be presented by a "movie"--a  sequence 
of plane diagrams of its cross-sections with $B^3\times k$ for various $k\in [0,1].$ 
The sequence starts with a diagram $D_0$ of $T_0$ and ends with a diagram $D_1$ 
of $T_1.$ To such a sequence we associate a homomorphism between complexes of 
factorizations $C(D_0)$ and $C(D_1)$ (with certain grading shifts thrown in), and 
show that, up to null-homotopies and  multiplications by nonzero rationals, the 
homomorphism does not depend on the movie presentation of $S$ (as long as 
 $D_0$ and  $D_1$ are the first and the last slides).   In this way, for each positive $n,$ we obtain 
an invariant of tangle cobordisms (trivial invariant if $n=1,$ and a variation of  
the one in [Kh4] if $n=2$). 

\vspace{0.1in} 

If $L$ is a link, $C(\Gamma),$ for any resolution $\Gamma$ of $D,$ 
is a 2-periodic complex of graded $\Q$-vector spaces. 
It has cohomology groups only in degree equal to the parity of 
the number of components in $L.$ Thus, after removing 
contractible summands, $C(\Gamma)$ reduces to a graded vector 
space, and $C(D)$ to a complex of graded vector spaces. Its 
cohomology groups 
$$H_n(D)= \oplusop{i,j\in \Z} H_n^{i,j}(D)$$ 
do not depend, up to isomorphism, on the choice of the projection of 
$L.$ The invariant of a link cobordism is a homomorphism between these  
cohomology groups, well-defined up to rescaling by nonzero rational numbers. 
Due to this, $H_n(D)$ are canonically (up to rescaling) associated to $L,$ 
rather than just to its diagram $D.$ Thus, notation $H_n(L)$ is justified. 
The Euler characteristic of $H_n(L)$ is the polynomial 
$P_n(L)$: 
$$ P_n(L) = \sum_{i,j\in \Z} (-1)^i q^j \dim_{\Q} H_n^{i,j}(L).$$ 

\vspace{0.1in} 

$H(\Gamma)$ is the cohomology of the complex $C(\Gamma)$ which 
has countable dimension as a $\Q$-vector space. However, 
since $H(\Gamma)$ is finite-dimensional, it's nontrivial 
in finitely many degrees only, and an obvious upper bound 
on $|j|$  with $H^j(\Gamma)\not= 0$ is $ne$ where $e$ 
is the number of edges in $\Gamma$ (both oriented and unoriented). 
The complex $C(D)$ is built out of 
finitely many $H(\Gamma),$ and the differential, being the 
signed sum of maps $\chi_0,\chi_1,$ can be 
determined combinatorially as well. Therefore, there exists 
a combinatorial algorithm to find homology groups $H_n(L),$ 
given $L.$  

\vspace{0.06in} 

If $n=1,$ the factorization $C_t$ is contractible, and 
any graph with a wide edge has trivial homology. The 
homology of a circle is $\Q.$ This implies that 
$H_1(L)\cong \Q,$ for any link $L,$ with $\Q$ in bidegree $(0,0).$  

\vspace{0.06in} 

When $n=2$ and $\Gamma$ is closed (has no boundary points), 
its cohomology groups $H(\Gamma)$ are isomorphic to $A_2^{\otimes k},$ 
where $A_2$ is the cohomology ring of the 2-sphere and $k$ the 
number of circles in the diagram given by deleting all wide edges 
of $\Gamma.$ The equivalence of $H_2$ with a version 
of homology theory in [Kh1] follows easily from this observation. 
Specifically, for each link $L$ there is an isomorphism 
 $$ H_2^{i,j}(L) \cong \mc{H}^{i,-j}(L^!)\otimes_{\Z}\Q,$$  
where $L^!$ is the mirror image of $L,$ and  
$\mc{H}$ is the homology theory defined in [Kh1, Section 7] 
(programs computing $\mc{H}(L)$ were implemented by Dror Bar-Natan [BN1] 
and Alexander Shumakovitch [Sh]).

\vspace{0.06in} 

We conjecture that $H_3$ is isomorphic to 
the homology theory constructed in [Kh3], after the latter 
is tensored with $\Q.$ Both theories are doubly-graded, have 
the same Euler characteristic, and are built out of similar 
long exact sequences. The language of foams used in [Kh3] and [R] 
should extend from $n=3$ to all $n$ and lead to a better 
understanding of $H_n.$   

\vspace{0.07in} 

Our categorification of polynomials $P_n(L)$ comes with several 
obvious caveats: 
\begin{itemize} 
 \item Homology groups $H_n(L)$ are finite-dimensional $\Q$-vector 
spaces rather than finitely-generated abelian groups (as in [Kh1], 
[Kh3] for $n=2,3.$)  
 \item Defining $H_n(L)$ requires a choice of a plane diagram. A 
more intrinsic definition would be most welcome.
 \item Our invariants of link and tangle cobordisms are 
projective (defined up to overall multiplication by a non-zero 
 rational number). 
 \item Reconstructing the Ozsv\'ath-Szab\'o-Rasmussen theory (the $n=0$ 
case) using this approach would require additional ideas, lacking 
at the moment.
\end{itemize}

 \vspace{0.1in} 

Construction of $C(\Gamma)$ involves tensoring factorizations 
$$R \stackrel{a}{\lra} R \stackrel{b}{\lra} R$$ 
for various $a$'s and $b$'s. We do an elementary study of these tensor 
products in Section~\ref{sec-koszul}. Section~\ref{sec-potent} contains 
a review of matrix factorizations and their properties. In 
Section~\ref{sec-functors} we explain how to view factorizations as 
functors and present the identity functor via a factorization. 
These two sections also introduce general framework for diagrammatical 
interpretation of matrix factorizations. Section~\ref{sec-homog} 
treats graded factorizations and explains how to modify 
the material of previous sections to cover this case. 
Section~\ref{sec-planar} is the computational core of the paper. We 
define factorization $C(\Gamma)$ assigned to a planar graph $\Gamma,$ 
construct morphisms $\chi_0,\chi_1$ between graphs $\Gamma^0,\Gamma^1$
in figure~\ref{pair2}, and prove direct sum decompositions of 
factorizations $C(\Gamma)$ that lift skein relations in 
figure~\ref{graphrel}. In Section~\ref{sec-tangle} we associate 
a complex of factorizations $C(D)$ to a plane giagram $D$ of a tangle and state 
Theorem~\ref{main-theorem}. This theorem, claiming 
the invariance of $C(D)$ under Reidemeister moves in a suitable 
category of complexes of factorizations, is proved in 
Section~\ref{sec-invar}. In Section~\ref{sec-corners} we use matrix 
factorizations to construct 2-dimensional topological quantum field 
theories with corners. The next section deals with functoriality 
of our tangle invariant, extending it to tangle cobordisms. In the last 
section we briefly outline an approach to categorification of quantum 
invariants of links colored by arbitrary fundamental representations of $sl(n).$ 

\vspace{0.1in} 

{\bf Acknowledgments:} M.K. is indebted to Ragnar-Olaf 
Buchweitz, Igor Burban, and Paul Seidel for introducing
him to matrix factorizations, and to Igor Burban for sharing 
 M.K.'s fascination with 2-periodic triangulated 
categories. 
 L.R. would like to thank Anton Kapustin for explanations 
regarding topological 2-dimensional models and his work with 
 Yi Li. We are grateful to Dror Bar-Natan and Greg Kuperberg for 
interesting discussions. While writing this paper we were 
partially supported by NSF grants DMS-0104139 and  DMS-0196131.


\section{The cyclic Koszul complex} \label{sec-koszul}

Let $R$ be a Noetherian commutative ring. The Koszul complex 
$R(a_1,\dots, a_m)$ associated to a sequence $a_1, \dots, a_m,$
where $a_i\in R,$ is the tensor product (over $R$) of complexes 
 \begin{equation*} 
  0 \lra R \stackrel{a_i}{\lra} R \lra 0 
 \end{equation*} 
for $i=1,\dots, m.$ A sequence $(a_1, \dots, a_m)$ is called 
$R$-regular iff $a_i$ is not a zero divisor in the quotient 
ring $R/(a_1, \dots, a_{i-1})R$ for each $i=1, \dots, m,$
and $R/(a_1, \dots, a_m)R\not= 0.$   
If the sequence is $R$-regular, its Koszul complex has 
cohomology in the rightmost degree only. For the converse 
to be true it suffices for $R$ to be local, see [E1, Section 17].  

The Koszul complex can be written in a more intrinsic 
way. Let $N$ be a free $R$-module of rank $m$ and 
$a:R\to N$ an $R$-module homomorphism. The Koszul complex 
of $a$ is the complex of exterior powers of $N$: 
 \begin{equation*} 
 0 \lra R \lra N \lra \Lambda^2 N \lra \Lambda^3 N \lra \dots,  
 \end{equation*} 
with the differential being the exterior product with $a(1)\in N.$ 
Choosing a basis of $N$ and writing $a(1)$ in the coordinates 
as $(a_1, \dots, a_m),$ we recover the Koszul 
complex of this sequence. 

\begin{definition} Let $w\in R.$ An $(R,w)$-duplex consists of 
two $R$-modules $M^0,M^1$ and module maps 
$$ M^0 \stackrel{d}{\lra} M^1 \stackrel{d}{\lra} M^0$$  
such that $d^2=w$ (as endomorphisms of $M^0$ and $M^1$). 
\end{definition}

In other words, $d^2m = wm$ for any $m\in M^0 \oplus M^1.$ 
"Duplex" is a term from [FKS], where it was used in a more 
general situation, with $R$ a (possibly noncommutative) ring and 
$w$ its central element. 

A homomorphism $f:M\lra N$ of duplexes is a pair of 
homomorphisms $f^0:M^0 \lra N^0$ and $f^1: M^1 \lra N^1$ that make 
the diagram below commute. 
 \begin{equation*}
   \begin{CD}
     M^0 @>{d^0}>> M^1 @>{d^1}>> M^0   \\
     @V{f^0}VV      @V{f^1}VV   @V{f^0}VV \\
     N^0 @>{d^0}>> N^1 @>{d^1}>> N^0
   \end{CD}
 \end{equation*}
The category of $(R,w)$-duplexes and duplex homomorphisms is 
abelian. 

  A homotopy $h$ between homomorphisms $f,g:M\lra N$ of duplexes 
is a pair of maps $h^i:M^i \lra N^{i-1}$ (with indices understood 
modulo $2$) such that 
 $$f-g = h d_M + d_N h.$$
Null-homotopic morphisms constitute an ideal in the category 
of duplexes. We call the quotient category by this ideal 
the category of duplexes up to homotopies, or simply the 
\emph{homotopy category of duplexes}. This category is triangulated. We denote 
the shift functor by $\langle 1 \rangle.$ 

An $(R,0)$-duplex $M$ is equivalent to a 2-complex (periodic 
complex with period two) of $R$-modules. We denote the cohomology 
of an $(R,0)$-duplex $M$ by $H(M)\cong H^0(M)\oplus H^1(M).$

\begin{definition} A factorization (or matrix factorization) is 
an $(R,w)$-duplex such that $M^0,M^1$ are free $R$-modules. 
\end{definition} 
The homotopy category of $(R,w)$-factorizations is triangulated. 

To a pair of elements $a,b\in R$ we associate an $(R,ab)$-factorization 
 $$R \stackrel{a}{\lra}R \stackrel{b}{\lra} R,$$ 
denoted $\{ a,b\}.$  

If $\bba=(a_1, \dots ,a_m)$ and $\bbb=(b_1, \dots ,b_m)$ are two 
sequences of elements of $R,$ we consider the tensor 
product factorization $\{ \bba, \bbb\} \define \otimes_i \{ a_i, b_i\}
,$ where the tensor product is over $R.$ We call $\{ \bba,\bbb\}$ 
\emph{the Koszul factorization} of the pair $(\bba,\bbb).$ We say that a pair 
$(\bba,\bbb)$ is \emph{orthogonal} if 
  $$\bba\bbb\define \sum_i a_i b_i=0.$$ 
The tensor product $\{ \bba, \bbb\}$ is a 2-complex iff 
$(\bba,\bbb)$ is orthogonal, in which case we call it the periodic 
(or cyclic) Koszul complex of $(\bba,\bbb).$  

\emph{Remark:} We allow the case $m=0.$ Then $\bba$ and $\bbb$ are empty 
sequences, $w=0,$  and the Koszul factorization is $R\lra 0 \lra R.$ 

\emph{Remark:} Factorizations $\{ \bba,\bbb\}$ were defined 
in [BGS] and used there in the classification of finite CM-type 
singularities. 

\begin{prop} If the entries of $\bba$ and $\bbb$ generate $R$ as an 
$R$-module, the factorization $\{ \bba,\bbb\}$ is contractible 
(its identity endomorphism is null-homotopic). 
\end{prop} 

\emph{Proof:} Write $1= \sum_i(a'_i a_i + b'_i b_i).$ The pair 
$(\bba',\bbb'),$ where $\bba'=(a'_1, \dots, a'_m)$ and $\bbb=
 (b'_1,\dots,b'_m)$ 
defines a homotopy between the identity and the 
zero endomorphisms of $\{\bba,\bbb\}.$  $\square$ 

\begin{corollary} If $(\bba,\bbb)$ is orthogonal and the entries of $\bba$ 
and $\bbb$ generate $R$ as an $R$-module, the 2-complex $\{ \bba, 
  \bbb\}$ is acyclic. 
\end{corollary} 

\vspace{0.1in} 

Let $I_{\bba,\bbb}\subset R$ be the ideal of $R$ generated by 
the entries of $\bba$ and $\bbb.$ 

\begin{prop} $H(\{ \bba,\bbb\})$ is an $R/I_{\bba,\bbb}$-module. 
\end{prop} 

\emph{Proof:} Multiplications by  $a_i$ and $b_i$ 
are null-homotopic endomorphisms of $\{a_i,b_i\},$ and, therefore, 
of $\{\bba, \bbb\}.$ 

$\square$ 

\vspace{0.1in}  

Let $\bba^i,\bbb^i$ be the sequences obtained from $\bba$ and $\bbb$ 
by omitting $a_i$ and $b_i,$ respectively. The factorization 
$\{ \bba,\bbb\} $ is the "total factorization" of the bifactorization
  \begin{equation*} 
  \{ \bba^i, \bbb^i \} \stackrel{a_i}{\lra} \{ \bba^i,\bbb^i\} 
  \stackrel{b_i}{\lra} \{ \bba^i, \bbb^i\}. 
  \end{equation*} 
If $b_i$ is a nonzerodivisor in $R,$ the second map is injective 
 (as map between $R$-modules), and
 the direct sum of the middle $\{\bba^i,\bbb^i\}$ and its image 
 under the differential is contractible. Furthermore, 
 suppose that $(\bba,\bbb)$ is orthogonal. Then the quotient of 
 $\{\bba,\bbb\}$ by this contractible subcomplex is isomorphic to 
 $\{\bba^i,\bbb^i\}_{R'},$ the periodic Koszul complex of the pair 
$(\bba^i,\bbb^i)$ in the quotient ring $R'\define R/(b_i).$ The 
quotient map induces an isomorphism on cohomology, 
  \begin{equation*} 
    H(\{ \bba,\bbb\}) \cong H(\{ \bba^i, \bbb^i\}_{R'}),
  \end{equation*} 
the isomorphism is that of $R'$-modules. 

\vspace{0.1in} 

\begin{corollary} \label{cor-reg} 
If $(\bba,\bbb)$ is orthogonal and $\bbb$ is $R$-regular then  
 \begin{equation*} 
  H^0(\{\bba,\bbb\}) \cong R/(b_1, \dots, b_m) \hspace{0.1in} 
  \mathrm{and} \hspace{0.1in} H^1(\{\bba,\bbb\})=0. 
 \end{equation*} 
\end{corollary}  

\vspace{0.1in} 

Likewise, if $a_i$ is a nonzerodivisor in $R$ and $(\bba,\bbb)$ is 
orthogonal, the quotient map 
 \begin{equation*} 
  \{\bba, \bbb\}\langle 1\rangle \lra \{\bba^i,\bbb^i\}_{R/(a_i)}
 \end{equation*}  
induces an isomorphism on cohomology. 

\vspace{0.1in} 

Motivated in part by Corollary~\ref{cor-reg}, we introduce 

\begin{definition} An orthogonal pair $(\bba,\bbb)$ is called 
homologically $R$-regular if
 $H^0(\{\bba,\bbb\})\not= 0$ and $H^1(\{ \bba,\bbb\})=0.$
\end{definition} 

If $(\bba,\bbb)$ is homologically $R$-regular then $I_{\bba,\bbb}$ 
is a proper ideal in $R.$ 

\vspace{0.1in} 

\emph{Example:} Let $m=1.$ An orthogonal pair $(a_1,b_1)$ is 
homologically $R$-regular iff $a_1R=\mathrm{Ann}(b_1)$ and 
$b_1R$ is  a proper subset of $\mathrm{Ann}(a_1).$ 
Furthermore, when $R$ is a domain, $(a_1,b_1)$ is homologically 
$R$-regular iff $a_1=0$ and $b_1$ is not invertible. 

\vspace{0.1in} 

\emph{Example:} Let ${\bf 0}=(0,\dots, 0).$ The pair $({\bf 0},\bbb)$ 
is homologically $R$-regular iff the Koszul complex of $\bbb$ has 
cohomology in the rightmost degree only. 

\vspace{0.1in} 

Suppose there is a subset $J\subset \{ 1, \dots, m\}$ with 
even number of elements such that the sequence 
 $(h_1, \dots, h_m)$ is $R$-regular, where $h_i=a_i$ if $i\in J$ 
 and $h_i=b_i$ otherwise. Then $(\bba,\bbb)$ is homologically 
 $R$-regular and $H^0(\{\bba, \bbb\})\cong A/(h_1, \dots, h_m).$ 

\vspace{0.1in} 

We can likewise define the notion of homologically 
$M$-regular pair $(\bba,\bbb)$ for any $R$-module $M,$ and even 
for any $(R,w)$-duplex $M$. For the latter, 
we require $\bba\bbb=-w, H^0(M\otimes_R \{\bba,\bbb\})\not= 0$ and 
$H^1(M\otimes_R \{ \bba,\bbb\})=0.$ 

\vspace{0.2in} 

Let $N$ be a finitely-generated free $R$-module and 
$\alpha: R \to N, \beta: N \to R$ be $R$-module maps. Consider 
the factorization $\{ \alpha, \beta\}$ given by   
 \begin{equation} 
   \Lambda^{even} N \xrightarrow{\wedge \alpha+ \neg \beta} 
  \Lambda^{odd} N \xrightarrow{\wedge \alpha + \neg \beta} 
  \Lambda^{even} N
 \end{equation} 
where $\wedge \alpha$ is the wedge product with $\alpha(1),$ 
 $$ \wedge \alpha : \Lambda^i N \lra \Lambda^{i+1} N,$$ 
 $\neg \beta$ is the contraction with $\beta,$ 
 $$ \neg \beta: \Lambda^i N \lra \Lambda^{i-1} N,$$ 
and 
 $$ \Lambda^{even}N = \oplusop{i} \Lambda^{2i}N, 
  \hspace{0.2in} 
    \Lambda^{odd}N = \oplusop{i} \Lambda^{2i+1}N. $$ 

 $\{ \alpha, \beta\}$ is a 2-complex iff $\beta\alpha=0,$ 
in which case we say that $\alpha$ and $\beta$ are orthogonal. 
Choosing a basis of $N$, we can write 
$\alpha=\bba=(a_1, \dots, a_m)^T$ and 
 $\beta=\bbb=(b_1, \dots, b_m).$ Then 
 $\{ \alpha,\beta\}$ is isomorphic to $\{\bba, \bbb\}.$  

\vspace{0.1in} 

Thus, with each element $\gamma\in N\oplus N^{\ast}$ we 
associate a factorization 
$$\{ \gamma\} \define \{ \alpha, \beta\}$$ 
where $\gamma= \alpha + \beta$ and $\alpha\in N, \beta\in N^{\ast}.$ 

Suppose $g:N\to N$ is an automorphism of the $R$-module 
$N.$ Write the composition $\beta\alpha$ as 
$\beta \gamma^{-1}\gamma \alpha.$ The factorizations 
$\{ \alpha,\beta\}$ and $\{ \gamma\alpha, \beta\gamma^{-1}\}$ are 
isomorphic, and $GL(N,R)$ acts on the set of pairs $(\alpha,\beta)$ 
preserving isomorphism classes of factorizations $\{ \alpha, \beta\}.$  

We let the group $H=(\Z_2)^{\times m}$ act on pairs 
by permuting $a_i$ with $a_i$ in  
 $(a_1,a_1), \dots, (a_m,a_m),$ for $1\le i \le m.$  

For $\sigma\in H,$ factorizations $\{ \bba,\bbb\}$ and 
$\{ \sigma(\bba,\bbb)\}$ are isomorphic if $\sigma$ is an even 
permutation in $H\subset \mathbb{S}_{2m}.$ If $\sigma$ is 
odd, factorization $\{ \sigma(\bba,\bbb)\}$ is isomorphic to 
 $\{ \bba, \bbb\}\langle 1\rangle.$ 

The natural $R$-linear pairing between $N^{\ast}$ and $N$ induces 
a symmetric bilinear form on $N\oplus N^{\ast}.$ 
Let $G'\subset GL(N\oplus N^{\ast},R)$ be the subgroup 
 generated by $GL(N,R)$ and $H$ 
that act as described above. The action of $G'$ preserves 
the inner product on $N\oplus N^{\ast}$ (with values in $R$). 
Let $G$ be the subgroup of $G'$ that consists of all products 
$g_1 \sigma_1 \dots g_r \sigma_r$ 
of elements from $GL(N,R)$ and $H$ (over all $r$), 
such that  $\sigma_1\dots \sigma_r$ is an even permutation. 

The following result is clear. 

\begin{prop} $\{\gamma\}$ and $\{ g\gamma\}$ are isomorphic 
factorizations for any $g\in G$ and $\gamma\in N\oplus N^{\ast}.$  
\end{prop} 

The cyclic Koszul complex can be thought of as the square root of 
the Koszul complex. Namely, 
$$ \Hom_R(\{ \bba, \bbb\}, \{ \bba, \bbb \} ) \cong 
 \{ \bba, \bbb\}\otimes_R \{ -\bbb, \bba \} \cong R(\bba, \bbb),$$ 
where $R(\bba,\bbb)$ is the Koszul complex of the length $2m$ 
sequence given by concatenating $\bba$ and $\bbb,$ with the 
grading collapsed from $\Z$ to $\Z_2.$


\section{Potentials, isolated singularities, and matrix 
 factorizations} \label{sec-potent} 

{\bf Potentials and their Jacobian algebras.} Start with a finite set 
$x=\{ x_1, \dots, x_k\}$ of variables and let 
$$R=\Q[[x]]\define\Q[[x_1,\dots, x_k]]$$ 
be the algebra of power series in $x_1, \dots, x_k$ with rational 
coefficients. Denote by $\mf{m}$ the unique maximal ideal 
in $R$ (it's generated by $x_1, \dots, x_k$).   
 
We say that  a polynomial $w=w(x)\in \mf{m}^2$ is a \emph{potential} if 
 the algebra $R/I_w$ is finite-dimensional, where $I_w\subset R$ 
is the ideal generated by partial derivatives 
$\partial_i w\define \partial_{x_i}w.$ Any 
polynomial $w\in R$ defines an algebraic map $\C^n\to \C.$ This map 
has an isolated singularity at $0$ iff $w$ is a potential. 

A polynomial $w\in \mf{m}^2$ is a potential 
iff $\partial_1w, \dots, \partial_kw$ is a regular sequence 
in $R,$ that is, $\partial_i w$ is not a zero divisor in 
$R/(\partial_1 w,\dots, \partial_{i-1}w)R$ for each $1\le i \le k.$ 

The quotient algebra $R_w \define R/I_w$ is called the local algebra 
of the singularity and the Jacobian (or Milnor) algebra of $w.$ Its 
dimension is the Milnor number of the singularity. 

 \begin{prop} The Jacobian algebra $R_w$ is symmetric. 
 \label{prop-symmetric} \end{prop} 
A non-degenerate $\C$-linear trace $\Tr :R_w\otimes_{\Q}\C\to \C$ 
is given by the residue formula 
 \begin{equation} 
  \Tr(a) = \frac{1}{(2\pi i)^k} 
  {{\mathop{\int}\limits_{|\partial_j w|=\epsilon }}}
  \frac{a dx_1 \dots dx_k}{\partial_1 w\dots \partial_k w} 
 \end{equation} 
(this is proved in [GH, Chapter 5] and [AGV]). Therefore,
the $\C$-algebra $R_w\otimes_{\Q}\C$ is symmetric. From 
[NN,Theorem 5] we deduce that $R_w$ is a symmetric 
$\Q$-algebra.  $\square$ 

  \vspace{0.1in} 

{\bf Exterior sum of potentials.} Let $x$ and $y$ be two disjoint finite sets of 
variables. Define the exterior sum of potentials $w_1(x)$ and 
$w_2(y)$ as the potential $w_1(x)+w_2(y)$ in variables $x$ and $y.$ 
The Jacobian algebra of the exterior sum of potentials is the tensor 
product of Jacobian algebras: 
 
 \begin{equation*} 
  R_{w_1(x)+ w_2(y)} \cong R_{w_1(x)} \otimes_{\Q} R_{w_2(y)}. 
 \end{equation*} 

 \vspace{0.1in} 

{\bf Matrix factorizations.} A \emph{matrix factorization} of a potential $w$ is a 
collection of two free $ A$-modules $M^0,M^1$ and two module maps 
$d^0: M^0 \lra M^1$ and $d^1:M^1\lra M^0$ such that 
 \begin{equation*} 
   d^0 d^1 = w \cdot \Id_{M^1}, \hspace{0.2in} 
   d^1 d^0 = w \cdot \Id_{M^0}
 \end{equation*} 
We will often write a matrix factorization as 
  \begin{equation*} 
    M^0 \stackrel{d}{\lra} M^1 \stackrel{d}{\lra} M^0
  \end{equation*} 
Modules $M^0$ and $M^1$ are not required to have finite rank. Since 
$w$ is invertible in the field of fractions $\Q((x))$ of $R,$ and 
the rank of a free $ R$-module $N$ can be defined as the dimension 
of the $\Q((x))$-vector space $N\otimes_R\Q((x)),$ we see that 
$M^0$ and $M^1$ have equal ranks (when the ranks are finite). In 
general, $rk(M^0)=rk(M^1),$ understood as equality of ordinals. 
We call $rk(M^0)$ the rank of $M.$ 

 \vspace{0.1in} 

A homomorphism $f:M\lra N$ of factorizations is a pair of 
homomorphisms $f^0:M^0 \lra N^0$ and $f^1: M^1 \lra N^1$ that make 
the diagram below commute. 
 \begin{equation*}
   \begin{CD}
     M^0 @>{d^0}>> M^1 @>{d^1}>> M^0   \\
     @V{f^0}VV      @V{f^1}VV   @V{f^0}VV \\
     N^0 @>{d^0}>> N^1 @>{d^1}>> N^0
   \end{CD}
 \end{equation*}
Denote the set of homomorphisms from $M$ to $N$ by 
$\Hom_{\MF}(M,N).$ It's an $ R$-module, with 
the action $a(f^0,f^1)=(af^0,af^1),$ for $a\in  R.$ 

Let $\MF^{all}_w$ be the category whose objects are matrix 
factorizations and morphisms are homomorphisms of factorizations. 
This category is additive and $ R$-linear. Direct sum of matrix 
factorizations $M$ and $N$ is defined in the obvious way:  
 \begin{equation*} 
  (M\oplus N)^i = M^i \oplus N^i, 
 \hspace{0.2in} d^i_{M\oplus N} = d^i_M + d^i_N.
 \end{equation*} 
The shifted factorization $M\langle 1\rangle$ is given by 
 \begin{eqnarray*} 
   M\langle 1\rangle^i & =  & M^{i+1}, 
   \hspace{0.2in} i =0,1 \hspace{0.1in} 
  \mathrm{mod} \hspace{0.1in} 2, \\
  d_{M\langle 1\rangle}^i & =  &  - d_{M}^{i+1}. 
 \end{eqnarray*} 
$\langle 1\rangle $ is a functor in $\MF^{all}_w$ whose square is 
the identity, $\langle 2 \rangle \cong \mathrm{Id}.$  

 \vspace{0.05in}  

\emph{Example:} If the set of variables $x$ is empty, then 
$R=\Q, f=0,$ and matrix factorization is a pair of vector 
spaces and a pair of maps between them such that their composition 
in any order is $0.$ We call such data a \emph{2-complex} (short for 
 \emph{2-periodic complex}).  

Matrix factorizations first appeared in the work of David Eisenbud 
[E1], who related them to maximal Cohen-Macaulay modules over 
isolated hypersurface singularities (also see [Y1, Chapter 7], [S], 
and references therein). A module $N$ over a commutative Noetherian 
local ring is called \emph{maximal Cohen-Macaulay} if the depth of $N$ equals 
the Krull dimension of the ring. Given a matrix factorization $M,$ 
the $R/(f)$-module $\mathrm{Coker}(d^1)$ is maximal Cohen-Macaulay, 
and any maximal Cohen-Macaulay module over $R/(f)$ can be presented
in this way (see the above references for details). For a 
reader-friendly treatment of the background concepts leading to 
maximal Cohen-Macaulay modules we suggest the book [E2]. 

\vspace{0.1in} 

 \begin{figure} [htb] \drawing{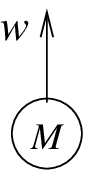}
  \caption{A matrix factorization} \label{fg-factorization}
 \end{figure}

\begin{figure} [htb] \drawing{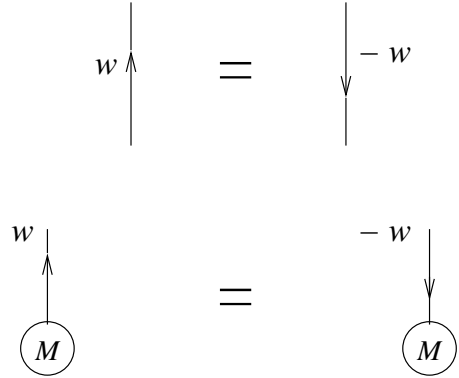}
  \caption{Orientation reversal} \label{fg-oreverse}
 \end{figure}

{\bf Graphical notation.} We denote a matrix factorization $M$ with potential 
$w$ as in figure~\ref{fg-factorization}. We allow reversing 
orientation of the arc attached to $M$ simultaneously with changing 
potential $w$ to $-w,$ see figure~\ref{fg-oreverse}. 

\begin{figure} [htb] \drawing{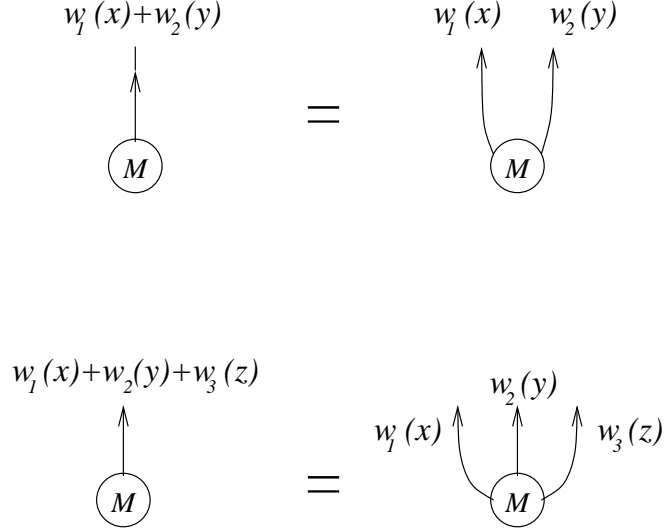}
  \caption{Diagrams for a factorization over exterior sum of 
 potentials}\label{fg-sum}
 \end{figure}

If $w$ is the exterior sum of potentials, we may introduce several 
arcs at $M,$ one for each summand, see figure~\ref{fg-sum}. 

{\bf Factorizations of finite rank.} 
Let $\MF'_w$ be the category of finite rank matrix 
factorizations. It's a full subcategory in $\MF^{all}_w.$ 
Choose bases in free $ R$-modules $M^0$ and $M^1$. The maps 
$d^0, d^1$ can be written as $m\times m$ matrices $D_0,D_1$ with 
coefficients in $ R.$ These matrices satisfy the equations  
 \begin{equation}  
   D_0D_1 = w\cdot\Id, \hspace{0.2in} D_1D_0= w\cdot\Id
 \end{equation} 
(of course, any one of these equations implies the other). 
Alternatively, we can describe this factorization by a 
$2m\times 2m$ 
matrix with off-diagonal blocks $D_0$ and $D_1$: 
  \begin{equation*}  
   D = \left( \begin{array}{cc} 0 &  D_1   \\  D_0  & 0 \end{array}
       \right),  \hspace{0.3in} D^2 = w \cdot \Id.  
  \end{equation*} 
Matrix description of objects in $\MF'_w$ extends to infinite rank 
factorizations. Matrices $D_0$ and $D_1$ then have infinite rank, 
but each of their columns has only finitely many non-zero entries.  

If factorizations $M$ and $N$ are written in matrix form, 
$M=(D_0,D_1),$ and $N=(D_0', D_1'),$ a homomorphism $f:M\to N$ is a 
pair of matrices $(F_0,F_1)$ such that $F_1 D_0 = D_0' F_0$ and 
$F_0 D_1 = D_1' F_1$ (note, though, that the two equations are 
equivalent). 

 \vspace{0.1in} 

{\bf Homotopies of factorizations.}  A homotopy $h$ between maps $f,g:M\lra N$ of 
factorizations is a pair of maps $h^i:M^i \lra N^{i-1}$ such that 
 $$f-g = h d_M + d_N h.$$
Null-homotopic morphisms constitute an ideal in the category 
$\MF^{all}_w.$ Let $\HMF^{all}_w$ be the quotient category by this 
ideal. It has the same objects as $\MF^{all}_w,$ but fewer 
morphisms: 
 \begin{equation*} 
   \Hom_{\HMF}(M,N) \define \Hom_{\MF}(M,N)/
   \{\mathrm{null}-\mathrm{homotopic}\hspace{0.1in} 
   \mathrm{morphisms} \}.  
 \end{equation*} 
Choose bases in $M^0$ and $M^1$ and write $d$ as a matrix $D.$
Differentiating the equation $D^2=w$ with respect to $x_i$ we get 
 \begin{equation*}  
   D (\partial_i D) + (\partial_i D) D = \partial_i w.
 \end{equation*} 
Therefore, the multiplication by $\partial_i w$ endomorphism 
of $M$ is homotopic to $0,$ and we obtain 

\begin{prop} \label{action-fd} The action of $ R$ on  $\Hom_{\HMF}
  (M,N)$ factors through the action of the Jacobian algebra $R_w.$
\end{prop} 

Let $\HMF'_w$ be the category with objects--finite rank matrix 
factorizations and morphisms--homomorphisms of factorizations modulo 
those homotopic to $0.$ This is a full subcategory  in $\HMF^{all}_w.$ 

 \vspace{0.05in} 

The free $R$-module $\Hom_R(M,N)$ is a 2-complex 
 \begin{equation*} 
  \Hom^0_R(M,N) \stackrel{d}{\lra} 
  \Hom^1_R(M,N) \stackrel{d}{\lra}\Hom^0(M,N), 
 \end{equation*} 
where   
 \begin{eqnarray*} 
  \Hom^0_R(M,N) & = & \Hom_R(M^0,N^0) \oplus    \Hom_R(M^1,N^1), \\
 \Hom^1_R(M,N) & = & \Hom_R(M^0,N^1) \oplus    \Hom_R(M^1,N^0), 
 \end{eqnarray*} 
and 
 \begin{equation*} 
   (df)m = d_N(f(m))+(-1)^i f (d_M(m)), \hspace{0.1in} 
   \mathrm{for} \hspace{0.1in} f\in \Hom^i_R(M,N). 
 \end{equation*} 
Denote the cohomology of this 2-complex by  
 \begin{equation*} 
   \Ext(M,N) = \Ext^0(M,N)\oplus \Ext^1(M,N). 
 \end{equation*} 
It is clear from definitions that  
  \begin{eqnarray} 
    \Ext^0(M,N) & \cong & \Hom_{HMF}(M,N), \\
    \Ext^1(M,N) & \cong & \Hom_{HMF}(M,N\langle 1\rangle). 
  \end{eqnarray} 
 
\begin{prop} \label{ext-fin} $\Ext(M,N)$ is a  finite-dimensional 
 $R_w$-module if $M$ and $N$ are finite rank factorizations. 
\end{prop} 
\emph{Proof:} We need to show that $\Hom_{\HMF}(M,N)$ is 
finite-dimensional. The latter is an $ R$-module quotient of 
$\Hom_{\MF}(M,N),$ which is a submodule of the free finite rank 
$R$-module $\Hom_R(M,N).$ Since $ R$ is Noetherian, 
subquotients of finitely-generated $ R$-modules are 
finitely-generated. The action of $ R$ on $\Hom_{\HMF}(M,N)$ 
factors through the action of the Jacobian ring $R_w,$ by 
proposition~\ref{action-fd}. Therefore, $\Hom_{\HMF}(M,N)$ is 
finite-dimensional, being finitely-generated over a 
finite-dimensional algebra. $\square$  

 \vspace{0.1in} 

\emph{Remark:} This proposition and many that follow fail if $w\in R$ 
is not a potential. However, $(R,w)$-factorizations for such 
degenerate $w$ will make important intermediate appearances, 
as our main examples come from tensor products of factorizations 
with degenerate $w$ (when $R/(w)$ does not have an 
isolated singularity at $0$).  
 
\begin{definition} A polynomial $p\in \mf{m}^2$ is called a 
\emph{degenerate potential} if the Jacobian algebra $R_w$ is 
infinite-dimensional. 
\end{definition} 

Any polynomial in $\mf{m}^2$ is either a potential or a 
degenerate potential (but never both). Unless specified 
otherwise, $w$ denotes a potential. 
 
 \vspace{0.1in} 

{\bf The Tyurina algebra.} The Tyurina algebra $R^T_w$ is defined as 
the quotient of the Jacobian algebra $R_w$ by the ideal generated 
by $w.$ The multiplication by $w$ is homotopic to $0$ in 
 $\Hom_{\MF}(M,M),$ for any factorization $M,$ therefore, 
 $\Hom_{\HMF}(M,N)$ and $\Ext(M,N)$ are modules over the Tyurina 
algebra (this is a slight extension of proposition~\ref{action-fd}). 

 \vspace{0.1in} 

{\bf Cohomology of factorizations} The quotient $M/\mf{m}M$ is a 2-complex of vector 
spaces (the square of the differential is $0$ since $w\in \mf{m}$): 
 \begin{equation*} 
   M^0/\mf{m}M^0 \stackrel{d}{\lra} M^1/\mf{m}M^1 
   \stackrel{d}{\lra} M^0/\mf{m}M^0. 
 \end{equation*} 
Denote the cohomology of this 2-complex by $H(M)$ (it is 
$\Z_2$-graded, $H(M)= H^0(M) \oplus H^1(M)$) and call it 
\emph{the cohomology of}  $M.$ Cohomology of factorizations is a functor from 
$\MF^{all}_w$ and $\HMF^{all}_w$ to the category  of $\Z_2$-graded 
$\Q$-vector spaces and grading-preserving linear maps. 

\begin{prop} \label{contract} The following conditions on 
$M\in \MF^{all}_w$ are equivalent 
 \begin{enumerate} 
  \item $H(M)= 0.$ 
  \item $H^0(M) = 0.$ 
  \item $H^1(M) = 0.$ 
  \item $M$ is isomorphic to the zero factorization in 
  the category $\HMF^{all}_w.$ 
  \item $M$ is isomorphic in $\MF^{all}_w$ to (possibly infinite) 
  direct sum of 
  \begin{equation}\label{con-one}  
     R\stackrel{1}{\lra}  R \stackrel{w}{\lra}  R 
  \end{equation} 
  and  
    \begin{equation}\label{con-two} 
      R\stackrel{w}{\lra}  R \stackrel{1}{\lra}  R.  
    \end{equation} 
   \end{enumerate} 
  \end{prop} 

\emph{Proof:} implications $5 \to 4 \to 1 \to 2 $ and $1\to 3$ are 
obvious. It suffices to establish $2 \to 5.$ Choose $ R$-module 
bases in $M^0$ and $M^1$ and write factorization maps in matrix 
form. If one of the entries does not lie in the maximal ideal 
$\mf{m}\subset  R,$ this entry is invertible, and a change of bases 
makes the matrices block-diagonal with two blocks, one of which is 
either $(1,w)$ or $(w,1).$ Therefore, $M$ is isomorphic to the 
direct 
sum of some factorization $N$ and either (\ref{con-one}) or 
(\ref{con-two}). Applying Zorn's lemma, we conclude that any 
factorization $M$ decomposes into direct sum 
$M\cong M_{es}\oplus M_c$ where $M_c$ is a direct sum of 
factorizations (\ref{con-one}) or (\ref{con-two}) and $M_{es}$ does 
not have any invertible elements in its matrix presentation 
(equivalently, $M_{es}$ does not contain any summands isomorphic to 
(\ref{con-one}) or (\ref{con-two})). Since 
$M^0_{es}\cong V\otimes_{\Q}  R,$ as $ R$-module, for some vector 
space $V,$ we have  
  \begin{equation*} 
    H^0(M) \cong H^0(M_{es}) \cong V. 
  \end{equation*} 
Therefore, if $M$ satisfies property $2,$ then $V=0$ 
and $M_{es}=0.$ Proposition follows.  $\square$ 

A factorization is called \emph{contractible} if it satisfies one of 
the five equivalent conditions in Proposition~\ref{contract}. 

Define the dimension $\dim(M)$ of factorization $M$ 
as the dimension of the $\mathbb{Q}$-vector space $H(M).$ 
We have $\dim(M)\le 2\hspace{0.05in}\mathrm{rk}(M).$ If $M$ is 
a finite rank factorization, the equality holds iff $M$ does not 
contain contractible summands. 

 \vspace{0.05in} 

\emph{Example:} A factorization of rank $1$ has the form 
 \begin{equation*} 
    R \stackrel{a_0}{\lra}  R \stackrel{a_1}{\lra}  R
 \end{equation*} 
for $a_0,a_1 \in  R$ such that $w=a_0a_1.$ 
Denote this factorization by $N.$  Then 
 \begin{equation*} 
  \Hom_{\MF}(N,N) \cong  R \hspace{0.1in} \mathrm{and} 
 \hspace{0.1in}  \Hom_{\HMF}(N,N)\cong  R/(a_0,a_1),
 \end{equation*} 
the quotient of $ R$ by the ideal generated by $a_0$ and $a_1.$ 
Factorization $N$ is contractible iff either $a_0$ or $a_1$ is 
invertible (does not lie in $\mf{m}$). It's nontrivial when 
both $a_0$ and $a_1$ belong to the maximal ideal (then 
$H(M)\cong \Q\oplus \Q$). This is only possible if $k$, the number 
of variables in the set $x,$ is at most $2$ (since we need 
$w=a_0a_1$ to be non-degenerate).  
  
\begin{prop} \label{why-der} The  following properties of a 
 morphism $f:M\to N$ of factorizations are equivalent:  
  \begin{enumerate} 
  \item  $f$ is an isomorphism in the homotopy category 
    of factorizations. 
  \item $f$ induces an isomorphism between cohomologies 
    of $M$ and $N.$ 
  \end{enumerate} 
\end{prop} 
\emph{Proof:} the implication  $1\to 2$ is obvious. To prove $2\to 1,$ 
choose decompositions $M\cong M_{es}\oplus M_c$ and 
$N\cong N_{es}\oplus N_c.$ The composition
 \begin{equation*}  
  f_{es}\hspace{0.05in}: \hspace{0.1in}  
  M_{es} \stackrel{i}{\lra} M \stackrel{f}{\lra} N 
  \stackrel{p}{\lra} N_{es}, 
 \end{equation*} 
where $i$ and $p$ are the obvious inclusion and projection, induces 
an isomorphism $H(f_{es}): H(M_{es}) \cong H(N_{es}).$ Given two free
$R$-modules $L_1$ and $L_2,$ an $R$-module map $L_1\to L_2$ which 
induces an isomorphism of quotients $L_1/\mf{m}L_1 \cong 
 L_2/\mf{m}L_2$ is an isomorphism of $R$-modules. Therefore, 
 $f_{es}:M_{es}\lra N_{es}$ is an isomorphism of $ R$-modules, 
and $f$ is an isomorphism in the homotopy category of 
factorizations.  $\square$ 

\emph{Remark:} We see from this proposition that, although $\HMF^{all}_w$ 
was defined as a "homotopy" category, it also has the flavor of a 
"derived" category: in the latter, a morphism which induces an 
isomorphism in cohomology is an isomorphism. For an explanation and 
generalizations of this phenomenon see Buchweitz [B] and Orlov [O]. 

 \vspace{0.05in} 

Proof of proposition~\ref{why-der} implies  

\begin{corollary} \label{cor-unique} A decomposition 
 $M\cong M_{es}\oplus M_c$ of $M$ into a factorization without 
 contractible summands and a contractible factorization is unique 
 up to an isomorphism. 
\end{corollary}  

We call $M_{es}$ \emph{the essential summand} of $M.$ 

\begin{corollary}\label{cor-fin}  Any factorization $M$ with 
 finite-dimensional cohomology is a direct sum of a finite rank 
 factorization and a contractible factorization. 
\end{corollary} 

 \vspace{0.05in} 

Let $\MF_w$ be the category whose objects are factorizations with 
finite-dimensional cohomology and morphisms are homomorphisms of 
factorizations. Let $\HMF_w$ be the quotient of $\MF_w$ by the 
ideal of null-homotopic morphisms. Any finite rank factorization 
has finite-dimensional cohomology, therefore, we have (full and 
faithful) inclusions of categories $\MF'_w \subset \MF_w$ and 
$\HMF'_w\subset \HMF_w.$ Corollaries \ref{cor-unique} and 
\ref{cor-fin} imply  
 
\begin{corollary} The inclusion $\HMF'_w\subset \HMF_w$ is 
  an equivalence of categories. 
\end{corollary} 

Although we will be doing most of our work in categories $\HMF_w,$ 
for various $w,$ the other five categories provide a useful 
supporting framework.  To reduce the confusion of six 
different categories, we arranged them into a table. 

\vspace{0.1in} 

 \begin{tabular}{|l|c|c|c|}  \hline
        &     &          &  factorizations with \\
        & all & finite   &  finite-dimensional  \\
        & factorizations & rank & cohomology \\  \hline 
  all           &               &           &   \\
  homomorphisms &  $\MF^{all}_w$ &  $\MF'_w$  & $\MF_w$  \\ 
                &               &           &       \\  \hline
  modulo        &               &           &       \\
 homotopic to $0$  & $\HMF^{all}_w$ &  $\HMF'_w$ & 
  $\cong$ \hspace{0.1in} $\HMF_w$ \hspace{0.3in}  \\
                &               &           &   \\  \hline      
  \end{tabular} 

  \vspace{0.1in} 

Categories in the bottom row are triangulated, and the two 
rightmost categories in this row are equivalent. Even though 
most of the time our constructions start with finite rank 
factorizations, our functors take them to infinite rank 
factorizations, but with finite-dimensional cohomology. The 
latter factorizations lie in $\HMF_w$ (they can be reduced to 
finite rank factorizations, but this operation, sometimes 
necessary on concrete factorizations, is awkward from the 
categorical viewpoint). At the same time, the collection of 
categories $\HMF_w,$ over various $w,$ is closed under all 
functors that are essential in our work, and is a natural and 
thrifty choice.   

 \vspace{0.1in} 

{\bf Dualities.}  The free $ R$-module $M^{\ast}= \Hom_R(M, R)$
admits a factorization 
 \begin{equation*} 
   (M^0)^{\ast}  \stackrel{(d^1)^{\ast}}{\lra} 
   (M^1)^{\ast}  \stackrel{(d^0)^{\ast}}{\lra} 
   (M^0)^{\ast}  
 \end{equation*} 
An inclusion of factorizations $M\subset M^{\ast\ast}$ is an 
isomorphism if $M$ has finite rank, and an isomorphism in 
the homotopy category if $M$ has finite-dimensional cohomology.
Thus, $M\to M^{\ast}$ is a contravariant equvalence 
in categories $\MF_w', \HMF_w',$ and $\HMF_w.$ 

 Assume that $M$ has finite rank, and choose $R$-module 
 bases in $M^0, M^1.$ The factorization is described by 
 two matrices $(D_0,D_1).$ The dual factorization $M^{\ast}$ 
 is described by transposed matrices $(D_1^t,D_0^t).$ 

 \vspace{0.1in} 

For $M\in \MF^{all}_w$ denote by $M_-$ the factorization 
  \begin{equation*} 
    M^0 \stackrel{-d^0}{\lra} M^1 \stackrel{d^1}{\lra} 
    M^0. 
  \end{equation*} 
this assignment extends to an equivalence of categories 
of factorizations with potentials $w$ and $-w.$ 

 \vspace{0.05in} 
 
Let $M_{\bullet}= (M^{\ast})_-.$ This assignment is 
a contravariant functor between categories of 
factorizations with potentials $w$ and $-w.$ 
Observe that 
$$ \{ \bba,\bbb \}_{\bullet} \cong \{-\bbb,\bba \}$$  

{\bf Excluding a variable.} 

Let $R=\Q[[x_1,\dots, x_m]]$ and $R'=\Q[[x_2, \dots, 
 x_m]]\subset R.$ Suppose that $w\in R'$ is a potential, 
and $w=\bba\bbb$ for a pair $\{ \bba,\bbb\}$ as in 
Section~\ref{sec-koszul}, where $a_i,b_i\in R.$ 
Suppose $b_i-x_1\in R',$ for 
some $i.$ Let $c=b_i-x_1,$ and $\bba^i,\bbb^i$ be the sequences 
obtained from $\bba$ and $\bbb$ by omitting $a_i$ and $b_i.$ 
Let $\psi:R\to R'$ be the homomorphism $\psi(x_j)=x_j,$ for 
$j\not= 1,$ and $\psi(x_1)=-c$ (so that $\psi(b_i)=0$). 
Note that $\psi(w)=w.$ 

Let $\psi(\bba^i),\psi(\bbb^i)$ be the sequences obtained by 
applying $\psi$ to every entry of $\bba^i$ and $\bbb^i.$ Then 
$\{ \psi(\bba^i),\psi(\bbb^i) \}$ is an $(R',w)$-factorization. 
By treating $R$ as an $R'$-module, we can view $\{\bba,\bbb\}$ 
as an $(R',w)$-factorization (of infinite rank). Let 
$$f: \{ \bba,\bbb \} \lra \{ \psi(\bba^i),\psi(\bbb^i) \}$$ 
 be the following homomorphism of $(R',w)$-factorizations
$$\begin{array}{ccccc} \{ \bba^i,\bbb^i\} & \stackrel{a_i}{\lra}
  &   \{ \bba^i,\bbb^i\} & \stackrel{b_i}{\lra} &
  \{ \bba^i,\bbb^i\}         \\
 \downarrow\psi   &   & \downarrow  &   & \downarrow \psi  \\
  \{ \psi(\bba^i), \psi(\bbb^i) \} & \lra & 0 & \lra & 
   \{ \psi(\bba^i), \psi(\bbb^i) \}
   \end{array} $$ 
The top line is the factorization $\{\bba,\bbb\}$ written 
as the total factorization of a bifactorization. 

\begin{prop} \label{prop-exclude}
 $f$ is an isomorphism in the homotopy category of 
 $(R',w)$-factorizations. 
\end{prop} 

\emph{Proof:} It would suffice to show that $f$ induces an isomorphism 
on cohomology. Multiplication by $b_i$ is an injective 
endomorphism of the $R'$-module $\{ \bba^i,\bbb^i\}.$ 
Decompose 
 $$ \{  \bba^i,\bbb^i\} \cong b_i  \{  \bba^i,\bbb^i\} \oplus M,$$ 
where the decomposition is that of $R'$-modules, and $M$ consists 
of vectors with all coordinates in $R'$ in the standard 
$R$-module basis of $ \{  \bba^i,\bbb^i\}.$ The 
$R'$-subfactorization 
$$ b_i \{  \bba^i,\bbb^i\} \stackrel{a_i}{\lra} \{ \bba^i, 
 \bbb^i\} \stackrel{b_i}{\lra} b_i \{ \bba^i, \bbb^i\}$$ 
of $\{ \bba,\bbb\}$ is contractible, while 
$f(M)=\{ \psi(\bba^i),\psi(\bbb^i)\},$ as $R'$-modules. 
It immediately follows that $f$ induces an isomorphism on 
cohomology. $\square$ 

\vspace{0.1in} 

Let $R=\Q[[x_1,\dots, x_{k+m}]]$ and 
 $R'=\Q[[x_{k+1},\dots, x_{k+m}]]\subset R.$ Suppose that 
 $w\in R'$ is a potential, and 
$w=\bba\bbb$ for a pair $\{ \bba,\bbb\}$ with $a_j,b_j\in R.$ 
Assume that $b_1-x_1$ is a polynomial in 
$x_2, \dots, x_{k+m},$ and, more generally, $b_j-x_j,$ for 
each $j$ between $1$ and $k,$ is a polynomial in 
$x_{j+1}, \dots, x_{k+m},$ modulo the ideal generated by 
$b_1,b_2, \dots, b_{j-1}.$ Let $\psi$ be the homomorphism 
$R\to R'$ uniquely determined by the conditions 
$\psi(x_j)=x_j,$ for $j>k,$ and $\psi(b_j)=0$ for $j\le k.$ 
Let $\bba_k,\bbb_k$ be the sequences obtained from $\bba,\bbb$ 
by omitting $a_j,b_j$ for all $j\le k.$ 
Let $\psi(\bba_k),\psi(\bbb_k)$ be the sequences obtained by 
applying $\psi$ to every entry of $\bba_k$ and $\bbb_k.$ 
Then $\{ \psi(\bba_k), \psi(\bbb_k)\}$ is an 
$(R',w)$-factorization. Treating $R$ as an $R'$-module, we view 
$\{\bba,\bbb\}$ as an $(R',w)$-factorization (of infinite rank).

\begin{prop} \label{exclude-more}  $\{ \bba,\bbb\}$ 
 and $\{ \psi(\bba_k),\psi(\bbb_k)\}$ 
are isomorphic in the homotopy category of $(R',w)$-factorizations. 
\end{prop} 

Proof of the previous proposition generalizes to this setup  
without difficulty. $\square$

\vspace{0.15in} 

{\bf Frobenius structure.} The following duality theorem was 
proved by Ragnar-Olaf Buchweitz [B] in much greater generality. 

\begin{theorem} \label{thm-frobenius} 
Suppose that $w$ is a potential in even number
of variables, $w=w(x_1,\dots, x_{2k}).$ There exists a collection 
of trace maps 
$$ \mathrm{Tr}_M: \Hom_{\HMF}(M,M)\lra \Q,$$
indexed by objects of $\HMF_w,$ such that for any 
$M,N\in \mathrm{Ob}(\HMF_w)$ the composition 
$$\Hom_{\HMF}(M,N)\otimes \Hom_{\HMF}(N,M)\lra \Hom_{\HMF}(M,M)
\stackrel{\mathrm{Tr}_M}{\lra}\Q$$
is a nondegenerate bilinear pairing. 
\end{theorem}

For a proof see Theorem~7.7.5, Proposition~10.1.5, Example~10.1.6,
 and Corollary~10.3.3 in [B]. $\square$  

\vspace{0.06in} 

Although theorem~\ref{thm-frobenius} is not used explicitly in 
this paper, the Frobenius structure of categories $\HMF_w$ is 
implicit in several of our constructions and proofs. We expect 
that theorem~\ref{thm-frobenius} will become indispensable 
in further investigations of the interplay between matrix 
factorizations and link homology.


\section{Factorizations as functors} \label{sec-functors} 

{\bf Internal tensor product.} 
Let $M\in MF^{all}_w$ and $N\in MF^{all}_{-w}.$ The tensor 
product $M\otimes_R N$ is a 2-complex 
 \begin{equation*} 
   (M\otimes_R N)^0 \stackrel{d}{\lra} 
   (M\otimes_R N)^1 \stackrel{d}{\lra} 
    (M\otimes_R N)^0,
  \end{equation*} 
where 
  \begin{equation*} 
   (M\otimes_RN)^j = \oplusop{i\in \{0,1\}}
    M^i \otimes_R N^{j-i}
  \end{equation*} 
and  
  \begin{equation*} 
   d(m\otimes n)= d_M(m)\otimes n + (-1)^i m  
   \otimes d_N(n), \hspace{0.2in} m\in M^i. 
  \end{equation*} 
 
\begin{prop} $M\otimes_R N$ has finite-dimensional 
cohomology if $M$ and $N$ are factorizations with finite 
dimensional cohomology. 
\end{prop} 

\emph{Proof:} this is a special case of proposition~\ref{hmf-closed},  
stated and proved below. $\square$ 

\begin{prop}  If $M$ has finite rank, there is a natural 
isomorphism of 2-complexes 
  \begin{equation} \label{t-hom} 
    \Hom_R(M,N)\cong N \otimes_R M_{\bullet} 
  \end{equation} 
 \end{prop} 
 
\emph{Proof:} The usual isomorphism of $R$-modules on 
the left and right hand sides of (\ref{t-hom}) intertwines 
the differentials in these 2-complexes. $\square$  

\begin{corollary} If $M$ has finite-dimensional cohomology, there 
is a natural isomorphism of cohomology groups 
  \begin{equation} 
      \Ext(M,N) \cong \mathrm{H}(N\otimes_R M_{\bullet}), 
  \end{equation} 
and (in the homotopy category) of 2-complexes
 \begin{equation} 
    \Hom_R(M,N)\cong N \otimes_R M_{\bullet} 
  \end{equation} 
\end{corollary} 

\begin{figure} [htb] \drawing{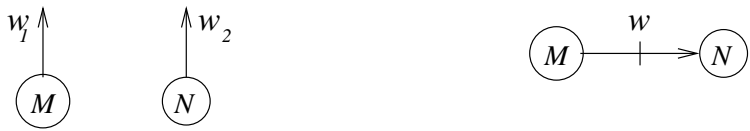}
  \caption{External and internal tensor products}
  \label{fg-tensor}
 \end{figure}

The internal tensor product $M\otimes N$ will be depicted by 
gluing the diagram's arcs, see figure~\ref{fg-tensor} right. 

 \vspace{0.1in} 

 {\bf Tensor product for sums of potentials.}  
To add dynamics to the world of matrix factorizations we 
need a large supply of functors between categories of 
factorizations over various isolated singularities $w.$ 
The tensor product (over $\Q[[y]]$) with a matrix factorization 
$M\in \MF_{w_1(x)-w_2(y)}$ is a functor from $\MF^{all}_{w_2(y)}$ 
to $\MF^{all}_{w_1(x)}.$ A slightly more general construction requires
three potentials $w_1,w_2,w_3$ and two factorizations 
$M\in \MF^{all}_{w_1(x)-w_2(y)},$   $N\in \MF^{all}_{w_2(y)- w_3(z)}.$ 
 We define their tensor product $M\otimes_y N$ by 
 \begin{equation} \label{eq-ten-gen}  
   (M\otimes_y N)^i = 
   \oplusop{j\in \{0,1\}} ( M^j \otimes_y N^{i-j}). 
  \end{equation} 
  and 
 \begin{equation*} 
  d(m\otimes n) = d_M(m) \otimes n + (-1)^i m \otimes d_N(n), 
 \hspace{0.15in} \mathrm{if} \hspace{0.15in} m\in M^i.  
 \end{equation*} 
Here $M^j\otimes_y N^{i-j}$ denotes the completed tensor 
product $M^j\widehat{\otimes}_{Q[[y]]} N^{i-j}$ (in the sense that 
the ring $\Q[[x,y,z]]$ of power series in variables $x,y,z$ is a 
completion of
$\Q[[x,y]]\otimes_{\Q[[y]]}\Q[[y,z]]$).  

 If $M$ or $N$ is contractible, so is their tensor product. 
 The tensor product can be viewed as a bifunctor 
   \begin{eqnarray*} 
    \MF^{all}_{w_1(x)-w_2(y)} \times \MF^{all}_{w_2(y)-w_3(z)}  
   & \lra & \MF^{all}_{w_1(x) - w_3(z)} \\
     \HMF^{all}_{w_1(x)-w_2(y)} \times \HMF^{all}_{w_2(y)-w_3(z)}  
   & \lra & \HMF^{all}_{w_1(x) - w_3(z)}
   \end{eqnarray*} 
 
The tensor product does not preserve the finite rank property 
(since $\Q[[x,y,z]]$ has infinite rank as a $\Q[[x,z]]$-module 
if the set of variables $y$ is nonempty). However, 

 \begin{prop}\label{hmf-closed} If $M$ and $N$ have finite 
 dimensional cohomology, so does their tensor product.
 \end{prop} 
  
 This proposition can be restated by saying that tensor product 
 restricts to bifunctors: 
   \begin{eqnarray} 
     \MF_{w_1(x)-w_2(y)} \times \MF_{w_2(y)-w_3(z)}  
   & \lra & \MF_{w_1(x) - w_3(z)}, \\
     \HMF_{w_1(x)-w_2(y)} \times \HMF_{w_2(y)-w_3(z)}  
   & \lra & \HMF_{w_1(x) - w_3(z)}. 
   \end{eqnarray}  
  
 \emph{Proof:} It suffices to show that $T=M\otimes_y N$ has 
finite-dimensional cohomology if $M$ and $N$ have finite 
rank. This cohomology $H(T)$ is a module over the Jacobian 
algebra $R_{f_2(y)},$ where $R=\Q[[y]],$ and a subquotient 
of the finitely-generated free $\Q[[y]]$-module $T/(x,z)T.$
Thus, $H(T)$ is a finitely-generated $R_{f_2(y)}$-module, and 
 necessarily has finite dimension.  $\square$ 

 \vspace{0.1in} 

\begin{figure} [htb] \drawing{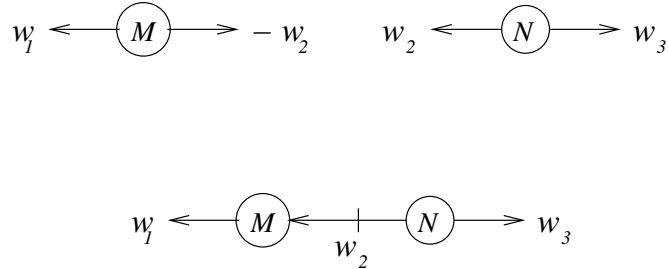}
 \caption{Graphical presentation of the tensor product over $\Q[[y]]$}
 \label{fg-tens1}
 \end{figure}

The tensor product of factorizations $M$ and $N$ will be 
depicted by joining the matching ends of their 
diagrams and placing a mark at the joint, as figure~\ref{fg-tens1} 
illustrates. 
Sometimes we will write $M$ as $M^x_y$ and $N$ as $N^y_z,$ 
and their tensor product as $M^x_y N^y_z.$ 

\begin{figure} [htb] \drawing{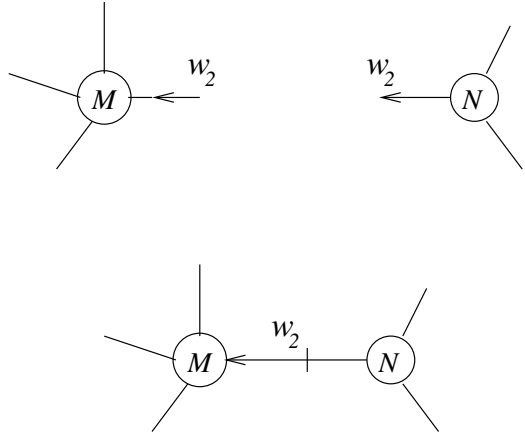}
 \caption{Tensor product over $\Q[[y]]$ when $w_1$ and $w_3$ 
 are sums}
 \label{fg-tens2}
 \end{figure}

If the summands $w_1(x)$ and $w_3(z)$ of the potentials for $M$ and 
$N$ are themselves sums of potentials, and we want to emphasize 
these decompositions, we will denote the tensor product 
$M\otimes_y N$ as in figure~\ref{fg-tens2}. 
 
\emph{Remark:} Matrix factorizations over exterior sums of 
potentials are studied in [EP], [HP], [P], and [Y2]. 

\vspace{0.1in} 

If a factorization $M$ has potential $w(x)-w(y)+w_1(z),$ 
the quotient $M/(x-y)M$ is a module 
over the ring $\Q[[x,y,z]]/(x-y).$ We treat $M/(x-y)M$
as a factorization, necessarily of infinite rank, with potential 
$w_1(z)$ over the ring $\Q[[z]]$ (were we to consider 
$M/(x-y)M$ over the larger ring $\Q[[x,z]],$ the potential 
would have been degenerate).  

\begin{prop} If factorization $M$ with potential 
$w(x)-w(y)+w_1(z)$ has finite-dimensional cohomology, so 
does $M/(x-y)M.$ 
\end{prop} 

Proof is similar to the one of proposition~\ref{hmf-closed}. $\square$ 

\begin{figure} [htb] \drawing{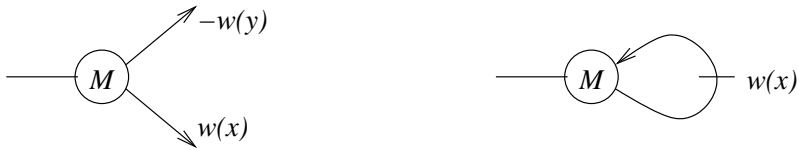}
 \caption{$M$ and its quotient $M/(x-y)M$}
 \label{fg-glue4}
 \end{figure}

We depict the 
quotient factorization by joining the $x$ and $y$ legs of $M,$ 
and placing a mark where the legs were joined, see figure~\ref{fg-glue4}. 
   
\vspace{0.1in} 

\begin{figure} [htb] \drawing{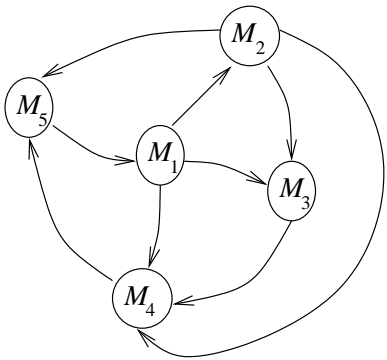}
 \caption{A closed diagram of factorizations; $d^2=0$ in the 
 tensor product, marks on the edges are omitted.}
 \label{fg-network}
 \end{figure}

A closed diagram of factorizations (as in figure~\ref{fg-network})
gives rise to a tensor product, the latter a two-periodic complex 
of modules over a suitable power series ring.  
Each oriented edge has a finite set of variables and a potential 
assigned to it, $d^2=0$ since the potentials cancel. 
If each factorization in the diagram has finite-dimensional 
cohomology, the complex will have finite-dimensional cohomology as 
well. The network does not even have to be  planar, and does 
not need to be embedded anywhere. In our paper, however, 
all such diagrams are going to be planar.  

 {\bf External tensor product.} 
When the intermediate set of variables is empty, the (completed) 
tensor product 
is over $\Q,$ and we call it \emph{the external tensor product} $M\otimes_{\Q} N.$
This operation was investigated by Yoshino [Y2]. The disjoint union 
of two diagrams denotes the external tensor product of corresponding 
factorizations, see figure~\ref{fg-tensor} left. 

 \vspace{0.1in} 

 {\bf Hiding minus signs in tensor products.} 
 Suppose $I$ is a finite set and 
 $M_a, a\in I,$ is a collection of factorizations (possibly with 
 degenerate potentials). To define the differential in the 
 tensor product $\otimesop{a\in I}M_a$ in a manifestly 
 intrinsic way, consider the Clifford ring $Cl(I)$ of the 
 set $I.$ It has generators $a\in I$ and relations 
 \begin{equation*}
  a^2=1, \hspace{0.1in}  ab+ba=0, \hspace{0.1in} a\not= b. 
 \end{equation*} 
As an abelian group, $Cl(I)$ has rank $2^{|I|},$ where $|I|$ is 
the cardinality of $I,$ and breaks down into direct sum  
   \begin{equation*} 
  Cl(I) = \oplusop{J\subset I} \Z_J.
  \end{equation*} 
 $\Z_J$ has generators--all ways to order elements of the 
set $J,$ and relations 
 \begin{equation*} 
 a \dots b c \dots e +  a\dots c b \dots e = 0
 \end{equation*} 
 for all orderings $a\dots b c \dots e$ of $J.$ The group $\Z_J$ 
 is isomorphic to $\Z,$ but there is no canonical isomorphism. 

$r_a^2=1$ where $r_a$ is the right multiplication by $a$
endomorphism of $Cl(I).$ For each $J\subset I$ which does not contain 
$a$ we have a 2-periodic sequence 
 \begin{equation*} 
   \Z_J \stackrel{r_a}{\lra}\Z_{J\sqcup \{ a\}} 
   \stackrel{r_a}{\lra} \Z_J . 
 \end{equation*} 
 Now define the tensor product factorization of $M_a$'s as 
 the sum, over all $J\subset I,$ of 
 \begin{equation*} 
  (\otimesop{a\in J}M^1_a)\otimes(\otimesop{b\in I\setminus J}M^0_b)
  \otimes_{\Z} \Z_J,  
 \end{equation*} 
with the differential 
 \begin{equation*} 
   d= \sum_{a\in I} d_a \otimes r_a 
 \end{equation*} 
 where $d_a$ is the differential in $M_a.$ We denote this tensor 
product by $\otimesop{a\in I} M_a.$ 

 When forming tensor products of factorizations, we'll use this 
 trick, and would need to assign labels to each 
term in a tensor product. If a label $a$ is assigned to 
a factorization $M,$ we write $M$ as 
 $$ M^0(\emptyset) \lra M^1(a)\lra M^0(\emptyset).$$


 {\bf Factorization of the identity functor.} 
 Start with the one-variable potential $w(x)=x^{n+1}$ and 
 consider the potential $w(x)-w(y)$ in two variables. Let  
 \begin{equation*}  
 \pi_{xy}=\frac{w(x)-w(y)}{x-y}=x^n + x^{n-1}y + \dots + y^n.
 \end{equation*} 
   Denote by $L^x_y$ the factorization 
 \begin{equation*} 
   R \stackrel{\pi_{xy}}{\lra} R \xrightarrow{x-y} R,  
 \end{equation*} 
where $R=\Q[[x,y]].$ We have  
 \begin{eqnarray} 
 \Hom_{\HMF}(L_y^x, L_y^x) & \cong & 
    \Q[x]/(x^n),  \\ 
 \Hom_{\HMF}(L_y^x, L_y^x\langle 1\rangle) & \cong & 0 . 
 \end{eqnarray} 
Indeed, $\Hom_R(L_y^x,L_y^x)$ is isomorphic to the Koszul 
complex of the sequence $(x-y, \pi_{xy}),$ with the grading collapsed 
from $\Z$ to $\Z_2.$ Regularity of this sequence implies the 
 formulas (of course, it's easy to check the above two 
formulas directly; for instance, the right hand side of 
the second formula is $0$ since $x-y$ and $\pi_{xy}$ are relatively 
 prime).  
If we assign label $(a)$ to $L^x_y,$ we can write 
 this factorization as  
 \begin{equation*} R(\emptyset) \stackrel{\pi_{xy}}{\lra} R(a) 
  \xrightarrow{x-y} R(\emptyset). 
 \end{equation*}
Here $R(\emptyset), R(a)$ are free $R$-modules of rank $1$ with 
basis vectors $1(\emptyset), 1(a).$ The differential takes 
$1(\emptyset)$ to $\pi_{xy}\cdot 1(a),$ and $1(a)$ to 
$(x-y)\cdot 1(\emptyset).$ 

\vspace{0.07in} 

\begin{figure} [htb] \drawing{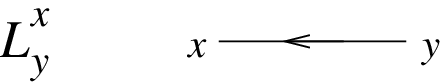}
  \caption{Arc with endpoints $x$ and $y$}
  \label{fg-arc1}
 \end{figure}

We depict $L^x_y$ by an arc oriented from $y$ to $x,$ 
see figure~\ref{fg-arc1}. If the potential assigned 
to an endpoint of a diagram has the form $x^{n+1},$ for 
some variable $x,$ we just write the variable at the endpoint, 
rather than the potential. 

 \begin{prop} \label{aisom} There is a natural isomorphism 
  \begin{equation*}  
   L_y^x M_z^y \cong M_z^x. 
  \end{equation*} 
 where $w_1(z)$ is any potential in variables $z,$ and $M_z^y$ any 
 factorization over $y^{n+1}-w_1(z).$
 \end{prop}  

(Notation $L_y^xM_z^y$ was explained several pages earlier.) 

\emph{Proof:} Assign labels $a,b,c$ to factorizations 
$L_y^x, M_z^y,$ and $M_z^x,$ respectively. 
The map of $w(x)-w_1(z)$ factorizations 
$$\tau_1: L_y^x \otimes_y M_z^y \lra M_z^x$$  
defined by taking $R(a)\otimes_y M^y_z$ to $0$ and 
$R(\emptyset) \otimes_y M^y_z\cong \Q[[x,y]]\otimes_y M^y_z$ 
onto $M^x_z$ by adding the relation $x=y$ induces an isomorphism on 
cohomologies of the two factorizations. $\square.$ 

This isomorphism is functorial in $M$ and implies 
 
\begin{corollary} Tensor product with $L_y^x$ is an invertible 
functor from the homotopy category of matrix factorizations 
with potential $y^{n+1}-w_1(z)$ to the homotopy category of 
matrix factorizations with potential $x^{n+1}-w_1(z),$ for any 
potential $w_1(z),$ and is isomorphic to the substitution functor 
that relabels $y$ into $x.$ $\square.$ 
\end{corollary} 

\begin{figure} [htb] \drawing{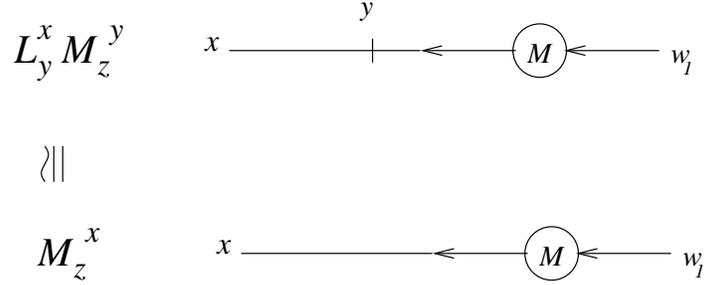}
  \caption{Removing a mark}
  \label{glue1}
 \end{figure}

We could say informally that the tensor product with $L_y^x$ 
is the identity functor. We depict $L_y^xM_z^y$ by gluing 
their "$y$"-endpoints and placing a mark at the gluing point, 
see figure~\ref{glue1} top. The proposition can be interpreted 
graphically as an isomorphism in figure~\ref{glue1}. 

\begin{figure} [htb] \drawing{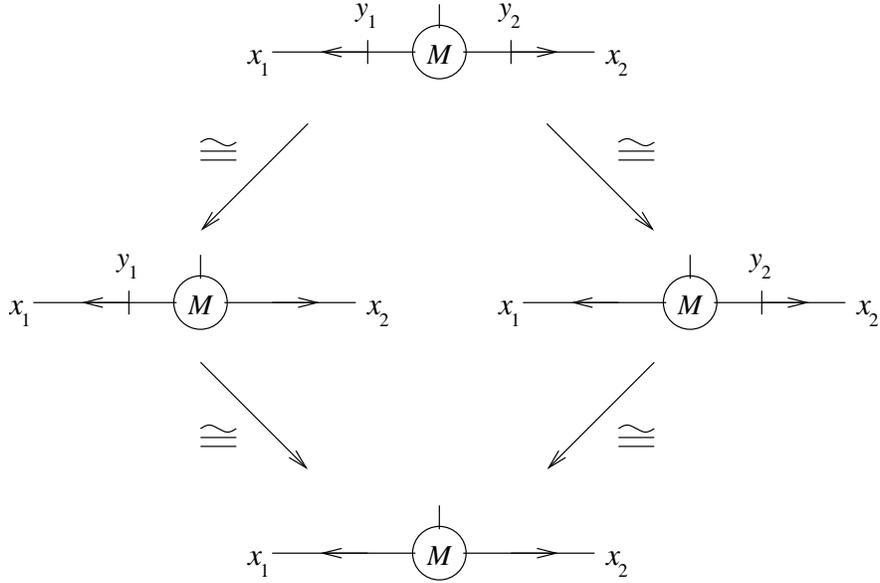}
  \caption{Commutative diagram of mark removals}
  \label{glue3}
 \end{figure}

\begin{prop} The diagram in figure~\ref{glue3} is commutative. 
 \end{prop} 

In other words, for two marks on different arcs, the order 
in which they are removed does not matter. Proof is a simple 
computation with maps $\tau_1.$ $\square$ 

 \begin{prop} \label{aisom2} There is a natural isomorphism 
  \begin{equation*}  
   M^z_x L^x_y\cong M^z_y. 
  \end{equation*} 
 where $w_1(z)$ is any potential in variables $z,$ and $M_z^x$ any 
 factorization over $w_1(z)-x^{n+1}.$
 \end{prop}  

\emph{Proof:}  Assign labels $a,b,c$ to factorizations 
$M^z_x, L^x_y,$ and $M_z^y,$ respectively. 
The map of $w_1(z)-y^{n+1}$ factorizations 
$$\tau_2: M^z_x \otimes_x L^x_y \lra M^z_y$$  
defined by taking $M^z_x\otimes_y R(b)$ to $0$ and 
$M^z_x \otimes_x R(\emptyset)\cong M^z_x \otimes_x \Q[[x,y]]$ 
onto $M^z_y$ by adding the relation $x=y$ induces an isomorphism on 
cohomologies of the two factorizations. $\square.$ 
 
\begin{figure} [htb] \drawing{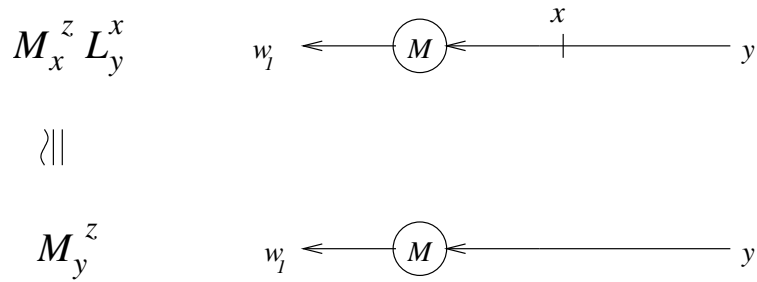}
  \caption{Removing a mark}
  \label{glue5}
 \end{figure}

The proposition can be interpreted 
graphically as an isomorphism in figure~\ref{glue5}. 

The isomorphism is functorial in $M$ and implies 
 
\begin{corollary} Tensor product with $L_y^x$ is an invertible 
functor between from the homotopy category of matrix factorizations 
with potential $w_1(z)-x^{n+1}$ to the homotopy category of 
matrix factorizations with potential $w_1(z)-y^{n+1},$ for any 
potential $w_1(z),$ and is isomorphic to the substitution functor 
that relabels $x$ into $y.$ $\square.$ 
\end{corollary} 

Commutative diagram~\ref{glue3} admits three other versions, with 
reversed orientation in one or two of the $x_1,x_2$ legs of $M$ 
and mark removal using the morphism $\tau_2$ instead of $\tau_1.$ 
Each of these three diagrams is commutative.  

\begin{figure} [htb] \drawing{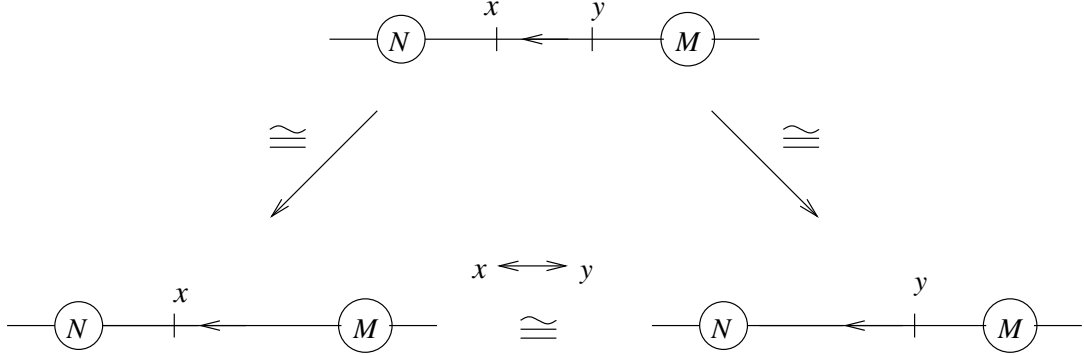}
  \caption{Yet another commutative diagram}
  \label{glue2}
 \end{figure}

In figure~\ref{glue2} we have two marks on an arc connecting 
factorizations $N$ and $M.$ We could remove the mark labelled $y$ 
using $\tau_1,$ or we could remove the other mark (via $\tau_2$), 
and then relabel $y$ into $x.$ 

\begin{prop} The two resulting morphisms from the top to the bottom 
left factorizations in figure~\ref{glue2} are equal. \end{prop} 

\emph{Proof:} direct computation. $\square$ 

There are other versions of figure~\ref{glue2}, with the orientation 
of the marked arc reversed, and with $N$ and $M$ being just one 
factorization (so that the arc starts and ends in the same 
factorization). Commutativity holds in all of these cases. 

\vspace{0.1in} 

We now specialize to the case when $M$ is itself an arc. 
The tensor product $L^x_y L^y_z$ is the factorization  
 \begin{equation} 
  \left( \begin{array}{c} R(\emptyset) \\ R(ab)  \end{array}
  \right) 
   \xrightarrow{
   \left( \begin{array}{cc} \pi_{xy}  & y-z \\
                            \pi_{yz}  & y-x \end{array}
   \right) } 
 \left( \begin{array}{c} R(a) \\ R(b)  \end{array} 
 \right) 
   \xrightarrow{
   \left( \begin{array}{cc} x-y  & y-z \\
                            \pi_{yz}  & -\pi_{xy} \end{array}
   \right) } 
  \left( \begin{array}{c} R(\emptyset) \\ R(ab)  \end{array}
  \right) 
\end{equation} 
 where $R=\Q[[x,y,z]],$ and we assigned labels $a,b$ to 
$L^x_y,L^y_z,$ respectively. The minus sign in front of $\pi_{xy}$  
 in the right $2\times 2$ matrix comes from the relation 
 $ba=-ab.$ 

$L^x_z$ is given by  
 \begin{equation*} 
   R'(\emptyset) \stackrel{\pi_{xz}}{\lra} R'(c) \xrightarrow{x-z}  R'(\emptyset)
 \end{equation*} 
 where $R'=\Q[[x,z]]$ (and notice label $c$).  

 We can specialize maps $\tau_1,\tau_2,$ introduced earlier, 
 to this case. 

 The map $\tau_1:  L^x_y L^y_z \lra L^x_z$ 
 is given by the pair of matrices  
  \begin{equation*} 
   \left(  (\phi_{y\to x}, 0) , (0,\phi_{y\to x}) \right),
  \end{equation*} 
 where $\phi_{y\to x}: R\to R'$ is the algebra homomorphism that 
 takes $x$ to $x,$ $z$ to $z,$ and $y$ to $x.$

 The map $\tau_2 :   L^x_y L^y_z \lra L^x_z$
 is given by the pair of matrices 
   \begin{equation*} 
   \left( ( \phi_{y\to z}, 0),(0,\phi_{y\to z}) \right),
   \end{equation*} 
 where $\phi_{y\to z}: R\to R'$ is the algebra homomorphism 
 that takes $x$ to $x,$ and $y,z$ to $z.$ 
  
 \begin{prop} Maps $\tau_1,\tau_2:  L^x_y L^y_z \lra L^x_z$
 are homotopic. 
 \end{prop} 

 Proof is straightforward. $\square$ 

 \vspace{0.1in} 

 From now on, an isomorphism means an isomorphism in the 
 homotopy category of factorizations, unless specified 
 otherwise.  Since $\tau_1,\tau_2$ are homotopic, they describe the same 
morphism in the homotopy category, denoted $\tau_y.$ 

Let $\tau'_y: L^x_z \lra  L^x_y L^y_z$ be given by the pair 
 of matrices  
 \begin{equation*} 
  \tau'_y = \left( \mtwobyone{1}{-e_{xyz}}, \mtwobyone{1}{1} 
  \right),  
 \end{equation*} 
 where 
 \begin{equation*} 
  e_{xyz}= \sum_{i+j+k=n-1} x^i y^j z^k.  
 \end{equation*} 
   
\begin{prop} \label{prop-inv} 
$\tau_y$ and $\tau'_y$ are mutually inverse isomorphisms 
 in the homotopy category of factorizations with potential 
 $w(x)-w(z).$ 
\end{prop} 

\emph{Proof:} it is clear that $\tau'_y$ is a homomorphism 
of factorizations, and $\tau_1\tau'_y$ is the identity endomorphism 
of $L^x_z(c).$ Proposition follows, since $\tau_y=\tau_1$ is an 
isomorphism of factorizations.  $\square$   

\begin{prop} \label{tau-assoc} 
$\tau$ is associative: there is an equality 
 $$\tau_z(\tau_y\otimes\Id)=\tau_y(\Id\otimes\tau_z)$$ 
of maps $L^x_y L^y_z L^z_w\lra L^x_w.$ 
\end{prop} 

\emph{Proof:} Direct computation. $\square$ 

\begin{corollary} $\tau'$ is associative. \end{corollary} 

\begin{corollary} \label{many-marks} 
For any $m$ and $k$ there is a canonical 
isomorphism of factorizations 
$$L^x_{z_1}L^{z_1}_{z_2}\dots L^{z_{m-1}}_{z_m} L^{z_m}_y \cong 
 L^x_{v_1}L^{v_1}_{v_2} \dots L^{v_{k-1}}_{v_k} L^{v_k}_y.$$ 
These isomorphisms are consistent. 
\end{corollary}  

\begin{figure} [htb] \drawing{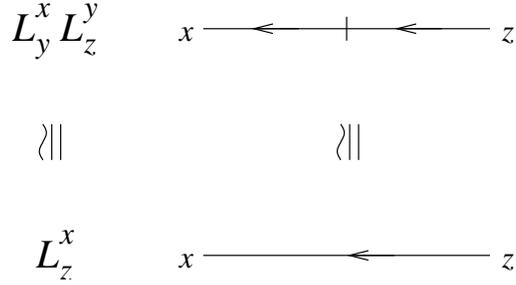}
  \caption{Adding or removing a mark}
  \label{fg-arc2}
 \end{figure}

$L^x_y L^y_z$ is depicted by two arcs glued together along matching 
enpoints, with a mark at the gluing point.
The morphism $\tau_y$ corresponds to removing a mark, and $\tau'_y$ to 
adding a mark, see figure~\ref{fg-arc2}.  

 \begin{figure} [htb] \drawing{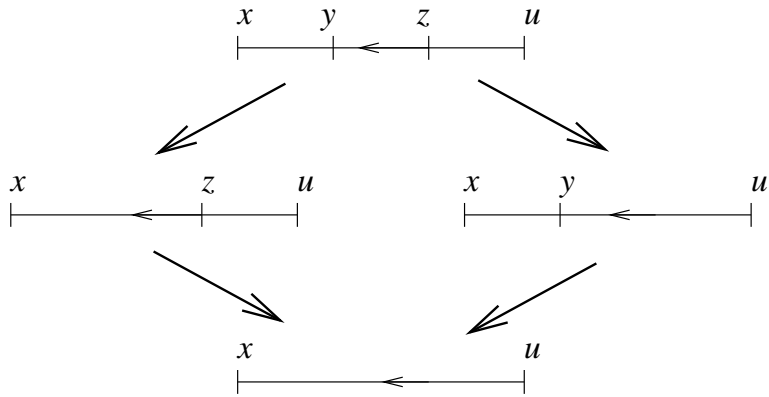}
  \caption{Associativity of mark removal}
  \label{assoc}
 \end{figure}

\begin{figure} [htb] \drawing{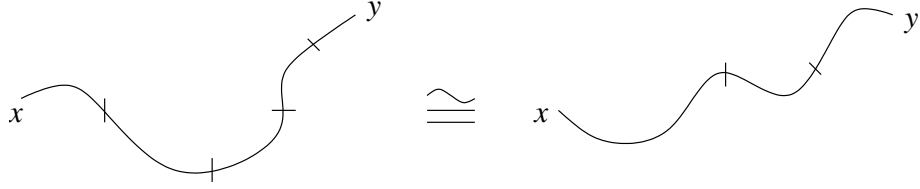}
  \caption{Canonical isomorphism of marked arcs}
  \label{fg-arc3}
 \end{figure}

Proposition~\ref{prop-inv} says that removing a marked point 
on an arc does not change the isomorphism class of a factorization, 
proposition~\ref{tau-assoc} says that arc removal is associative, 
while corollary~\ref{many-marks} asserts that two arcs, each with an 
arbitrary number of marks, are canonically isomorphic 
(figure~\ref{fg-arc3}). 
 
Denote by $L_x^x$ the quotient of $L_y^x$ by the relation 
$y=x.$  $L_x^x$ is a 2-complex of $\Q[[x]]$-modules 
$$\Q[[x]] \stackrel{\pi_{xx}}{\lra} \Q[[x]] \stackrel{0}{\lra} \Q[[x]],$$ 
and has cohomology only in degree $1.$ We depict $L_x^x$ by 
an oriented circle with one mark $x.$ 
 \begin{figure} [htb] \drawing{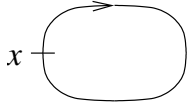}
  \caption{2-complex $L_x^x.$}
  \label{elxx}
 \end{figure}

Let 
$$A\define\mathrm{H}^{\ast}(\mathbb{CP}^{n-1},\Q)=\Q[X]/(X^n),$$
and  
$$\iota \hspace{0.1in} : \hspace{0.1in}\Q \lra A$$
be the unit map, $\iota(1)=1.$ Choose a non-zero rational 
number $\zeta,$ and define 
$$\varepsilon \hspace{0.1in} : \hspace{0.1in} A \lra \Q $$
as the trace map 
$$\varepsilon(X^{n-1})=\zeta, \hspace{0.1in} 
  \varepsilon(X^i)=0, \hspace{0.1in} \mathrm{if} \hspace{0.1in} 
  i \not= n-1.$$ 

We identify the Milnor ring $R_{w(x)}\cong \Q[x]/(x^n)$ with 
$A$ by taking $x^i\in R_w$ to $X^i \in A.$ To an oriented 
circle without marks we associate the 2-periodic 
complex $0 \lra A \lra 0,$ denoted $A\langle 1 \rangle.$ 

\begin{figure} [htb] \drawing{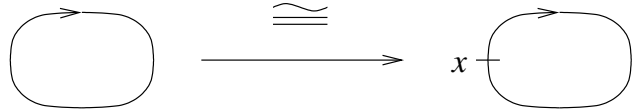}
  \caption{Adding a mark to a circle}
  \label{cisom}
 \end{figure}

We fix an isomorphism $\nu_x: A\langle 1\rangle \cong L_x^x$ of 
2-periodic complexes of vector spaces, up to homotopies, 
by taking $X^i\in A$ to $x^i\in \Q[[x]]\cong (L^x_x)^1,$
for $0\le i \le n-1.$ Graphically, this isomorphism means 
adding a mark to a circle without marks, see figure~\ref{cisom}.   
 
\begin{figure} [htb] \drawing{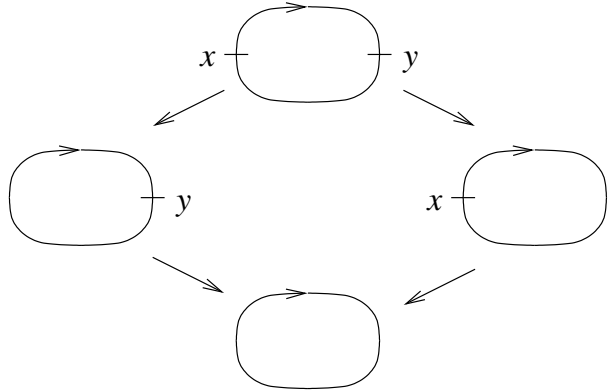}
  \caption{Commutativity of mark removal}
  \label{commutativity}
 \end{figure}

It's easy to see that the maps
 $$L^x_y L^y_x  \stackrel{\tau_y}{\lra} L_x^x  
 \stackrel{\nu^{-1}_x}{\lra} A\bracket{1} $$ 
 and 
 $$L^x_y L^y_x  \stackrel{\tau_x}{\lra} L_y^y  
 \stackrel{\nu^{-1}_y}{\lra} A\bracket{1} $$
are homotopic. This implies consistency between 
two ways to remove two marks from a circle, 
figure~\ref{commutativity}. 

\begin{figure} [htb] \drawing{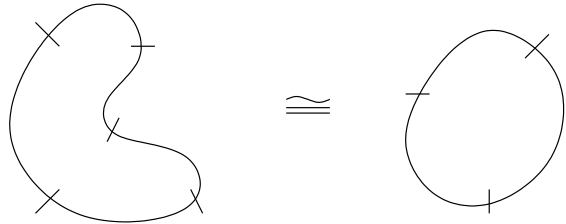}
  \caption{Canonical isomorphism of 2-periodic complexes assigned 
 to marked circles}
  \label{gencomm}
 \end{figure}

Combining with the associativity property for 
mark removal and addition, we conclude that 
the factorizations assigned to two circles with 
arbitrary number of marks are canonically 
isomorphic in the homotopy category of 2-complexes 
of $\Q$-vector spaces, figure~\ref{gencomm}. 

\emph{Remark:} 
This canonical isomorphism allows us to view $\iota$ and 
$\varepsilon$ as maps between $\Q$ and the 2-complex assigned 
to a circle with an arbitrary number of marks. 

\vspace{0.1in} 

\begin{figure} [htb] \drawing{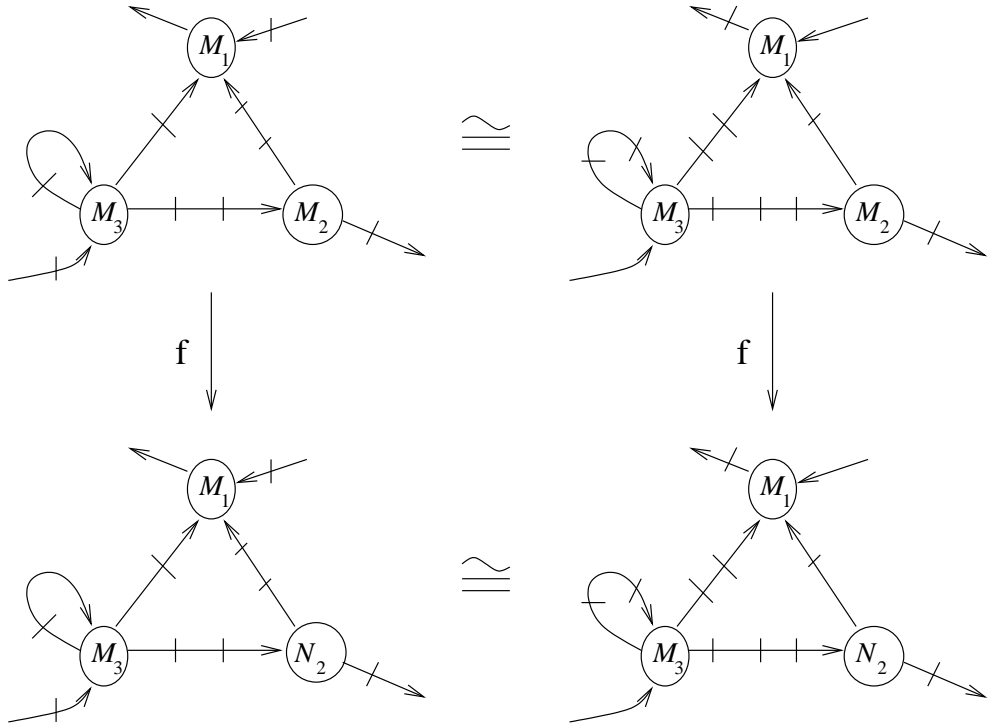}
 \caption{Canonical isomorphism of networks; a map 
 $f:M_2\lra N_2$ induces a homomorphism of factorizations assigned 
to whole networks; the diagram is commutative.}
  \label{network2}
 \end{figure}

Suppose we have a collection of factorizations $M_i, i\in I,$ 
such that each of them has potential $w_i$ which is a signed 
sum of $x_j^{n+1},$ for $j$ in a subset of $I.$ Let us tensor 
$M_i$'s and a number of arc factorization $L$ together in 
some way to produce a network of factorizations 
(as in figure~\ref{network2} top left). Each arc in the network 
has potential $x^{n+1}.$ We divide arcs into internal and external. 
External arcs are those with 
at least one loose endpoint. The resulting tensor product 
$Z_1=\otimes_{i} M_i$ is a factorization with potential which is a 
signed sum of $x_j^{n+1},$ over all loose endpoints $j.$ The sign 
is determined by the orientation of the network near $j.$

Suppose now we tensor $M_i$'s and several arc factorizations together 
so as to produce the same network but, possibly, with different 
marks (as in figure~\ref{network2} top right). Denote this 
tensor product by $Z_2.$ Our results imply 

\begin{prop}\label{pr-can-iso}
 Factorizations $Z_1$ and $Z_2$ are canonically 
isomorphic (in the homotopy category). The isomorphism is 
natural in $M_i$'s, in the following sense. If $f:M_r\lra N_r$ 
is a homomorphism of factorizations, for some $r,$ and $Z_1'$
(respectively, $Z_2'$) is obtained from $Z_1$ (respectively, $Z_2$) 
by substituting $N_r$ instead of $M_r$ in the tensor product, 
then the diagram 
\begin{equation*}
    \begin{CD}
      Z_1 @>{\cong}>> Z_2  \\
      @V{f}VV      @V{f}VV   \\
      Z_1' @>{\cong}>> Z_2' 
    \end{CD}
\end{equation*}
commutes (see figure~\ref{network2}). 
\end{prop} 

{\bf Realization of the identity functor (multi-variable case).} 
Given a multi-variable potential $w=w(x_1,\dots, x_k),$ we 
can write 
$$w(x)-w(y)=\sum_{i=1}^k w_i (x_i-y_i)$$ 
where 
$$w_i =\frac{w(y_1,\dots, y_{i-1},x_i,\dots, x_k)-
 w(y_1,\dots,y_i,x_{i+1},\dots, x_k)}{x_i-y_i}$$ 
are polynomials in $x$'s and $y$'s. Let $L$ be the 
tensor product (over $i$ from $1$ to $k$) of factorizations
$$ R \stackrel{w_i}{\lra} R \stackrel{x_i-y_i}{\lra} R$$ 
where $R=\Q[[x,y]].$ 

\begin{prop} \label{thick-arc}
For any potential $w_1(z)$ and any factorization
$M$ with potential $w(y)-w_1(z)$ the tensor product 
$L\otimes_y M$ is isomorphic in the homotopy category 
of factorizations to the factorization $M$ with $y$'s relabelled 
into $x$'s. The isomorphism is functorial in $M.$ 
\end{prop} 

To prove this result, apply proposition~\ref{prop-exclude} 
repeatedly to exclude $y_1,y_2,\dots, y_k.$  $\square$ 

Proposition~\ref{thick-arc} says, essentially, that 
tensoring with $L$ is the identity functor. 
Also, we have 
\begin{eqnarray*} 
 \Hom_{\HMF}(L,L) & \cong &  \Q[[x]]_{w(x)},  \\
 \Ext^1_{\HMF}(L,L) & \cong & 0 , 
\end{eqnarray*} 
i.e., the endomorphism ring of $L$ is the Jacobian algebra 
of $w.$ Note that theorem~\ref{thm-frobenius} applied to $L$ implies 
proposition~\ref{prop-symmetric}. 

The results of the previous subsection and of section~\ref{sec-corners} 
can be generalized from $x^{n+1}$ to arbitrary multi-variable 
potentials, but we postpone this analysis to a later paper.


\section{Homogeneous potentials and graded factorizations} 
\label{sec-homog} 

{\bf Homogeneous potentials and graded factorizations.} 
Suppose that each variable $x_i$ is given a positive integer 
degree $p_i$ and the potential $w$ is homogeneous of degree $p.$ 
Then $w$ belongs to the Jacobian ideal, 
 $$w= \frac{1}{p}\sum_i p_i x_i \frac{\partial f}{\partial x_i}$$ 
(Euler's formula), so that the Milnor and Tyurina algebra
are isomorphic.   

The ring of power series with coefficients in a field is local. 
After switching from rings to graded rings, the role of the 
power series rings is played by polynomial rings whose generators 
are in positive degrees only, since these rings have only one 
maximal homogeneous ideal. From now on in this paper we work with 
homogeneous potentials and switch from the power series ring 
to the ring of polynomials. From here on $R$ is a polynomial 
ring.  

A graded (or homogeneous) factorization with a homogeneous 
potential $w$ of degree $2(n+1)$ consists of free graded $R$-modules 
$M^0, M^1$ and degree $n+1$ homomorphisms $d^0, d^1$ such 
that $d^1d^0=w$ (this implies $d^0d^1=w$), 
 $$ M^0 \xrightarrow{d^0} M^1 \xrightarrow{d^1} M^0.$$   
A homomorphism of graded factorizations is required to have 
degree $0,$ while a homotopy should have degree $-n-1.$ Each of 
the six 
categories of factorizations described in section~\ref{sec-potent} 
has a graded version, denoted in lowercase letters. For instance, 
$\hmf_w$ is the homotopy category of graded factorizations of 
$w$ with finite-dimensional cohomology. This category is triangulated. 

We denote by $\{ m\}$ the grading shift up by $m.$ Factorization 
 $M\{ m\}$ has the form 
 \begin{equation*} 
    M^0\{ m\} \stackrel{d^0}{\lra} M^1\{ m\} 
 \stackrel{d^1}{\lra} M^0\{m\}. 
 \end{equation*} 
Cohomological shift functor $\langle 1 \rangle$ does not change the 
grading of $M^0,M^1,$ and commutes with the grading shift functor 
$\{ m\}.$  
  
All results and constructions of sections \ref{sec-koszul},
\ref{sec-potent} and \ref{sec-functors} easily extend to graded 
factorizations. There is no need to complete tensor products of 
factorizations in the graded case. Factorization $\{ a,b\}$ is 
defined for homogeneous $a,b$ with $\deg(a)+\deg(b)= 2(n+1),$ 
and has the form 
 $$ R \stackrel{a}{\lra} R\{n+1-\deg(a) \} \stackrel{b}{\lra} R.$$  
We shifted the degree of the middle $R$ so that the 
differentials would have degree $n+1.$ Then, 
$$ \{a,b\}\bracket{1}=\{b,a\}\{\deg(b)-n-1\}.$$ 
Notice two different uses of curly brackets: to denote Koszul 
factorizations and to denote shifts. 

When extending subsection "Realization of the identity functor 
(one-variable case)" of section~\ref{sec-functors} to the 
graded case, we give 
variables $x,y,z$ degree $2.$ Factorization $L^x_y$ has the form 
 $$ R\xrightarrow{\pi_{xy}} R\{ 1-n\} \xrightarrow{x-y} R.$$ 
The isomorphism in proposition~\ref{aisom} has degree $0,$ 
assuming $M$ is a graded factorization. 

\vspace{0.1in} 

Cohomology $H(M)$ of a homogeneous factorization is a 
$\Z_2\oplus \Z$-graded $\Q$-vector space, 
 $$ H(M)=\oplusop{i\in \{0,1\},j\in\Z} H^{i,j}(M).$$
We define the graded dimension of $M$ as 
$$\mathrm{gdim}(M)\define \sum_{j\in \Z, i\in \{ 0,1\}} 
  \mathrm{dim}H^{i,j}(M)\hspace{0.06in} q^j s^i.$$ 
Then, for external tensor product $M\otimes N,$ 
$$\mathrm{gdim}(M\otimes N) = \mathrm{gdim}(M)\hspace{0.04in}
 \mathrm{gdim}(N),$$ 
assuming the relation $s^2=1.$ The dimension of $M$ is the 
$q=s=1$ specialization of the graded dimension. 

\vspace{0.07in}

The algebra $A$ (cohomology of $\mathbb{CP}^{n-1}$), defined in the 
previous section, is naturally 
graded. We shift the grading of $A$ down by $n-1$ so that $\deg(1)=1-n,$
and, more generally, $deg(X^i)= 2i+1-n.$  
The unit map $\iota: \Q \to A$ and the trace map 
$\varepsilon: A\to \Q$ have degree $1-n,$ while the multiplication in 
$A$ has degree $n-1.$ 
To a circle without marks we assign the 2-complex 
$A\bracket{1}$ of graded vector spaces.

\vspace{0.06in} 

Let $M$ be a homogeneous factorization of finite rank. 
If the endomorphism ring $\mathrm{Hom}_{\mathrm{mf}}(M,M)$ 
contains a degree $0$ idempotent $e,$ we can decompose $M$ into 
a direct sum of factorizations, $M= eM \oplus (1-e)M.$ 
Suppose instead that the quotient ring $\Hom_{\hmf}(M,M)$ 
contains an idempotent $e.$ Let $I$ be the kernel of the quotient
map 
$$\Hom_{\mathrm{mf}}(M,M)\stackrel{f}{\lra}\Hom_{\hmf}(M,M),$$
and $J$ the ideal in $\Hom_{\mathrm{mf}}(M,M)$ of 
endomorphisms that induce the zero map on cohomology. 
Clearly, $I\subset J,$ and $J$ is a nilpotent ideal 
in $\Hom_{\mathrm{mf}}(M,M),$ since a degree $0$ 
endomorphism cannot have coefficients of arbitrary large 
degrees. Thus, $J^N=0,$ for some $N,$ and, therefore, 
$I^N=0.$ Nilpotent ideals have the lifting idempotents
property (exercise, or see [Be, Theorem 1.7.3]). There exists 
an idempotent 
$\widetilde{e}\in \Hom_{\mathrm{mf}}(M,M)$
that lifts $e,$ that is $f(\widetilde{e})=e.$ It allows 
us to decompose $M= \widetilde{e}M\oplus (1-\widetilde{e})M.$ 
We state this as 

\begin{prop} \label{prop-s-prop} 
 The category $\hmf_w$ has the splitting idempotents property. 
\end{prop} 

\emph{Remark:} The category $\HMF_w$ has this property as well, 
which can be seen by a slight modification of the above argument,
using that $\cap_{N} J^N=0$ and lifting $e$ recursively. 

An additive category is called \emph{Krull-Schmidt} if any object 
has the unique decomposition property. 
In other words, if $M \cong \oplusop{i\in I} M_i$ and 
$M\cong \oplusop{j\in J} N_j,$ for some sets $I,J$ and 
indecomposables $M_i,N_j,$ then there is a bijection $z:I\to J$ 
such that $M_i\cong N_{z(i)}.$ 

\begin{prop} \label{prop-krsh} The category $\hmf_w$ is 
Krull-Schmidt. \end{prop} 

Indeed, any object of $\hmf_w$ is isomorphic to a finite 
rank factorization. This and proposition~\ref{prop-s-prop} 
imply that the endomorphism ring of any indecomposable in 
 $\hmf_w$ is local. Proposition~\ref{prop-krsh} follows. $\square$ 

\vspace{0.06in} 

For graded factorizations $M,N$ with finite-dimensional 
cohomology there is an isomorphism of graded $R$-modules  
  \begin{equation*} 
  \Hom_{\HMF}(M,N) = \oplusop{i\in \Z} \Hom_{\hmf}(M\{i\}, N).
 \end{equation*}

Given a graded cyclic Koszul factorization $\{\bba,\bbb\},$ 
its graded dual is 
\begin{equation}
  \{ \bba,\bbb\}_{\bullet} \cong \{ -\bbb,\bba\}
  \cong \{\bba,-\bbb\} \bracket{m}\{ \sum_i \deg(a_i) -m(n+1)\},
\end{equation}  
where sequences $\bba$ and $\bbb$ have length $m.$  

 \vspace{0.1in} 

{\bf Quasi-homogeneous potentials.} A potential that is homogeneous with respect to 
some basis in $R$ is called \emph{quasi-homogeneous}. A quasi-homogeneous 
potential lies in its Jacobian ideal. From the work of K.Saito 
[Sa] we know that the converse is true: only quasi-homogeneous 
potentials lie in their Jacobian ideal (so that the Milnor and 
Tyurina algebras are isomorphic). To a mathematician, 
quasi-homogeneous singularities are distinguished 
 by the existence of the associated Frobenius manifold [BV], [D], 
and other special properties [AGV], [Di], [Ku]. 
From a string theorist's viewpoint, each 
quasi-homogeneous singularity gives rise to the rich structure of a 
(super) conformal 2D field theory [VW], [M], while an arbitrary 
isolated singularity only produces a 2D topological field theory, the 
latter equivalent to a commutative Frobenius algebra (a 0-dimensional 
Gorenstein ring). 

In the rest of the paper we are dealing exclusively with 
homogeneous potentials.


 \section{Planar graphs and factorizations} 
\label{sec-planar} 

{\bf Factorization from a planar graph.} 
We consider graphs of a particular kind embedded in a disk
(for an example see figure~\ref{diagbound}). 
A graph can have both unoriented and oriented edges, and oriented 
loops. Unoriented edges are called "wide" and depicted correspondingly. 
Any unoriented edge has two oriented edges entering it at one vertex, 
and two oriented edges leaving it at the other. Oriented edges might 
end on the boundary of the disc. Inside the disk we allow only trivalent 
vertices where a wide edge and a pair of oriented edges meet.  
An oriented edge is \emph{internal} if none of its endpoints is on 
the boundary of the disc. Otherwise, the edge is called \emph{external}. 
Any internal 
edge has one or more marks placed on it. We also treat boundary points 
as marks, this ensures that each external edge has a mark. 
Additional marks on external edges are allowed. An oriented loop 
may or may not have marks. 

If $\Gamma$ is such a graph, we denote by $m(\Gamma)$ the set of its
marks and by $\partial \Gamma$ the set of boundary points, the latter a 
subset of $m(\Gamma).$ If $i\in \partial\Gamma,$ the sign of $i,$
denoted $s(i),$ is $1$ if the edge at $i$ is oriented outward, 
and $-1$ if the edge is oriented inward. For instance, boundary 
points marked $1,2,7$ in figure~\ref{diagbound} have sign $1,$ 
and points $4,8,9$ have sign $-1.$ 

Let $R=\Q[x_i], i\in m(\Gamma),$ be the ring of polynomials in 
variables $x_i,$ over all marks $i,$ and $R'$ be its subring 
$\Q[x_i], i\in \partial\Gamma.$ We introduce a grading on $R$ and 
$R'$ by giving each $x_i$ degree $2.$ 

Assign to $\Gamma$ the potential 
$$w(\Gamma)= \sum_{i\in \partial \Gamma} s(i) x_i^{n+1}.$$ 

To $\Gamma$ we now associate a graded factorization 
$C(\Gamma)$ over the ring $R'$ with potential $w(\Gamma).$ 
First, to a wide edge $t$ bounded by marks $1,2,3,4$ as in 
figure~\ref{fg-xonefour} we assign a factorization with potential 
$$w_t=x_1^{n+1}+x_2^{n+1}-x_3^{n+1}-x_4^{n+1}.$$ 
Starting with formal variables $x,y,$ we can write $x^{n+1}+y^{n+1}$ 
as a polynomial in $x+y$ and $xy.$ Let $g$ be this polynomial, 
$$g(x+y, xy) = x^{n+1} + y^{n+1}.$$
Explicitly, 
 \begin{equation*} 
  g(s_1,s_2) = s_1^{n+1}+(n+1)\sum_{1\le i\le \frac{n+1}{2}} 
  \frac{(-1)^i}{i} \sbinom{n-i}{i-1} 
  s_2^i s_1^{n+1-2i}.   
 \end{equation*} 
$w_t$ can be written as
\begin{eqnarray*}
 w_t & = & g(x_1+x_2,x_1x_2)-g(x_3+x_4, x_3 x_4)  \\
  & = &  g(x_1+x_2,x_1x_2)-g(x_3+x_4,x_1x_2)+
 g(x_3+x_4,x_1x_2) - g(x_3+x_4, x_3 x_4)   \\
  & = &  \frac{g(x_1+x_2,x_1x_2)-g(x_3+x_4,x_1x_2)}{x_1+x_2-x_3-x_4}
 (x_1+x_2-x_3-x_4) \\
  & + & \frac{g(x_3+x_4,x_1x_2)-g(x_3+x_4,x_3x_4)}{x_1x_2-x_3x_4}
 (x_1x_2-x_3x_4).
\end{eqnarray*}
Let 
\begin{eqnarray}
 u_1 & = & u_1(x_1,x_2,x_3,x_4)=\frac{g(x_1+x_2,x_1x_2)- \label{eq-uone}
  g(x_3+x_4,x_1x_2)}{x_1+x_2-x_3-x_4}, \\
 u_2 & = & u_2(x_1,x_2,x_3,x_4)=\frac{g(x_3+x_4,x_1x_2) - 
   g(x_3+x_4,x_3x_4)}{x_1x_2-x_3x_4}. \label{eq-utwo} 
\end{eqnarray}
Note that $u_1$ and $u_2$ are polynomials, and 
$$w_t=u_1(x_1+x_2-x_3-x_4)+ u_2(x_1x_2-x_3x_4).$$ 
To $t$ we assign graded factorization
$$ C_t \define \{(u_1,u_2),(x_1+x_2-x_3-x_4,x_1x_2-x_3x_4)\} 
 \hspace{0.05in} \{ -1\}.$$ 
In other words, $C_t$ is the tensor product of graded factorizations
$$ R_t\xrightarrow{u_1} R_t\{1-n\}\xrightarrow{x_1+x_2-x_3-x_4}R_t$$
and 
$$ R_t\xrightarrow{u_2} R_t\{3-n\}\xrightarrow{x_1x_2-x_3x_4}R_t,$$
with the grading shifted down by $1.$ 

In general, a wide edge $t$ will be bounded by marks $i,j,k,l.$ 
Then $C_t$ is defined as above, with $1,2,3,4$ converted into 
$i,j,k,l.$ 

 \vspace{0.06in} 

If $\alpha$ is an arc in an oriented edge (or in an oriented loop) 
bounded by marks $i$ and $j$ and oriented from $j$ to $i,$ we denote 
by $L_j^i$ the factorization 
 $$ R_{\alpha} \stackrel{\pi_{ij}}{\lra} R_{\alpha} 
 \stackrel{x_i-x_j}{\lra} R_{\alpha}$$
where $R_{\alpha}=\Q[x_i,x_j]$ and 
 $$\pi_{ij}=\frac{x_i^{n+1}-x_j^{n+1}}{x_i-x_j}.$$ 
This factorization was introduced earlier, in 
section~\ref{sec-functors}, as the factorization assigned to 
an arc.  
 
\vspace{0.05in}

To an oriented loop  without marks we assign the 
2-complex $A\bracket{1}$ (see section~\ref{sec-functors}). 

\vspace{0.08in} 

Finally, we define $C(\Gamma)$ as the tensor product of 
$C_t,$ over all wide edges $t,$ of $L_j^i,$ over all 
$\alpha,$ and of $A\bracket{1},$ over all markless loops in $\Gamma.$ 
The tensor product is formed over suitable intermediate rings so 
that $C(\Gamma)$ is 
a free module of finite rank over $R.$ For instance, to form 
$C(\Gamma)$ for $\Gamma$ in figure~\ref{diagbound}, we first tensor 
$C_{t_1}$ with $L^3_5$ over the ring $\Q[x_3].$ Then tensor 
the result with $C_{t_2}$ over $\Q[x_5].$ In conclusion, tensor 
$C_{t_1}\otimes L^3_5\otimes C_{t_2}$ with $L^7_6$ over 
$\Q[x_6].$ 

\vspace{0.1in} 

{\bf Properties.} 
Clearly, $C(\Gamma)$ is a factorization with potential $w(\Gamma).$ 
We treat it as a graded factorization over the ring $R'$ of 
polynomials in boundary variables. $w(\Gamma)$ is a nondegenerate 
potential in this ring. If a graph has at least one internal mark, 
 $C(\Gamma)$ has infinite rank as an $R'$-module.  

\begin{prop} For any graph $\Gamma,$ factorization $C(\Gamma)$ 
lies in $\hmf_{w(\Gamma)}.$ 
\end{prop} 

In other words, $C(\Gamma)$ has finite-dimensional cohomology. 
This follows at once from results of section~\ref{sec-potent}. $\square$ 

\vspace{0.1in} 

Suppose that $\Gamma'$ is obtained from $\Gamma$ by placing a different 
collection of internal marks on oriented edges and loops of $\Gamma.$ 
Then the two graphs have the same potential $w(\Gamma')=w(\Gamma)$ 
assigned to them, and factorization $C(\Gamma')$ belongs to 
the category $\hmf_{w(\Gamma)}$ since the graphs share the same 
set of boundary points. 

\begin{prop} There is a canonical isomorphism in $\hmf_{w(\Gamma)}$ 
  $$ C(\Gamma') \cong C(\Gamma).$$ 
\end{prop}

 This is a special case of proposition~\ref{pr-can-iso}. \hspace{0.06in}
 $\square$

 \vspace{0.2in} 


{\bf Maps $\chi_0$ and $\chi_1.$} Consider graphs $\Gamma^0,\Gamma^1$ 
depicted in figure~\ref{pair2}. Factorization $C(\Gamma^0)$ is the 
tensor product of factorizations $L^1_4$ and $L^2_3,$ and is given by 
\begin{equation*}
  \left( \begin{array}{c}R(\emptyset) \\ R(ab)\{2-2n\}\end{array}
   \right)
 \xrightarrow{P_0} 
  \left( \begin{array}{c} R(a)\{1-n\} \\ R(b)\{1-n\}  \end{array}
 \right)
 \xrightarrow{P_1} 
 \left( \begin{array}{c}R(\emptyset) \\ R(ab)\{2-2n\}\end{array}
   \right)
\end{equation*}
where
\begin{equation*}
P_0=\left(\begin{array}{cc} \pi_{14}  & x_2-x_3 \\
      \pi_{23}  & x_4-x_1 \end{array} \right), \hspace{0.1in}
P_1=\left(\begin{array}{cc} x_1-x_4  & x_2-x_3 \\
      \pi_{23}  & -\pi_{14} \end{array} \right), 
\end{equation*} 
$$\pi_{ij}=\sum_{k=0}^n x_i^k x_j^{n-k},$$
and we assigned labels $a,b$ to $L^1_4$ and $L^2_3.$ 
Factorization $C(\Gamma^1)$ is  
\begin{equation*}
  \left(\begin{array}{c} R(\emptyset)\{-1\} \\ R(a'b')\{3-2n\} 
  \end{array} \right) \xrightarrow{Q_1}
 \left( \begin{array}{c} R(a')\{-n\} \\ R(b')\{2-n\}  \end{array}
  \right)             \xrightarrow{Q_2}
   \left(\begin{array}{c} R(\emptyset)\{-1\} \\ R(a'b')\{3-2n\} 
  \end{array} \right),
\end{equation*} 
with 
\begin{eqnarray*}
 Q_1 & = & \left(\begin{array}{cc} u_1  & x_1x_2-x_3x_4 \\
         u_2  & x_3+x_4-x_1-x_2 \end{array}\right), \\
 Q_2 & = & \left(\begin{array}{cc} x_1+x_2-x_3-x_4  & 
   x_1x_2-x_3x_4 \\ u_2  & -u_1 \end{array} \right). 
\end{eqnarray*}
We assigned labels $a',b'$ to the two factorizations (with degenerate
potentials) whose tensor product is $C(\Gamma^1).$ 

A map between $C(\Gamma^0)$ and $C(\Gamma^1)$ can be described 
by a pair of $2\times 2$-matrices $(U_0,U_1).$ 
Let $\chi_0: C(\Gamma^0)\lra C(\Gamma^1)$ be given by the pair 
\begin{eqnarray*}
U_0 & = & \left(\begin{array}{cc} x_4-x_2+\mu(x_1+x_2-x_3-x_4) & 
    0 \\ a_1  & 1 \end{array} \right),  \\
U_1 & = & \left(\begin{array}{cc} x_4 + \mu(x_1-x_4) &  
   \mu(x_2-x_3)-x_2 \\ -1 & 1 \end{array}\right),
\end{eqnarray*}
where 
$$a_1 = (\mu-1) u_2+\frac{u_1+x_1u_2-\pi_{23}}{x_1-x_4}$$ 
and $\mu\in \Z.$ Different choices of $\mu$ give homotopic 
maps. The map $\chi_0$ has degree $1.$ 

Let $\chi_1:C(\Gamma^1)\lra C(\Gamma^0)$ be 
given by the pair of matrices $(V_0,V_1)$:  
$$V_0=\left(\begin{array}{cc} 1 & 0 \\ a_2  & a_3 \end{array}
 \right), \hspace{0.2in} V_1=\left(
 \begin{array}{cc} 1 &  x_3+\lambda(x_2-x_3) \\
   1 & x_1+\lambda(x_4-x_1) \end{array}\right) $$
where 
$$a_2=\lambda u_2 + \frac{u_1+x_1u_2-\pi_{23}}{x_4-x_1}, 
 \hspace{0.2in} a_3=\lambda(x_3+x_4-x_1-x_2) + x_1 -x_3$$ 
and $\lambda\in \Z.$ Different choices of $\lambda$ give 
homotopic maps. The map $\chi_1$ has degree $1.$

The composition $\chi_1\chi_0$ is described by the pair 
$(V_0U_0, V_1U_1).$ Computing the products and specializing 
to $\mu=1-\lambda$ we get that 
\begin{equation}\label{eq-mult}
 V_0U_0=V_1U_1=(x_1-x_3+\lambda(x_3+x_4-x_1-x_2))\hspace{0.04in}
 \mathrm{I}, \end{equation}
where $\mathrm{I}$ is the identity $2\times 2$ matrix. 
Therefore, the composition $\chi_1\chi_0$ is homotopic 
to the multiplication by $x_1-x_3$ endomorphism of 
$C(\Gamma^0),$
 $$  \chi_1 \chi_0 = m(x_1-x_3).$$
This is obtained by setting $\lambda=0$ in (\ref{eq-mult}). Choosing 
instead $\lambda=1,$ we see that $\chi_1\chi_0$ is homotopic 
to multiplication by $x_4-x_2.$ There is no contradiction 
here, since the endomorphism of multiplication by 
$x_1+x_2-x_3-x_4$ is null-homotopic. 

Likewise, the composition $\chi_0\chi_1$ is homotopic 
to the multiplication by $x_1-x_3$ (and to the multiplication 
by $x_4-x_2$) endomorphism of $C(\Gamma^1).$ 

\vspace{0.1in}

\begin{figure} [htb] \drawing{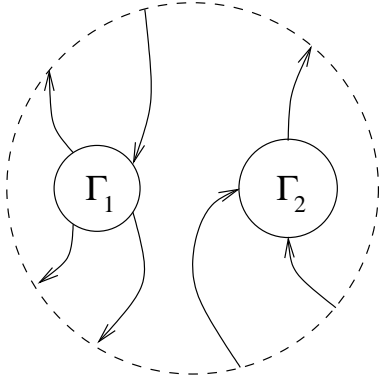}
 \caption{A disjoint union of $\Gamma_1$ and $\Gamma_2$} 
 \label{disjoint}
 \end{figure}

{\bf Disjoint union.} Given two graphs $\Gamma_1,\Gamma_2$ with potentials 
$w_1,w_2,$ the potential of their disjoint union 
$\Gamma_1\sqcup \Gamma_2$ is the exterior sum $w_1+w_2.$ 
An example of a disjoint union is depicted in figure~\ref{disjoint}. 
More often than not, this operation is not uniquely determined
by $\Gamma_1$ and $\Gamma_2,$ since we can place $\Gamma_2$ 
between any pair of adjacent exterior legs of $\Gamma_1.$ 
When $\Gamma_2$ has no exterior legs, we can place 
it inside any region of $\Gamma_1,$ including those not 
adjacent to the border of the disc. 

 \begin{prop} There is a canonical isomorphism in $\hmf_{w_1+w_2}$  
 \begin{equation*} 
   C(\Gamma_1 \sqcup \Gamma_2) \cong C(\Gamma_1) \otimes_{\Q} 
   C(\Gamma_2). 
 \end{equation*} 
 \end{prop}  

This is obvious from the definition of $C(\Gamma).$  $\square$ 

  \begin{corollary} If $\Gamma_2$ is a disjoint union of $\Gamma_1$ 
 and a loop, then  
  \begin{equation*} 
     C(\Gamma_2) \cong C(\Gamma_1)\bracket{1}\otimes_{\Q} A. 
  \end{equation*} 
  \end{corollary} 

\vspace{0.2in} 

 \begin{figure} [htb] \drawing{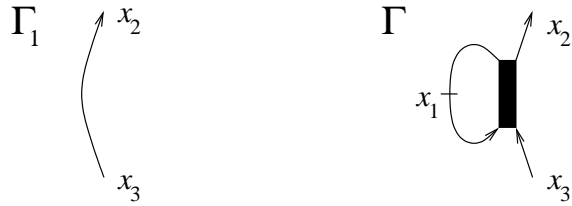}
 \caption{Graphs $\Gamma_1$ and $\Gamma$} \label{udigon1}
 \end{figure}

{\bf Direct sum decomposition I.} For graphs $\Gamma_1, \Gamma$ in figure~\ref{udigon1}, 
$C(\Gamma_1),C(\Gamma)$ are $(R',w)$-factorizations,
where $R'=\Q[x_2,x_3]$ and  $w=x_2^{n+1}-x_3^{n+1}.$ 
Notice that $C(\Gamma)$ has infinite rank as a 
$(R',w)$-factorization. It can also be viewed as a factorization 
over the larger ring $\Q[x_1,x_2,x_3]$ but with a degenerate 
potential. 

\begin{prop} \label{prop-digonm1} There is an isomorphism in $\hmf_w$ 
 $$C(\Gamma)\bracket{1} \cong 
 {\mathop{\oplus}\limits_{i=0}^{n-2}} C(\Gamma_1)\{ 2-n+2i\}. $$
\end{prop} 

\begin{figure} [htb] \drawing{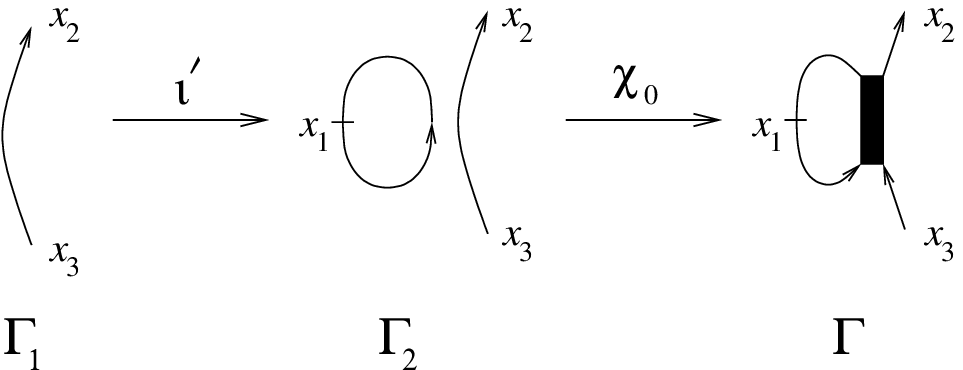}
 \caption{Map $\alpha$} \label{udigon2}
 \end{figure}

\emph{Proof:} Let $\Gamma_2$ be the disjoint union of $\Gamma_1$ 
and a circle with one mark, see figure~\ref{udigon2}. 
Define grading-preserving maps 
\begin{eqnarray*} 
 \alpha & : & C(\Gamma_1)\bracket{1} \lra C(\Gamma)\{n-2\},     \\
 \beta &  : & C(\Gamma)\{2-n\} \lra C(\Gamma_1)\bracket{1}    
\end{eqnarray*}
as follows. $\alpha$ is the composition (figure~\ref{udigon2}) 
  $$ C(\Gamma_1)\bracket{1} \stackrel{\iota'}{\lra}
   C(\Gamma_2)\{n-1\} \stackrel{\chi_0}{\lra} C(\Gamma)\{n-2\},$$
where $\iota'$ is the tensor product of the identity of 
$C(\Gamma_1)$ with the "unit" map $\iota.$

 \begin{figure} [htb] \drawing{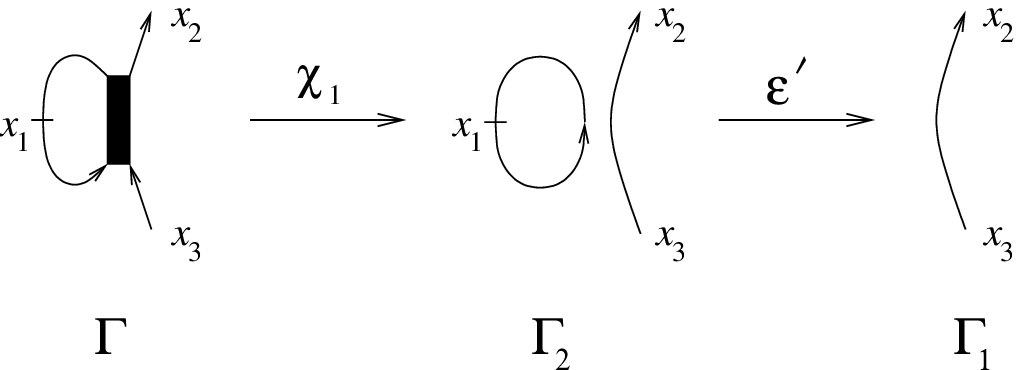}
 \caption{Map $\beta$} \label{udigon3}
 \end{figure}

$\beta$ is the composition (figure~\ref{udigon3}) 
 $$ C(\Gamma)\{2-n\} \stackrel{\chi_1}{\lra} C(\Gamma_2)\{1-n\} 
 \stackrel{\varepsilon'}{\lra} C(\Gamma_1)\bracket{1},$$
where $\varepsilon'$ is the product of the identity of $C(\Gamma_1)$ 
with the trace map $\varepsilon.$ 

Define grading-preserving maps $\alpha_i,\beta_i,$ for 
$0\le i \le n-2$ by 
 \begin{eqnarray*} 
  \alpha_i & : & C(\Gamma_1)\bracket{1} \lra C(\Gamma)\{n-2-2i\},\\
  \alpha_i & = & \sum_{j=0}^i m(x_1^j x_2^{i-j}) \alpha, \\
  \beta_i & : & C(\Gamma) \{ n-2-2i\} \lra C(\Gamma_1)\bracket{1},\\
  \beta_i & = &  \beta\hspace{0.03in} m(x_1^{n-i-2})    
 \end{eqnarray*} 
Here $m(x_1^j x_2^{i-j})$ denotes the endomorphism of $C(\Gamma)$ 
which is the multiplication by $x_1^j x_2^{i-j}.$ 

\begin{lemma} \label{lemma-delta1}
$\beta_j \alpha_i = \delta_{i,j}\zeta\mathrm{Id}$. 
\end{lemma}

\emph{Proof:} we have 
 \begin{eqnarray*}
 \beta_j \alpha_i & =& 
  \sum_{k=0}^i\beta m(x_1^{n-j-2}) m(x_1^k x_2^{i-k})\alpha
  =  \sum_{k=0}^i \varepsilon' \chi_1 \chi_0 
  m(x_1^{n-j-2}) m(x_1^k x_2^{i-k}) \iota'  \\
  & = & \sum_{k=0}^i \varepsilon' m(x_1-x_3)m(x_1^{n+k-j-2}x_2^{i-k}) 
  \iota'  \\
  & = & \sum_{k=0}^i \varepsilon' m(x_1-x_2)m(x_1^{n+k-j-2}x_2^{i-k})
  \iota'   \\
  & = &  \sum_{k=0}^i \varepsilon' m(x_1-x_2)m(x_1^{n+k-j-2}x_2^{i-k})
  \iota'   \\
  & = & \varepsilon' m(x_1^{n+i-j-1}-x_1^{n-j-2}x_2^{i+1}) \iota' 
   = \varepsilon' m(x_1^{n+i-j-1}) \iota' \\
  & = & \varepsilon(X^{n+i-j-1})\mathrm{Id} = 
    \delta_{i,j}\zeta \mathrm{Id}. 
  \end{eqnarray*} 
In the third equality we used that $\chi_1 \chi_0=m(x_1-x_3),$ and 
in the fourth that $m(x_2)=m(x_3)$ as endomorphisms of $C(\Gamma_2).$ 
Recall that $\zeta=\varepsilon(X^{n-1})$ is a nonzero rational number. 
$\square$ 

The lemma implies that the map 
\begin{equation} \label{eq-alphap}
\alpha'=\sum_{i=0}^{n-2}\alpha_i
 \hspace{0.2in} : \hspace{0.2in} 
  {{\mathop{\oplus}\limits_{i=0}^{n-2}}}C(\Gamma_1)\bracket{1} 
  \{ 2+2i-n\} \lra C(\Gamma)
\end{equation} 
is an inclusion of a direct summand, since $\beta'\alpha'=\mathrm{Id},$ 
where  
$$ \beta'=\sum_{i=0}^{n-2}\beta_i
 \hspace{0.2in} : \hspace{0.2in} C(\Gamma) \lra 
  {{\mathop{\oplus}\limits_{i=0}^{n-2}}}C(\Gamma_1)\bracket{1} 
  \{ 2+2i-n\}. $$ 
In particular, $\alpha'$ induces an injective map $H(\alpha')$ on 
cohomology 
of these factorizations. Factorization $C(\Gamma_1)$ has dimension $2,$
and the factorization on the left hand side of (\ref{eq-alphap}) 
has dimension $2(n-1).$ To finish the proof of 
proposition~\ref{prop-digonm1} it suffices to show that $H(\alpha')$ 
is bijective, which, in turn, follows from the following lemma. 

\begin{lemma} Factorization $C(\Gamma)$ has dimension $2(n-1).$ 
\end{lemma} 

\emph{Proof:} The quotient of $C(\Gamma)$ by the ideal $(x_2,x_3)$ 
is a 2-complex, and a free module of rank $4$ over 
the ring $R_1=\Q[x_1].$ It's the cyclic Koszul complex of the 
pair $\{((n+1)x_1^n,-(n+1)x_1^{n-1}),(0,0)\},$ since 
\begin{eqnarray*} 
 u_1|_{x_1=x_4, x_2=0,x_3=0} & = & (n+1)x_1^n, \\ 
 u_2|_{x_1=x_4, x_2=0,x_3=0} & = & -(n+1)x_1^{n-1}.  
\end{eqnarray*} 
By changing coordinates using a suitable automorphism $\gamma$ of 
$N=R_1\oplus R_1,$ as explained at the end of section~\ref{sec-koszul}, 
we see that this 2-complex is isomorphic to the cyclic Koszul complex 
 of the pair $\{ (x_1^{n_1},0),(0,0)\}.$ The complex 
$$\{x_1^{n-1},0\} \hspace{0.2in} : \hspace{0.2in} 
  R_1 \stackrel{x_1^{n-1}}{\lra} R_1 \stackrel{0}{\lra} R_1 $$ 
has cohomology of dimension $n-1.$ Tensoring it with $\{ 0,0\}$ 
doubles the dimension. Thus,    
 $$\dim H(C(\Gamma))=\dim H(\{ (x_1^{n-1},0),(0,0)\}) = 2(n-1).
 \hspace{0.15in}\square$$ 
 Proposition~\ref{prop-digonm1} follows. $\square$ 

 \vspace{0.2in} 

 \begin{figure} [htb] \drawing{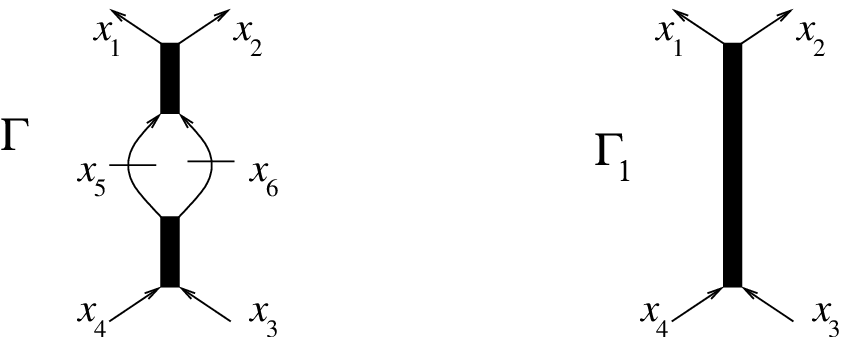}
 \caption{} \label{skein2}
 \end{figure}

 \vspace{0.2in} 

{\bf Direct sum decomposition II.} Consider factorizations $C(\Gamma),C(\Gamma_1),$ 
where $\Gamma,\Gamma_1$ are now graphs depicted in figure~\ref{skein2}. 
Both  $C(\Gamma),C(\Gamma_1)$ are $(R',w)$-factorizations, where 
 $R'=\Q[x_1,x_2,x_3,x_4]$ and 
$$w=x_1^{n+1}+x_2^{n+1}-x_3^{n+1}-x_4^{n+1}.$$
Factorization $C(\Gamma)$ has infinite rank. 

\begin{prop} \label{prop-tf} There is an isomorphism in 
$\mathrm{hmf}_w$ 
 \begin{equation*} 
    C(\Gamma) \cong C(\Gamma_1)\{ 1\} \oplus C(\Gamma_1)\{ -1\}.  
 \end{equation*} 
\end{prop} 

\emph{Proof:} Let $R=\Q[x_1,\dots,x_6],$ $s_1=x_5+x_6,$ $s_2=x_5x_6$ 
and $R_0=\Q[x_1,x_2,x_3,x_4,s_1,s_2].$ The ring $R$ is a free 
module of rank $2$ over its subring $R_0.$ 

$C(\Gamma_1)=\{ \bba,\bbb\} \{ -1\}$ where 
$\bba=(u_1,u_2), \bbb=(x_1+x_2-x_3-x_4, x_1x_2-x_3x_4),$ 
and $u_1=u_1(x_1,x_2,x_3,x_4), u_2=u_2(x_1,x_2,x_3,x_4).$ 

$C(\Gamma)=\{ \widetilde{\bba},\widetilde{\bbb}\}\{ -2\}$ where 
\begin{equation}\label{2by4matr} 
(\widetilde{\bba},\widetilde{\bbb})=\left(\begin{array}{cc} 
    u_1' & x_1+x_2-x_5-x_6 \\
    u_2' & x_1x_2 - x_5x_6 \\ u_1''  & x_5+x_6 -x_3-x_4 \\
    u_2'' & x_5x_6-x_3x_4 \end{array} \right)
\end{equation} 
and 
\begin{eqnarray*} u_1'= u_1(x_1,x_2,x_5,x_6), & 
 u_2'=u_2(x_1,x_2,x_5,x_6) \\
 u_1''=u_1(x_5,x_6,x_3,x_4), & u_2''=u_2(x_5,x_6,x_3,x_4).
\end{eqnarray*} 
The first column of (\ref{2by4matr}) lists entries of 
$\widetilde{\bba},$ the second lists entries of $\widetilde{\bbb}.$ 
All coefficients of $\widetilde{\bba}$ and $\widetilde{\bbb}$ lie 
in $R_0.$ Let $\{ \widetilde{\bba},\widetilde{\bbb}\}_0$ be the 
restriction of $\{ \widetilde{\bba},\widetilde{\bbb}\}$ 
to $R_0,$ i.e., the tensor product of factorizations 
 $$R_0\stackrel{a_i}{\lra} R_0 \stackrel{b_i}{\lra} R_0$$
where $a_i,b_i$ are entries of $\widetilde{\bba},\widetilde{\bbb}.$ 
We can decompose 
$$\{\widetilde{\bba},\widetilde{\bbb}\} \cong 
 \{\widetilde{\bba},\widetilde{\bbb} \}_0 \oplus 
\{\widetilde{\bba},\widetilde{\bbb} \}_0\{ 2\}$$
as an $(R_0,w)$-factorization, since $R\cong R_0 \oplus R_0\{2\}$ 
as graded $R_0$-modules.  Potential $w$ does not depend 
on variables $s_1,s_2$ in $R_0;$ the third entry of 
$\widetilde{\bbb}$ is $s_1-x_3-x_4$ and the fourth $s_2-x_3x_4.$ 
Applying proposition~\ref{exclude-more} to exclude $s_1$ and $s_2,$ 
we conclude that $\{\widetilde{\bba},\widetilde{\bbb} \}_0$ is 
isomorphic in the category $\hmf_w$ to $\{\bba,\bbb\}.$ 
Proposition~\ref{prop-tf} follows. $\square$

\vspace{0.2in} 

\begin{figure} [htb] \drawing{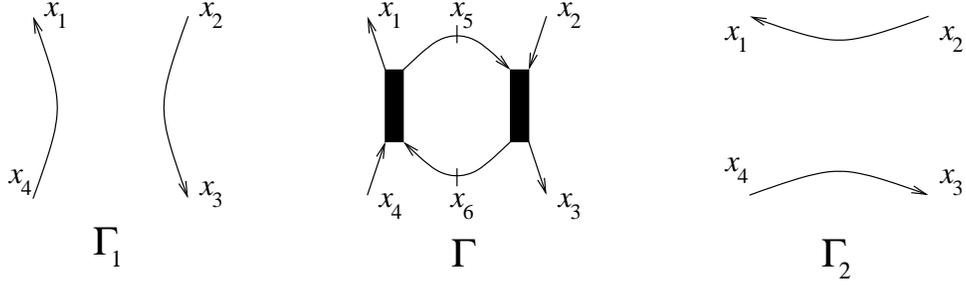}
 \caption{Graphs $\Gamma_1,\Gamma,\Gamma_2$} \label{square1}
 \end{figure}

{\bf Direct sum decomposition III.} 
Consider graphs $\Gamma, \Gamma_1, \Gamma_2$ in figure~\ref{square1}. 

$C(\Gamma),C(\Gamma_1),C(\Gamma_2)$ are homogeneous factorizations 
with the potential 
$$w= x_1^{n+1}-x_2^{n+1}+x_3^{n+1}-x_4^{n+1}$$ 
over the ring $R'=\Q[x_1,x_2,x_3,x_4].$ 

\begin{prop} \label{prop-dir-III}
 There is a direct sum decomposition in $\mathrm{hmf}_w$ 
 $$ C(\Gamma) \cong C(\Gamma_2) \oplus \biggl(
     {\mathop{\oplus}\limits_{i=0}^{n-3}}
 C(\Gamma_1)\langle 1 \rangle \{ 3-n+2i\}\biggr). $$
\end{prop} 

\begin{figure} [htb] \drawing{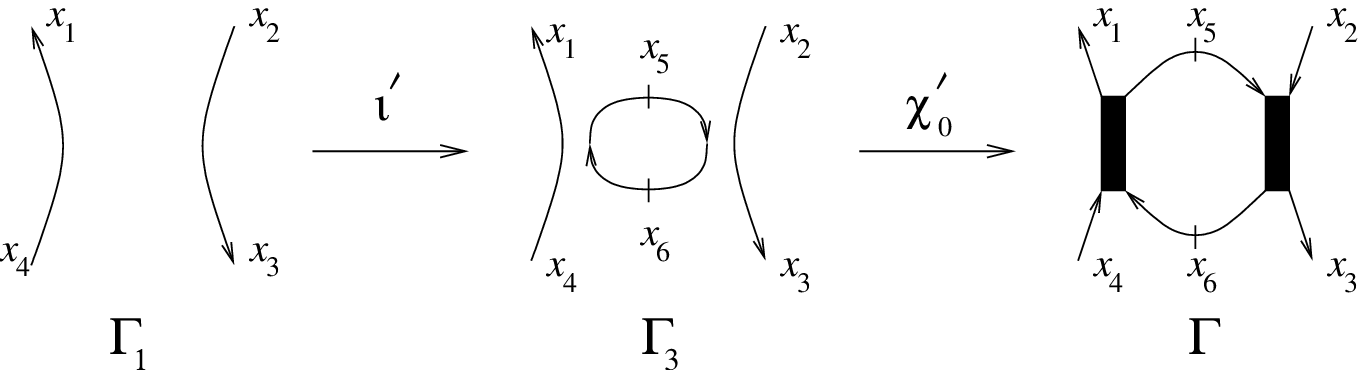}
 \caption{Map $\alpha$} \label{square2}
 \end{figure}

\emph{Proof:} define grading-preserving maps 
\begin{eqnarray*} 
\alpha & : & C(\Gamma_1)\langle 1 \rangle \lra 
 C(\Gamma)\{ n-3\},  \\
\beta & : & C(\Gamma)\{n-3\} \lra C(\Gamma_1)\langle 1 \rangle
\end{eqnarray*} 
as follows. $\alpha$ is the composition (see figure~\ref{square2})
$$C(\Gamma_1)\langle 1 \rangle \xrightarrow{\iota'} 
 C(\Gamma_3)\{ n-1\} \xrightarrow{\chi_0'} C(\Gamma)\{ n-3\}$$ 
where $\iota$ is the tensor product of the identity of $C(\Gamma_1)$ 
with the "unit" map $\iota$ from $\Q\langle 1\rangle$ to 
the factorization assigned to a circle with two marks. $\chi_0'$ is 
the composition of two $\chi_0$ maps. 

\begin{figure} [htb] \drawing{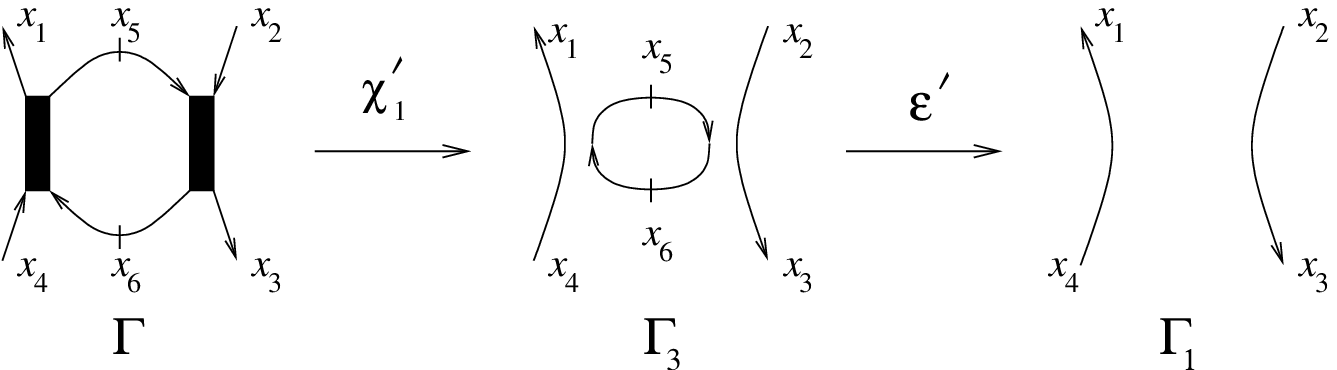}
 \caption{Map $\beta$} \label{square3}
 \end{figure}

$\beta$ is the dual composition(see figure~\ref{square3}), 
$$ C(\Gamma) \stackrel{\chi_1'}{\lra} C(\Gamma_3)\{ -2\} 
 \stackrel{\varepsilon'}{\lra} C(\Gamma_1)\langle 1\rangle \{ n-3\}$$
where $\varepsilon'$ is the tensor product of the identity of 
$C(\Gamma_1)$ with the trace map, and $\chi_1'$ is 
the composition of two $\chi_1$ maps. 

Define grading-preserving maps $\alpha_i, \beta_i,$ for $0\le i\le n-3,$
 by 
\begin{eqnarray*} 
 \alpha_i & : & C(\Gamma_1)\langle 1 \rangle \{3-n+2i\} \lra C(\Gamma)\\
  \alpha_i & = &  m(x_5^i) \alpha , \\
 \beta_i & : & C(\Gamma) \lra  C(\Gamma_1)\langle 1 \rangle \{3-n+2i\}\\
 \beta_i & = & \beta \sum_{a+b+c=n-3-i} m(x_1^a x_3^b x_5^c). 
 \end{eqnarray*}
 
\begin{lemma} $\beta_j\alpha_i = \delta_{i,j}\zeta \mathrm{Id}.$ 
\end{lemma} 

\emph{Proof:} A computation similar to the one in the proof of 
lemma~\ref{lemma-delta1}. We leave details to the reader. $\square$ 

\vspace{0.1in} 

Let $\alpha'=\sum_{i=0}^{n-3}\alpha_i$ and 
 $\beta'=\zeta^{-1}\sum_{i=0}^{n-3} \beta_i,$ considered as 
 grading-preserving maps between 
 $$N={{\mathop{\oplus}\limits_{i=0}^{n-3}}}C(\Gamma_1)\bracket{1} 
  \{3-n+2i \}$$
and $C(\Gamma).$  
Then $\beta'\alpha'=\mathrm{Id}_N,$ so that $\alpha'\beta'$ 
is an idempotent endomorphism of $C(\Gamma)$ in the category 
$\mathrm{hmf}_w.$ The splitting idempotents property in 
$\mathrm{hmf}_w$ (proposition~\ref{prop-s-prop}) implies that 
 $C(\Gamma)\cong N\oplus M$ for some graded factorization $M.$ 

\begin{lemma} \label{lemma-gdim-0}
$\gdim\hspace{0.05in} C(\Gamma)=sq^{-1}[n-1](1+sq^{1-n})(1+sq^{3-n}).$ 
\end{lemma}

\emph{Proof:} $\gdim\hspace{0.05in} C(\Gamma)$ is the graded 
 dimension of the 
cohomology of the 2-complex $C(\Gamma) / \mf{m} C(\Gamma),$ 
where $\mf{m}$ is the maximal ideal $(x_1,x_2,x_3,x_4)$ of 
$R'=\Q[x_1,x_2,x_3,x_4].$ This quotient is a free module of 
rank $16$ over the ring $R_1=\Q[x_5,x_6],$ and, after the shift by 
$\{2\},$ the cyclic Koszul complex (over $R_1$) of the 
 pair
$$(\bba,\bbb) = \left( \begin{array}{cc} u_1' &  x_5 -x_6  \\
            u_2'  &  0   \\
            u_1'' & x_6-x_5  \\
            u_2'' & 0  \end{array} \right),$$
where 
\begin{eqnarray*} 
 u_1'=u_1(x_1,x_5,x_6,x_4)|_{x_1=x_4=0},  &   & 
 u_1'' = u_1(x_3,x_6,x_5,x_2)|_{x_2=x_3=0},   \\
 u_2'=u_2(x_1,x_5,x_6,x_4)|_{x_1=x_4=0},  &   & 
 u_2'' = u_2(x_3,x_6,x_5,x_2)|_{x_2=x_3=0}. 
\end{eqnarray*} 
We can modify this pair using transformations described at the 
end of section~\ref{sec-koszul} without changing its graded dimension. 
For instance, we can add the first entry of the second column to the
third, simultaneosly with subtracting the third entry of the first 
column from the first entry of the same column. This corresponds 
to a suitable change of basis in the free $R_1$-module $R_1^{\oplus 4}.$ 
The resulting pair is 
$$ \left( \begin{array}{cc} u_1'-u_1'' &  x_5 -x_6  \\
            u_2'  &  0   \\
            u_1'' &  0   \\
            u_2'' & 0  \end{array} \right).$$ 
Moreover, $u_1'-u_1''=0,$ since $\{\bba,\bbb\}$ is a 2-complex, and 
  $x_5-x_6$ is not a zero divisor. Using proposition~\ref{prop-exclude}, 
we can cross out the first row of the above matrix and quotient the 
rest of the entries by the relation $x_5=x_6$ without changing 
the cohomology and graded dimension of this 2-complex. Thus, it 
suffices to find the graded dimension of the pair 
$$ (\bba_1,\bbb_1) = \left( \begin{array}{cc} u_2'|_{x_6=x_5}  &  0   \\
  u_1''|_{x_6=x_5} & 0 \\ u_2''|_{x_6=x_5} & 0  \end{array} \right)$$
over the ring $\Q[x_5].$ 
A direct computations tells us that 
\begin{eqnarray*}
  u_2'|_{x_6=x_5} & = & (n+1)x_5^n, \\
  u_1''|_{x_6=x_5} & = & -(n+1)x_5^{n-1},  \\     
  u_2''|_{x_6=x_5} & = & -(n+1)x_5^{n-1}. 
\end{eqnarray*} 
Adding suitable multiples of the last entry of the first column
to other entries of this column and dividing by $n+1,$ 
we simplify the 2-complex to 
$$ \left( \begin{array}{cc}   0 &  0   \\
            0 &  0   \\
            x_5^{n-1} & 0  \end{array} \right).$$
The graded dimension of this 2-complex is the product of graded 
dimensions of its rows (everything now is over the ring $\Q[x_5]$), 
since all but one rows consist of zeros. The entries $(n+1)x_5^n$ 
and $-(n+1)x_5^{n-1}$ have degrees $2n$ and $2n-2,$ respectively. 
Therefore, the graded dimension of the first row 
is $1+sq^{1-n},$ while the graded dimension of the second row 
is $1+sq^{3-n}.$ The graded dimension of $\{x_5^{n-1},0\}$ is 
 $s(q^{n-1}+q^{n-3}+\dots + q^{3-n})= sq [n-1].$ The graded 
dimension of $C(\Gamma)$ is $q^{-2}$ times the product of these 
three graded dimensions. The lemma follows. $\square.$ 

\vspace{0.15in} 

The graded dimension of $C(\Gamma_1)$ is the product of graded 
dimensions of factorization assigned to its two arcs. Thus, 
   $$ \gdim\hspace{0.05in} C(\Gamma_1)= (1+sq^{1-n})^2.$$ 
Likewise, 
   $$ \gdim\hspace{0.05in} C(\Gamma_2)= (1+sq^{1-n})^2. $$
Recall that earlier we decomposed $C(\Gamma)\cong N \oplus M.$ 
Therefore, 
   $$ \gdim\hspace{0.05in} C(\Gamma) = \gdim\hspace{0.05in} N + 
    \gdim\hspace{0.05in} M. $$
Clearly,  
$$\gdim\hspace{0.05in} N = s [n-2]\hspace{0.03in} 
 \gdim\hspace{0.05in} C(\Gamma_1).$$ 

\begin{lemma} $C(\Gamma_2)$ and $M$ have the same graded 
dimension.
\end{lemma} 

This follows from lemma~\ref{lemma-gdim-0} and the above equations. 
In the computation we use that $s^2=1.$ $\square$  

\begin{lemma} \label{lemma-grdim3} $\gdim\hspace{0.05in} 
 \Ext_{\HMF}(C(\Gamma_2),C(\Gamma))= 
 q^{2n-2} [2][n][n-1].$ 
\end{lemma} 
The ext groups here are computed in the category $\HMF_w,$ 
and are naturally $\Z\oplus \Z_2$-graded, since the 
factorizations are homogeneous. The graded dimension $\gdim$ 
is the Poincare polynomial of a $\Z\oplus \Z_2$-graded 
vector space. 

\begin{figure} [htb] \drawing{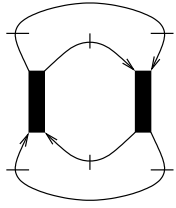}
   \caption{Graph $\Gamma_4$} \label{gamma4}
 \end{figure}

To prove the lemma, we start with the isomorphism  
$$\Ext_{\HMF}(C(\Gamma_2),C(\Gamma)) \cong H(C(\Gamma)\otimes_{R'} 
 C(\Gamma_2)_{\bullet}) \{2n-2\}$$ 
of $\Z\oplus \Z_2$-graded vector spaces, 
implying the equality of graded dimensions of the two sides. 
The 2-complex whose cohomology is computed on the right hand side 
(without the shift) is isomorphic to the 2-complex $C(\Gamma_4),$ where 
$\Gamma_4,$ shown in figure~\ref{gamma4}, is given by coupling 
$\Gamma_2$ to $\Gamma.$  

Applying propositions~\ref{prop-tf} and \ref{prop-digonm1}, we find 
that $\gdim\hspace{0.05in} C(\Gamma_4)= [2][n][n-1].$ 
 Lemma~\ref{lemma-grdim3} follows. $\square$ 

Furthermore, 
\begin{equation*} \gdim \hspace{0.05in} 
  \Ext_{\HMF}(C(\Gamma_2),C(\Gamma_1))=q^{2n-2}[n]s.
\end{equation*}  
Indeed, the left hand side, up to the shift $\{ 2n-2\},$ is 
isomorphic to the cohomology of the factorization assigned 
to the graph given by coupling $\Gamma_1$ and $\Gamma_2.$ This 
graph is a circle. Variable $s$ appears because 
$\Ext^1$ is non-trivial. 

Since factorization $N$ is the direct sum of $C(\Gamma_1)$'s 
with shifts, we can compute the graded dimension of 
$\Ext_{\HMF}(C(\Gamma_2),N),$ and, hence, of the ext groups 
$\Ext_{\HMF}(C(\Gamma_2),M),$ as the difference of two 
graded dimensions. The answer is 
\begin{equation} 
\gdim \hspace{0.05in}\Ext_{\HMF}(C(\Gamma_2),M) = q^{2n-2}[n]^2. 
\end{equation} 
Similarly, we compute that 
\begin{equation} 
 \gdim \hspace{0.05in}\Ext_{\HMF}(M,C(\Gamma_2)) = q^{2n-2}[n]^2. 
\end{equation} 
To avoid the last computation, we could invoke 
 theorem~\ref{thm-frobenius} 
and its graded version (see [B, Example~10.1.6]) to show 
that the graded dimensions of the two ext groups are equal. 

Since $q^{2n-2}[n]^2\in 1+ q\Z[q],$ we obtain the following 
corollary. 

\begin{corollary} \label{space-q} 
The vector space of morphisms between $C(\Gamma_2)$ 
and $M$ in $\mathrm{hmf}_w$ is isomorphic to $\Q.$  
\end{corollary} 

In other words, any two non-trivial grading-preserving morphisms 
from $C(\Gamma_2)$ to $M$ are multiples of each other. 
Same goes for morphisms in the opposite direction. 

After stripping off contractible summands from $M$ we obtain 
a graded factorization of rank $2$ which has the form 
$$ R'\oplus R'\{2-2n\} \lra R'\{1-n\}\oplus R'\{1-n\} 
 \lra R' \oplus R'\{2-2n\}. $$
Assume from now on that $M$ has been reduced to this minimal 
form. Then we strengthen the corollary, by noticing the absence 
of homotopies of degree $-1-n$ between $C(\Gamma_2)$ 
and $M.$ We conclude 
that the space of grading-preserving morphisms between 
$C(\Gamma_2)$ and $M$ in $\mathrm{mf}_w$ is one-dimensional 
(i.e. even before throwing out null-homotopic morphisms). 

Let $\theta$ be a nontrivial grading-preserving morphism 
from $C(\Gamma_2)$ to $M,$ and $\theta'$ a nontrivial 
grading-reserving morphism from $M$ to $C(\Gamma_2).$ Each 
of these morphisms is unique up to scaling by a non-zero 
rational number. 

\begin{lemma} The composition $\theta'\theta$ is non-zero. 
\end{lemma} 

\emph{Proof:} Recall that a morphism between two factorizations
can be described by a pair of matrices $(F_0,F_1)$ 
with polynomial coefficients. Define the rank of a morphism as 
the rank of $F_0,$ treated as a matrix over the quotient field 
of the polynomial ring. Matrices $F_0$ and $F_1$ have the same 
rank.  

Factorization $C(\Gamma_2)$ has the form 
\begin{equation*}
  \left( \begin{array}{c}R \\ R\{2-2n\}\end{array}
   \right)
 \xrightarrow{P_0} 
  \left( \begin{array}{c} R\{1-n\} \\ R\{1-n\}  \end{array}
 \right)
 \xrightarrow{P_1} 
 \left( \begin{array}{c}R \\ R\{2-2n\}\end{array}
   \right)
\end{equation*}
where 
\begin{equation*}
P_0=\left(\begin{array}{cc} \pi_{12}  & x_3-x_4 \\
      \pi_{34}  & x_2-x_1 \end{array} \right), \hspace{0.1in}
P_1=\left(\begin{array}{cc} x_1-x_2  & x_3-x_4 \\
      \pi_{34}  & -\pi_{12} \end{array} \right).  
\end{equation*}

Factorization $M$ has the form 
\begin{equation*}
  \left( \begin{array}{c}R \\ R\{2-2n\}\end{array}
   \right)
 \xrightarrow{V_0} 
  \left( \begin{array}{c} R\{1-n\} \\ R\{1-n\}  \end{array}
 \right)
 \xrightarrow{V_1} 
 \left( \begin{array}{c}R \\ R\{2-2n\}\end{array}
   \right)
\end{equation*}
where $V_0,V_1$ are matrices with homogeneous entries (of degrees
equal to degrees of the matching entries in $P_1,P_2$) such that 
$V_1V_0=w\cdot\mathrm{Id}.$ 

Note that if $\theta'\theta$ is homotopic to $0,$ then it is actually 
zero, since there are no homotopies of degree $-n-1$ from $C(\Gamma_2)$ 
to itself. Assume that $\theta'\theta=0.$ This is only possible 
if the ranks of both $\theta'$ and $\theta$ are equal to $1.$ 
Let $\theta$ be given by a pair of matrices 
$$\Theta_0 = \left( \begin{array}{cc} f_{11} & f_{12} \\ f_{21} & f_{22} 
  \end{array} \right), \hspace{0.2in} 
\Theta_1 = \left( \begin{array}{cc} g_{11} & g_{12} \\ g_{21} & g_{22} 
  \end{array} \right).$$ 
Since $\theta$ has zero degree, we see that all entries of $\Theta_1$ 
as well as diagonal entries of $\Theta_0$ are rationals, $f_{12}=0,$ 
and $f_{21}$ is a polynomial of degree $2n-2.$  

Assume that $f_{22}\not= 0.$ Then we can rescale a basis vector of $M^0$ 
so that $f_{22}=1.$ Since $\Theta_0$ has rank $1,$ we see that 
$f_{11}=0.$ By changing the basis in $M^1,$ if necessary, 
we can assume $g_{21}=0.$ Let 
$$V_0=\left( \begin{array}{cc} v_{11} & v_{12} \\ v_{21} & v_{22} 
  \end{array} \right).$$
The equation $\Theta_1 P_0= V_0\Theta_0$ implies 
$$v_{22}f_{21}=g_{22}\pi_{34}, \hspace{0.2in} 
  v_{22}=g_{22}(x_2-x_1).$$ 
This forces $g_{22}=0,$ since $\pi_{34}$ is not divisible by $x_2-x_1.$ 
Then $v_{22}=0,$ which implies that $\mathrm{det}(V_0)=-v_{12}v_{21}.$
This determinant must divide $w^2$ (which is the determinant of 
$V_1V_0.$ This is impossible since $v_{12}$ is a homogeneous linear 
function of $x_1,x_2,x_3,x_4$ and 
$$w=x_1^{n+1}-x_2^{n+1}+x_3^{n+1}-x_4^{n+1}.$$  

Therefore, $f_{22}=0.$ The equation $\Theta_1 P_0= V_0\Theta_0$ implies
$$ g_{11}(x_3-x_4)+g_{12}(x_2-x_1)=0, \hspace{0.15in} 
  g_{21}(x_3-x_4)+g_{22}(x_2-x_1)=0.$$
Since $g_{ij}$ are rational numbers, they are all zeros. 
This is a contradiction, since $\Theta_1$ should have rank $1.$  
\hspace{0.07in} $\square$  

\hspace{0.1in} 

The lemma implies that $\theta'\theta$ is a non-zero multiple 
of the identity morphism of $C(\Gamma_2).$ In addition, we know 
that $C(\Gamma_2)$ and $M$ have the same graded dimension. Therefore, 
 $C(\Gamma_2)$ and $M$ are isomorphic in $\mathrm{hmf}_w.$ 
Proposition~\ref{prop-dir-III} follows. $\square$ 

\vspace{0.3in} 

\begin{figure} [htb] \drawing{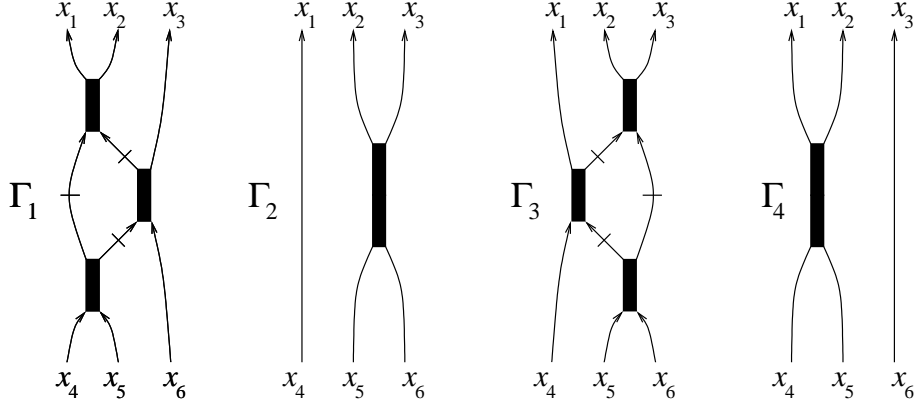}
 \caption{Graphs $\Gamma_1,\Gamma_2,\Gamma_3,\Gamma_4$} \label{fourgr}
 \end{figure}

{\bf Direct sum decomposition IV.} 
Consider graphs $\Gamma_i, 1\le i \le 4,$ depicted in 
 figure~\ref{fourgr}. Factorizations $C(\Gamma_i)$ have potential  
$$w=x_1^{n+1}+x_2^{n+1}+x_3^{n+1}-x_4^{n+1}-x_5^{n+1}-x_6^{n+1}.$$ 
We view them as factorizations over the ring $R'=\Q[x_1,\dots, x_6].$

 \begin{prop} \label{isom-dir-4} There is an isomorphism in $\hmf_w$ 
  $$ C(\Gamma_1)\oplus C(\Gamma_2)\cong C(\Gamma_3)\oplus C(\Gamma_4).$$ 
 \end{prop} 

\emph{Proof:} if $n=1,$ $C(\Gamma_i)=0$ for all $i,$ and the proposition 
is trivial. Assuming $n>1,$ we introduce another factorization. Let 
 \begin{eqnarray*} 
  s_1 & = & x_1+x_2+x_3, \\
  s_2 & = & x_1x_2+x_1x_3+x_2x_3, \\
  s_3 & = & x_1x_2x_3, \\
  s_4 & = & x_4+x_5+x_6, \\
  s_5 & = & x_4x_5 + x_4x_6 + x_5x_6, \\
  s_6 & = & x_4x_5x_6. 
 \end{eqnarray*} 
 Let $R_s\subset R'$ be the subring generated by $s_1, \dots, s_6.$ 
It's isomorphic to the polynomial ring in $s_1, \dots, s_6.$ As 
an $R_s$-module, $R'$ is free of rank $36.$  
Define a 3-variable polynomial $h$ by  
 $$h(s_1,s_2,s_3) = x_1^{n+1} + x_2^{n+1} + x_3^{n+1}.$$ 
Then 
  \begin{eqnarray*} 
  w & = & h(s_1, s_2,s_3) - h(s_4,s_5,s_6)  \\
    & = & \frac{h(s_1,s_2,s_3)-h(s_4,s_2,s_3)}{s_1-s_4} (s_1-s_4) + \\ 
   &  & \frac{h(s_4,s_2,s_3)-h(s_4,s_5,s_3)}{s_2-s_5} (s_2-s_5) + \\ 
   &  & \frac{h(s_4,s_5,s_3)-h(s_4,s_5,s_6)}{s_3-s_6} (s_3-s_6) \\
   & = & v_1 \alpha_1 + v_2 \alpha_2 + v_3 \alpha_3,   
  \end{eqnarray*} 
where 
 \begin{eqnarray*} 
  \alpha_1 & = & s_1 - s_4 = x_1+x_2+x_3-x_4-x_5-x_6, \\
  \alpha_2 & = & s_2 - s_5 = x_1x_2+x_1x_3+x_2x_3-x_4x_5-x_4x_6-x_5x_6, \\
  \alpha_3 & = & s_3 - s_6 = x_1x_2x_3-x_4x_5x_6, \\
    v_1    & = &  \frac{h(s_1,s_2,s_3)-h(s_4,s_2,s_3)}{s_1-s_4}, \\
    v_2    & = &  \frac{h(s_4,s_2,s_3)-h(s_4,s_5,s_3)}{s_2-s_5}, \\
    v_3    & = &  \frac{h(s_4,s_5,s_3)-h(s_4,s_5,s_6)}{s_3-s_6}. 
 \end{eqnarray*} 
 $v_1,v_2,v_3$ are polynomials in $x_1, \dots, x_6.$ The quotient 
of $R_s$ by the ideal $(\alpha_1,\alpha_2,\alpha_3)$ is isomorphic 
to $\Q[s_1,s_2,s_3].$ The images of $v_1,v_2,v_3$ in this 
quotient ring have the form 
$$ \frac{\partial h(s_1,s_2,s_3)}{\partial s_1}, \hspace{0.2in}  
\frac{\partial h(s_1,s_2,s_3)}{\partial s_2}, \hspace{0.2in}
\frac{\partial h(s_1,s_2,s_3)}{\partial s_3}. $$

It's well-known that the quotient of $\Q[s_1,s_2,s_3]$ by 
the ideal generated by these derivatives is isomorphic to 
the cohomology ring of the Grassmannian of 3-dimensional subspaces 
 of $\C^n$ (see [MS, page 127], for example). Thus, there is 
an algebra isomorphism 
\begin{equation} \label{cong-grass} 
  R_s/(\alpha_1,\alpha_2,\alpha_3,v_1,v_2,v_3) \cong 
  \mathrm{H}^{\ast}(\mathrm{Gr}(3,n), \Q). 
\end{equation}  
In the degenerate case $n=2$ both sides are zero. 

Likewise, the quotient ring 
$ R'/(\alpha_1,\alpha_2,\alpha_3,v_1,v_2,v_3) $ 
is the cohomology ring of the configuration space 
$$\{ N_1 \subset N_2 \subset N_3 \supset N_2'\supset N_1' | 
 \hspace{0.06in} \dim N_i = \dim N_i'=i, \hspace{0.05in} 
   N_i,N_i'\subset \C^n \}.$$

\begin{lemma} The sequence $(\alpha_1,\alpha_2,\alpha_3,v_1,v_2,v_3)$ 
 and each of its permutations are regular sequences in rings $R_s$ 
 and $R'.$ 
\end{lemma} 

\emph{Proof:} The length of the sequence equals the number of generators 
of polynomial rings $R_s$ and $R'.$ Since the quotient rings by 
the ideal generated by this sequence are finite-dimensional, the 
sequence is regular. \hspace{0.06in} $\square$ 

\vspace{0.1in} 

Let $\Upsilon$ be the following $(R',w)$-factorization: 
 $$\Upsilon \define \{ (v_1,v_2,v_3), (\alpha_1,\alpha_2,\alpha_3) \} 
  \{ -3 \}.$$ 
$\Upsilon$ is the tensor product of factorizations 
  \begin{eqnarray*} 
   R' \xrightarrow{v_1} & R'\{1-n\} & \xrightarrow{\alpha_1} R', \\
   R' \xrightarrow{v_2} & R'\{3-n\} & \xrightarrow{\alpha_2} R', \\
   R' \xrightarrow{v_3} & R'\{5-n\} & \xrightarrow{\alpha_3} R',
  \end{eqnarray*} 
shifted by $\{-3\}.$ One could think of $\Upsilon$ as the factorization 
assigned to the diagram in the lower left corner of figure~\ref{triple}. 

Proposition~\ref{isom-dir-4} will follow from 
 
 \begin{prop} \label{core-isom} There is an isomorphism in $\hmf_w$ 
   $$ C(\Gamma_1) \cong \Upsilon \oplus C(\Gamma_4).$$ 
 \end{prop}  

 Indeed, factorization $C(\Gamma_1)$ turns into $C(\Gamma_3),$ 
 and $C(\Gamma_2)$ into $C(\Gamma_4)$ if we permute $x_1$ with $x_3$ 
 and $x_4$ with $x_6.$ Factorization $\Upsilon$ is invariant under 
 this operation, since $\alpha_i$'s and $v_i$'s are, so that  
 proposition~\ref{core-isom} also implies the isomorphism 
  $$ C(\Gamma_3) \cong  \Upsilon \oplus C(\Gamma_2).$$ 
 Then both sides in the equation in proposition~\ref{isom-dir-4} 
 are isomorphic to $\Upsilon\oplus C(\Gamma_2)\oplus C(\Gamma_4).$ 
 Thus, it suffices to establish proposition~\ref{core-isom}. 

 \vspace{0.1in} 

 Proof of proposition~\ref{core-isom} occupies the next 11-12 pages. 
 
\begin{lemma} $C(\Gamma_4)$ is a direct summand of $C(\Gamma_1).$ 
 \end{lemma} 

\begin{figure} [htb] \drawing{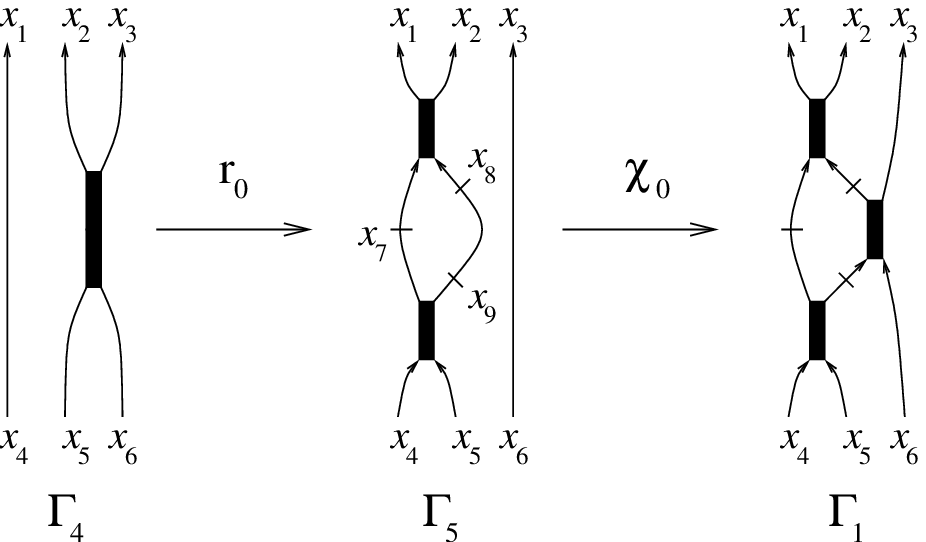}
 \caption{Map $g_0$} \label{mapgnull}
 \end{figure}

\emph{Proof:} Let 
$$g_0=\chi_0 r_0, \hspace{0.2in}
  C(\Gamma_4) \stackrel{r_0}{\lra} C(\Gamma_5) \stackrel{\chi_0}{\lra}
  C(\Gamma_1)$$ 
be the composition depicted in figure~\ref{mapgnull}.  

Proposition~\ref{prop-tf} implies that $C(\Gamma_5)\cong C(\Gamma_4)\{-1\} 
 \oplus C(\Gamma_4)\{1\},$ and $r_0$ is defined to be the inclusion of 
$C(\Gamma_4)$ into $C(\Gamma_5)$ as the direct summand 
 $C(\Gamma_4)\{-1\}.$ Maps $r_0$ and $\chi_0$ have degrees $-1$ and $1,$ 
 respectively, so that $g_0$ is grading-preserving. 

\begin{figure} [htb] \drawing{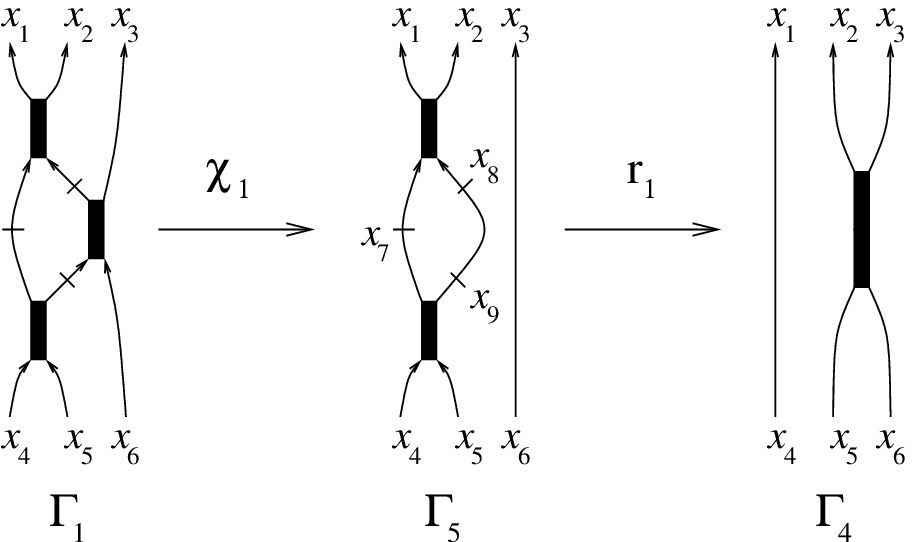}
 \caption{Map $g_1$} \label{mapgone}
 \end{figure}

Let 
$$g_1 = r_1 \chi_1, \hspace{0.2in} 
 C(\Gamma_1) \stackrel{\chi_1}{\lra} C(\Gamma_5) \stackrel{r_1}{\lra} 
 C(\Gamma_4)$$ 
be the composition depicted in figure~\ref{mapgone}.  

$r_1$ is the projection of $C(\Gamma_5)\cong C(\Gamma_4)\{-1\} \oplus 
 C(\Gamma_4)\{1\}$ onto the second direct summand (the direct sum 
decomposition is the one obtained in the proof of proposition~\ref{prop-tf}). 

We compute 
$$g_1g_0 = r_1 \chi_1 \chi_0 r_0= r_1 \hspace{0.03in} m(x_8-x_6) r_0 = 
  r_1 \hspace{0.03in} m(x_8) r_0,$$ 
since $r_1 \hspace{0.03in} m(x_6) r_0=0.$ The proof of 
proposition~\ref{prop-tf} implies that $r_1\hspace{0.03in} m(x_8) r_0$ 
is a nonzero multiple of the identity morphism. Therefore, $g_0$ is 
split injective, and  $C(\Gamma_4)$ is a direct 
summand of $C(\Gamma_1).$ \hspace{0.1in}  $\square$ 

\vspace{0.1in} 

Let $M$ be the complement of $C(\Gamma_4)$ in $C(\Gamma_1),$ 
   $$ C(\Gamma_1) \cong M \oplus C(\Gamma_4).$$

 \begin{lemma} $\gdim\hspace{0.05in} C(\Gamma_1) = 
 q^{-3}(1+q^2)(1+sq^{1-n})(1+sq^{3-n})^2$ if $n>2.$ \end{lemma} 

 \emph{Proof} is similar to that of lemma~\ref{lemma-gdim-0}.
 Up to a shift, $C(\Gamma_1)$ is the Koszul factorization 
 $\{\bba,\bbb\}$ for some length $6$ sequences whose entries 
 are polynomials in $x_1, \dots, x_6$ and in three other variables 
 (say $y_1,y_2,y_3$) assigned to the three internal marks of 
 $\Gamma_1.$ We can exclude two of these variables using 
 proposition~\ref{prop-exclude}, and then specialize to 
 $x_1=\dots=x_6=0,$ since we are computing cohomology. Details are 
 left to the reader. $\square$ 

 \vspace{0.1in}  

 In particular, if $n>2,$ the dimension of $C(\Gamma_1)$ is $16.$ 
 If $n=2,$ factorization $C(\Gamma_1)$ has dimension $8,$ 
 which is also the dimension of $C(\Gamma_4).$ Hence, $M=0$ if $n=2,$ 
 and there is an isomorphism $C(\Gamma_1)\cong C(\Gamma_4).$
 Factorization $\Upsilon$ is contractible ($v_3=3$ if $n=2$), and 
 propositions~\ref{core-isom}, \ref{isom-dir-4} follow. 
 From now on we assume $n>2.$ 

 It's easy to see that 
 \begin{eqnarray*}  
 \gdim \hspace{0.05in} C(\Gamma_4) & = & q^{-1} 
  (1+sq^{1-n})^2(1+sq^{3-n}), \\
 \gdim \hspace{0.05in} \Upsilon & = & q^{-3}(1+sq^{1-n})(1+sq^{3-n})
  (1+sq^{5-n}). 
 \end{eqnarray*}   
 Therefore, 
   $$ \gdim \hspace{0.05in} M =\gdim \hspace{0.05in} C(\Gamma_1) - 
   \gdim \hspace{0.05in} C(\Gamma_4)=\gdim \hspace{0.05in} \Upsilon.$$
  Thus, $M$ and $\Upsilon$ have the same graded dimension. 

 \begin{lemma}  
 $$ \gdim \hspace{0.05in} \Ext(\Upsilon,C(\Gamma_4)) =
    q^{3n-3}[2][n][n-1][n-2].$$ 
 \end{lemma} 

 \emph{Proof:}  The ext groups in the lemma are isomorphic to 
 the cohomology of the 2-complex 
$$C(\Gamma_4)\otimes \Upsilon_{\bullet}\cong \{ \bba,\bbb\} \bracket{1}\{3n-9\},$$ 
where 
$$ (\bba,\bbb)= \left( \begin{array}{cc} \pi_{36} & x_3 -x_6 \\
                    u_1(x_1,x_2,x_4,x_5) & x_1+x_2 -x_4 -x_5   \\
                    u_2(x_1,x_2,x_4,x_5) & x_1 x_2 - x_4 x_5  \\
                      \alpha_1   &   v_1 \\
                     \alpha_2   & v_2 \\
                     \alpha_3 & v_3 \end{array} \right)  $$
Pass to the subring $R''$ of $R'=\Q[x_1, \dots, x_6]$ generated by 
$$x_1'=x_1+x_2, \hspace{0.1in} x_2'=x_1x_2, \hspace{0.1in} x_3, 
  \hspace{0.1in} x_4'=x_4 + x_5, \hspace{0.1in} x_5'=x_4 x_5, \hspace{0.1in} x_6,$$
and denote by $\{\bba,\bbb\}_0$ the cyclic Koszul complex of the pair $(\bba,\bbb)$ viewed 
as a pair in $R''$ (this is possible since all coefficients of $\bba$ and $\bbb$ lie in $R''$). 
Then  
 $$ \gdim \hspace{0.02in}  H(\{\bba,\bbb\}) = q^2 [2]^2 \gdim \hspace{0.02in} 
  H(\{\bba,\bbb\}_0).$$ 
 We use the first three lines of $(\bba,\bbb)_0$ to exclude $x_6,x_4',$ and $x_5',$
 and denote by $(\bba^1,\bbb^1)$ the pair obtained by crossing out the first three 
lines of $(\bba,\bbb)_0$ and passing to the quotient ring of $R''$ by relations 
 $x_6=x_4'=x_5'=0.$ Then $\bba^1=(0,0,0),$  $\bbb^1=(-v_1,-v_2,-v_3),$  
$$ H(\{\bba,\bbb\}_0)\cong H(\{\bba^1,\bbb^1 \}) \cong H((0,0,0),(v_1,v_2,v_3)),$$ 
and the lemma follows easily, since $(v_1,v_2,v_3)$ is a regular sequence in 
this quotient ring $R'''=R''/(x_6,x_4',x_5'),$ and 
 $R'''/(v_1,v_2,v_3)$  is isomorphic to the cohomology ring of the partial flag variety 
  $$ \{ N_2\subset N_3 | \dim(N_i)=i, \hspace{0.1in} N_3\subset \C^n\}.$$ 
 
$\square$ 
 
 \begin{lemma}  
 $$ \gdim \hspace{0.05in} \Ext(\Upsilon,C(\Gamma_1)) =
    q^{3n-3}[2]^3[n][n-1][n-2].$$
 \end{lemma} 
 
\emph{Proof:} The ext groups in the lemma are isomorphic to 
 the cohomology of the 2-complex $C(\Gamma_1)\otimes \Upsilon_{\bullet}.$ 
This is a complex of free $R^1$-modules, where $R^1$ is the ring of polynomials 
in $x_1,\dots, x_6$ and three "internal" variables corresponding to the three marks in 
$\Gamma_1.$ Passing to a suitable subring (over which $R^1$ is a free rank $4$ module) 
and excluding several variables, one 
can show that $C(\Gamma_1)\otimes \Upsilon_{\bullet}$ is isomorphic (up to contractible 
complexes) to the 
direct sum of four copies of $C(\Gamma_1)\otimes \Upsilon_{\bullet}$ with shifts, 
implying 
   $$ \gdim \hspace{0.05in} \Ext(\Upsilon,C(\Gamma_1)) = [2]^2 \hspace{0.02in}
    \gdim \hspace{0.05in} \Ext(\Upsilon,C(\Gamma_4)).$$ 

 $\square$ 

The lemmas imply 
  \begin{eqnarray*} 
  \gdim \hspace{0.05in} \Ext(\Upsilon,M) & = & 
    \gdim \hspace{0.05in} \Ext(\Upsilon,C(\Gamma_1))  -  
   \gdim \hspace{0.05in} \Ext(\Upsilon,C(\Gamma_4))  \\ 
   & = & q^{3n-3}[2][3][n][n-1][n-2]. 
  \end{eqnarray*} 
 Since $q^{3n-3}[2][3][n][n-1][n-2]\in 1 + q\Z[q],$ the space of 
  degree $0$ morphisms from $\Upsilon$ to $M$ is one-dimensional. 
  Let $\theta$ be a nontrivial (not null-homotopic) 
 degree $0$ morphism from $\Upsilon$ to $M.$

 \begin{lemma} \label{theta-isom} 
  $\theta$ is an isomorphism of factorizations. \end{lemma} 

 \emph{Proof:} Factorization $\Upsilon$ has the form 

\begin{equation*} 
 \left( \begin{array}{c} R'\{-3\} \\ R'\{5-2n\} \\ R'\{3-2n\} \\
  R'\{1-2n\} \end{array} \right) \xrightarrow{\Lambda_0} 
  \left( \begin{array}{c} R'\{2-n\} \\ R'\{-n\} \\ R'\{-2-n\} \\
  R'\{6-3n\} \end{array} \right) \xrightarrow{\Lambda_1}  
 \left( \begin{array}{c} R'\{-3\} \\ R'\{5-2n\} \\ R'\{3-2n\} \\
  R'\{1-2n\} \end{array} \right), 
\end{equation*}
where 
\begin{equation*} 
  \Lambda_0 = \left( \begin{array}{cccc} v_3 & -\alpha_2 & -\alpha_1 & 
  0   \\  v_2 & \alpha_3 & 0 & -\alpha_1  \\
  v_1 & 0 & \alpha_3 & \alpha_2 \\ 0 & v_1 & -v_2 & v_3 
   \end{array} \right), \hspace{0.2in} 
  \Lambda_1 = \left( \begin{array}{cccc} \alpha_3 & \alpha_2 & 
  \alpha_1 & 0 \\ -v_2 & v_3 & 0 & \alpha_1 \\
  -v_1 & 0 & v_3 & -\alpha_2 \\ 0 & -v_1 & v_2 & \alpha_3 
  \end{array} \right). 
\end{equation*} 

We can assume that the factorization $M=M^0\oplus M^1$ has no 
contractible summands. Since it has the same graded dimension 
as $\Upsilon,$ we can write $M$ as  
\begin{equation} \label{m-shifts} 
 \left( \begin{array}{c} R'\{-3\} \\ R'\{5-2n\} \\ R'\{3-2n\} \\
  R'\{1-2n\} \end{array} \right) \xrightarrow{\Pi_0} 
  \left( \begin{array}{c} R'\{2-n\} \\ R'\{-n\} \\ R'\{-2-n\} \\
  R'\{6-3n\} \end{array} \right) \xrightarrow{\Pi_1}  
 \left( \begin{array}{c} R'\{-3\} \\ R'\{5-2n\} \\ R'\{3-2n\} \\
  R'\{1-2n\} \end{array} \right)
\end{equation} 
for some $4\times 4$ matrices $\Pi_0,\Pi_1$ with homogeneous 
entries whose degrees match degrees of the entries of $\Lambda_0,$
 $\Lambda_1.$ Let $\Pi_0 = (a_{ij})_{i,j=1}^4$ 
and $\Pi_1=(b_{ij})_{i,j=1}^4.$ The entry 
$a_{14}$ is a degree $n+1$ homomorphism from  $R'\{1-2n\}$ to 
$R'\{2-n\},$ and, necessarily a rational number. If $a_{14}\not= 0,$ 
then $M$ contains a contractible summand, which is a contradiction. 
Therefore, $a_{14}=0.$ Likewise, $b_{14}=0$ if $n>4.$ 

 We can write morphism $\theta$ via a 
pair of matrices 
$$\Theta_0 = (c_{ij})_{i,j=1}^4, \hspace{0.2in} 
  \Theta_1= (f_{ij})_{i,j=1}^4.$$ 
Since $\theta$ is a morphism of factorizations, the diagram 
  \begin{equation*}
   \begin{CD}
     \Upsilon^0 @>{\Lambda_0}>> \Upsilon^1 @>{\Lambda_1}>> \Upsilon^0   \\
     @V{\Theta_0}VV      @V{\Theta_1}VV   @V{\Theta_0}VV \\
     M^0 @>{\Pi_0}>> M^1 @>{\Pi_1}>> M^0
   \end{CD}
 \end{equation*}
is commutative, and the following relations hold 
\begin{eqnarray} \label{eq-m-1} 
  \Pi_0 \Theta_0 & = & \Theta_1 \Lambda_0,   \\ 
  \Pi_1 \Theta_1 & = & \Theta_0 \Lambda_1. \label{eq-m-2} 
\end{eqnarray} 
$\theta$ is grading-preserving. If $n>4,$ sequences of degree 
shifts $(-3,5-2n,3-2n,1-2n)$ and $(2-n,-n,-2-n,6-3n)$ as 
strictly decreasing, and matrices $\Theta_0,\Theta_1$ 
are lower-triangular, $c_{ij}=f_{ij}=0$ for $i<j,$
since there are no grading-preserving homomorphisms from $R'\{k_1\}$ 
to $R'\{k_2\}$ if $k_2>k_1.$ If $n=4,$ the sequences are not 
strictly decreasing, but changing bases in $M^0$ and $M^1,$ if 
necessary, we can assume that $\Theta_0$ and $\Theta_1$ are 
lower-triangular in this case too. When $n=3,$ we cannot assume right 
away that $\Theta_0$ and $\Theta_1$ are lower-triangular. We 
postpone considering cases $n=3$ and $n=4$ until later, and from here on 
until two paragraphs after the proof of lemma~\ref{f33is0} restrict 
to the case $n>4.$ 

$\theta$ is an isomorphism iff both $\Theta_0$ and $\Theta_1$ 
are invertible over the ring $R'.$ Since they are lower-triangular 
with rational diagonal entries, they are invertible iff all 
diagonal entries are non-zero. Moreover, $\Theta_1$ is invertible 
iff $\Theta_0$ is invertible. This is due to the equation 
  $\Pi_0\Theta_0 \Lambda_0^{-1}= \Theta_1$ and the property 
that determinants of $\Pi_0$ and $\Lambda_0$ are nonzero multiples 
 of $w^2$ (to verify this property, use that $w$ is an 
 irreducible polynomial, $\Lambda_0 \Lambda_1 = \Pi_0 \Pi_1=w \mathrm{I},$
 and $\Pi_0$ has coefficients in the same degrees as $\Lambda_0$). 

Thus, $\theta$ is an isomorphism iff $\Theta_0$ is invertible. 

\vspace{0.1in} 

Matrix equation (\ref{eq-m-1}) 
can be rewritten as $16$ equations on the entries of matrices 
$\Pi_0,\Theta_0,\Theta_1.$ Three of these equations have the form 
 \begin{eqnarray*}
   a_{24}c_{44} & = & - \alpha_1 f_{22}, \\
   a_{34}c_{44} & = & -\alpha_1 f_{32} +\alpha_2f_{33}, \\
   a_{44}c_{44} & = & -\alpha_1 f_{42} + \alpha_2 f_{43} + v_3 f_{44}. 
 \end{eqnarray*}

\begin{lemma} \label{c44not0} $c_{44}\not= 0.$ \end{lemma} 

To prove the lemma, assume to the contrary. Then $c_{44}=0$ and 
the above equations reduce to 
  \begin{eqnarray*}
    0 & = & - \alpha_1 f_{22} \\
    0 & = & -\alpha_1 f_{32} +\alpha_2f_{33} \\
    0 & = & -\alpha_1 f_{42} + \alpha_2 f_{43} + v_3 f_{44}. 
 \end{eqnarray*}
Since $f_{22},f_{33},f_{44}\in \Q,$ and the sequence 
$(\alpha_1,\alpha_2,v_3)$ is regular, the equations imply 
  $$ f_{22}=f_{32}=f_{33}=f_{44}=0.$$ 
The last equation among the three above simplifies to 
  $$ f_{42}\alpha_1 = f_{43}\alpha_2.$$ 
 Since $\alpha_1 $ and $\alpha_2$ are relatively prime, 
  $$f_{42}=\alpha_2 z_0, \hspace{0.3in} f_{43}=\alpha_1 z_0,$$ 
 for some polynomial $z_0.$ 

 Three of the equations in (\ref{eq-m-2}) are 
  \begin{eqnarray*} 
    b_{24}f_{44} & = & \alpha_1 c_{22}, \\
    b_{34}f_{44} & = & \alpha_1 c_{32}-\alpha_2 c_{33}, \\
    b_{44}f_{44} & = & \alpha_1 c_{42}-\alpha_2 c_{43}+\alpha_3 c_{44}. 
  \end{eqnarray*} 
Since $f_{44}=0,$ and $(\alpha_1,\alpha_2,\alpha_3)$ is a 
 regular sequence, we derive 
  $$ c_{22}=c_{32}=c_{33}=c_{44}=0,$$
 and 
 $$ c_{42}=\alpha_2 z_1, \hspace{0.3in} c_{43}=\alpha_1 z_1,$$ 
 for some polynomial $z_1.$ 
 
 The remaining terms of matrix equations (\ref{eq-m-1}), (\ref{eq-m-2}) 
 imply  
 \begin{eqnarray*} 
 f_{21}=-a_{24}z_1, & f_{31}=-a_{34}z_1, & f_{41}=\alpha_3 z_0 - 
   a_{44}z_1, \\
 c_{21}= b_{24}z_0, & c_{31}=b_{34}z_0, & c_{41}=b_{44}z_0-v_3 z_1. 
 \end{eqnarray*} 
 Let $Z_0,Z_1$ be $4\times 4$-matrices with the only nonzero 
entry $z_0,$ (respectively, $z_1$), placed on the intersection 
of the fourth row and the first column. Then 
 $$ \Theta_0 =\Pi_1 Z_0 - Z_1 \Lambda_0,\hspace{0.2in} 
    \Theta_1 =Z_0 \Lambda_1 - \Pi_0 Z_1. $$ 
We see that $\theta$ is homotopic to $0,$ through the homotopy 
$(Z_0,-Z_1).$ This contradicts our assumption on $\theta.$
Therefore, $c_{44}\not=0,$ and lemma~\ref{c44not0} is true. $\square$  
  
\begin{lemma} \label{f44not0} $f_{44}\not=0.$ \end{lemma} 

\emph{Proof:} An argument in the proof of the previous lemma shows 
 that $f_{44}=0$ implies $c_{44}=0.$ \hspace{0.1in} $\square$ 

\vspace{0.1in} 

Since $f_{44}\not=0 $ and $c_{44}\not= 0,$ 
by rescaling the last basis vector in $M^0$ and in $M^1,$ we can 
assume 
 $$ f_{44}=c_{44}=1.$$ 

To finish the proof of lemma~\ref{theta-isom} in the case $n>3,$ we 
assume to the contrary: $\theta$ is not an isomorphism.  

\begin{lemma}\label{f11is0} $f_{11}=0$ if $\theta$ is not 
an isomorphism. \end{lemma}

\emph{Proof:} Assume $f_{11}\not= 0.$ By changing the first basis vector in 
$M^1,$ if necessary, we can reduce to the case 
 $$f_{11}=1, \hspace{0.2in} f_{21}=f_{31}=f_{41}=0. $$ 
The entry $(1,3)$ of equation (\ref{eq-m-1}) is 
 $$ a_{13}c_{33}= - \alpha_1 f_{11}= - \alpha_1.$$ 
Therefore $c_{33}\not= 0,$ and we assume, without loss of 
generality (by changing the third basis vector in $M^0$), that 
 $$c_{33}=1, \hspace{0.2in} c_{43}=0.$$
Then $a_{13}=-\alpha_1,$ and equation (\ref{eq-m-1}), entry 
$(1,2),$ simplifies to 
$$-\alpha_2 = a_{12}c_{22} - \alpha_1 c_{32}.$$
Since $\alpha_2$ is not divisible by $\alpha_1,$ we know that 
$c_{22}\not=0.$ Changing the second basis vector in $M^0,$ we 
assume 
 $$c_{22}=1, \hspace{0.2in} c_{32}=0, \hspace{0.2in} 
   c_{42}=0, \hspace{0.2in} a_{12}=-\alpha_2.$$
Since $c_{22}=c_{33}=c_{44}=1$ and $\Theta_0$ is not invertible 
(since we assumed that $\theta$ is not an isomorphism), necessarily 
 $c_{11}=0.$ 
 
The entry $(1,1)$ of (\ref{eq-m-1}) can be written as 
 $$v_3 = -\alpha_2 c_{21} - \alpha_1 c_{31}.$$ 
Contradiction, since $v_3$ is not in the ideal 
 $(\alpha_1,\alpha_2).$ Lemma~\ref{f11is0} follows. $\square$

\begin{lemma} \label{c11is0} $c_{11}=0$ if $\theta$ is 
not an isomorphism \end{lemma} 

Assume otherwise ($c_{11}\not= 0$) and change the first basis 
vector of $M^0$ so that 
 $$c_{11}=1, \hspace{0.2in} c_{21}=c_{31}=c_{41}=0.$$ 
From the previous lemma we know that $f_{11}=0.$ 
Entry $(1,3)$ of the equation (\ref{eq-m-2}) reads 
$\alpha_1 = b_{13}f_{33}.$ We see that $f_{33}\not=0,$ and 
after modifying the third basis vector of $M^1$ we assume 
 $$f_{33}=1, \hspace{0.2in} f_{43}=0, \hspace{0.2in} b_{13}=\alpha_1.$$ 
Entry $(1,2)$ of the equation (\ref{eq-m-2}) reads 
$$ \alpha_2 = b_{12}f_{22} + \alpha_1 f_{32}.$$ 
Since $\alpha_2$ does not factor, $f_{22}\not=0,$ and 
we can reduce to the case 
 $$f_{22}=1, \hspace{0.2in} f_{32}=f_{42}=0, \hspace{0.2in} 
 b_{12}=\alpha_2.$$ 
Entry $(1,1)$ of equation (\ref{eq-m-2}) simplifies to 
 $$ \alpha_3 = \alpha_2 f_{21} + \alpha_1 f_{31}.$$ 
Contradiction, since $\alpha_3$ is not in the ideal generated by 
$\alpha_1$ and $\alpha_2.$ Lemma~\ref{c11is0} follows. $\square$ 

\vspace{0.07in}  

By now, we have reduced our consideration to the case 
 $$f_{44}=c_{44}=1, \hspace{0.2in} f_{11}=c_{11}=0.$$ 

\begin{lemma} $f_{33}=0$ if $\theta$ is not an isomorphism. 
  \label{f33is0} \end{lemma} 

\emph{Proof:} if $f_{33}\not=0,$ we can assume $f_{33}=1,$ 
  $f_{43}=0.$ Entry $(1,3)$ of the equation (\ref{eq-m-2}) becomes 
$b_{13}=0.$ 
Note that $b_{12}\not=0;$ otherwise $b_{11}$ (which is a 
polynomial of degree $3$ in $x_i$'s) would be 
the only nonzero entry in the first row of $\Pi_1,$ and a divisor 
of $\mathrm{det}\hspace{0.03in} \Pi_1,$ which is proportional to 
$w^2.$ Since $w$ is irreducible, this is impossible, so that 
 $b_{12}\not= 0.$ 

Entry $(1,2)$ of the equation (\ref{eq-m-2}) simplifies to 
 $b_{12}f_{22}=0.$ Thus, $f_{22}=0.$ Entry $(1,1)$ of the 
same equation reduces to $b_{12}f_{21}=0.$ Thus, $f_{21}=0.$ 

Next, entries $(2,3)$ and $(2,4)$ of (\ref{eq-m-2}) reduce to 
 $$ b_{24}=\alpha_1 c_{22}, \hspace{0.2in} b_{23}=\alpha_1 c_{21}.$$ 
Entry $(2,2)$ simplifies to 
 $$b_{23}f_{32} + b_{24}f_{42} = \alpha_2 c_{21} + v_3 c_{22}.$$
Since $v_3$ is not in the ideal generated by $\alpha_1,\alpha_2,$
we know that $c_{22}=0.$ Then $b_{24}=\alpha_1 c_{22}=0 $ 
(entry $(2,4)$ of (\ref{eq-m-2})). Entry $(2,3)$ tells us that 
$b_{23}= \alpha_1 c_{21},$ while entry $(2,2)$ reduces to 
 $\alpha_2 c_{21} = \alpha_1 c_{21} f_{32}.$ Since $\alpha_2$ 
is not divisible by $\alpha_1,$ we have $c_{21}=0$ and $b_{23}=0.$ 

Switching to the equation (\ref{eq-m-1}), and looking at the 
entry $(4,4),$ we get $a_{44}=v_3 - \alpha_1 f_{42},$ while entry 
 $(4,3)$ simplifies to 
 $$ -v_2 - \alpha_1 f_{41} = a_{43} c_{33} + 
 (v_3-\alpha_1 f_{42})c_{43}. $$ 
Since $v_2$ does not lie in the ideal $(\alpha_1,v_3),$ we have 
$c_{33}\not= 0.$ Then we can assume $c_{33}=1, c_{43}=0.$ 
Entries $(1,3)$ and $(2,3)$ tells us that $a_{13}=a_{23}=0.$ 

To summarize, we've shown that 
 $$a_{13}=a_{23}=a_{24}=0.$$ 
Also, $a_{14}=0.$ Therefore, the matrix $\Pi_0$ is block lower-diagonal. 
It's determinant is a non-zero multiple of $w^2.$ The determinant of 
 $\Pi_0$ is 
divisible by the determinant of its lower diagonal $2\times 2$ block 
 $$ \left( \begin{array}{cc} a_{33} & a_{34} \\ a_{43} & a_{44} 
   \end{array} \right) = 
  \left( \begin{array}{cc} \alpha_3-\alpha_1 f_{31} & 
  \alpha_2-\alpha_1 f_{32} \\ -v_2-\alpha_1 f_{41} & v_3 - \alpha_1 
  f_{42} 
   \end{array} \right)$$
The only possibility is for this determinant to be a nonzero multiple 
 of $w,$ leading to the equation 
 $$ ( \alpha_3-\alpha_1 f_{31})(  v_3 - \alpha_1 f_{42} ) + 
  (\alpha_2 -\alpha_1 f_{32})(v_2+\alpha_1 f_{41}) = \lambda w,$$ 
for $\lambda \in \Q.$ Using that $w=\alpha_1 v_1 + 
  \alpha_2 v_2 + \alpha_3 v_3,$ we reduce the equation to 
 $$\alpha_1(-v_1-\alpha_3 f_{42} - v_3 f_{31} + \alpha_1 f_{31}f_{42} 
 +\alpha_2 f_{41} -v_2 f_{32} - \alpha_1 f_{32}f_{41})= (\lambda-1)w.$$ 
Since $w$ is not divisible by $\alpha_1,$ we have $\lambda=1,$ 
and the equation becomes
 $$-v_1-\alpha_3 f_{42} - v_3 f_{31} + \alpha_1 f_{31}f_{42} 
 +\alpha_2 f_{41} -v_2 f_{32} - \alpha_1 f_{32}f_{41}=0.$$ 
Contradiction, since $v_1$ does not belong to the ideal generated by 
$\alpha_1,$ $\alpha_2,$ $\alpha_3,$ $v_2,$ $v_3.$ Lemma~\ref{f33is0}  
follows. $\square$ 

\vspace{0.1in} 

Entries $(4,3)$ and $(4,4)$ of the equation (\ref{eq-m-2}) 
now become 
 \begin{eqnarray*} 
  b_{44} f_{43} & = & \alpha_1 c_{41} + v_3 c_{43} + v_2, \\
  b_{44}  & = & \alpha_1 c_{42} - \alpha_2 c_{43} + \alpha_3, 
 \end{eqnarray*} 
implying 
 $$ v_2 = f_{43}( \alpha_1 c_{42} - \alpha_2 c_{43} + \alpha_3) 
      - \alpha_1 c_{41} - v_3 c_{43}.$$ 
This is impossible, since $v_2$ does not lie in 
 the ideal $(\alpha_1,\alpha_2,\alpha_3, v_3).$ 

Therefore,  $\theta$ is invertible and $\Upsilon$ is isomorphic 
to $M$ if $ n>4.$ Lemma~\ref{theta-isom} and 
propositions~\ref{core-isom}, \ref{isom-dir-4} 
follow in the case $n>4.$          
 
\vspace{0.2in} 

If $n=4,$ the sequences of degree 
shifts in $M^0,M^1$ (see formula (\ref{m-shifts}) are 
$(-3,-3,-5,-7)$ and $(-2,-4,-6,-6).$ These sequences are decreasing, 
and matrices $\Theta_0,\Theta_1$ 
are block lower-triangular. By changing bases in $M^0$ and $M^1,$ if 
necessary, we can assume that $\Theta_0$ and $\Theta_1$ are 
lower-triangular. The entry $b_{14}$ of $\Pi_1$ is a homogeneous 
linear polynomial in $x$'s. If $b_{14}=0,$ 
lemmas~\ref{c44not0}-\ref{f33is0} hold for $n=4$ as well (with 
the simplification that $z_0=0$ in the proof of lemma~\ref{c44not0}), 
and we are done. Assume now that $b_{14}\not= 0.$ Then 

\begin{itemize} 
 \item $f_{44}=0$ (from equation (\ref{eq-m-2}), entry $(1,4)$), 
 \item $c_{22}=0$ (equation (\ref{eq-m-2}), entry $(2,4)$), 
  \item $c_{33}=0$ (equation (\ref{eq-m-2}), entry $(3,4)$), 
  \item $c_{32}=0$ (equation (\ref{eq-m-2}), entry $(3,4)$), 
  \item $c_{44}=0$ (equation (\ref{eq-m-2}), entry $(4,4)$), 
  \item $f_{22}=0$ (equation (\ref{eq-m-1}), entry $(2,4)$), 
   \item $f_{11}=0$ (equation (\ref{eq-m-1}), entry $(1,3)$), 
   \item $f_{33}=0$ (equation (\ref{eq-m-1}), entry $(3,4)$), 
   \item $f_{32}=0$ (equation (\ref{eq-m-1}), entry $(3,4)$). 
\end{itemize} 
 Furthermore, $b_{14}f_{43}= \alpha_1 c_{11}$ (equation (\ref{eq-m-2}), 
 entry (1,3)), and $b_{14}f_{42}=\alpha_2 c_{11}$ (equation 
 (\ref{eq-m-2}), entry $(1,2)$). Since $\alpha_1$ and $\alpha_2$ 
 are relatively prime, this is only possible if $c_{11}=0,$ 
which, in turn, implies $f_{43}=0, f_{42}=0.$ 
Next, $c_{21}=0$ (equation (\ref{eq-m-2}), entry $(2,2)$), and 
 $c_{31}=0$ (equation (\ref{eq-m-2}), entry $(3,2)$). 
From the remaining equations we derive 
\begin{eqnarray*} 
 c_{41} = z v_3, & c_{42}=-z \alpha_2, & c_{43}=-z \alpha_1, \\
 f_{21} = z a_{24}, & f_{31}=z a_{34}, & f_{41}= z a_{44}, 
\end{eqnarray*} 
for some rational number $z.$ Therefore, $\theta$ is null-homotopic, 
through the homotopy $(0,Z),$ where $Z$ is the $4\times 4$ matrix 
with the only non-zero entry $z$ in the lower left corner. 
Contradiction. Lemma~\ref{theta-isom} and 
propositions~\ref{core-isom}, \ref{isom-dir-4} follow 
in the case $n=4.$ 
 
\vspace{0.2in} 

If $n=3,$ degree conditions force coefficients 
$$a_{14}, b_{23}, c_{21}, c_{14}, c_{23},c_{24},c_{34}, 
  f_{12},f_{13},f_{14},f_{23}, f_{43}$$ 
to be zero. By changing bases in $M^0,M^1,$ if necessary,  
we can assume that $c_{13}=f_{24}=0.$ 

\begin{lemma} \label{c44not0n3} $c_{44}\not= 0$ if $n=3.$ 
\end{lemma} 

\emph{Proof:} Assume to the contrary, $c_{44}=0.$ Three of the 
equations (\ref{eq-m-1}) reduce to 
 \begin{eqnarray*} 
    0 & = & - \alpha_1 f_{22}, \\
    0 & = & -\alpha_1 f_{42}+v_3 f_{44}, \\
    0 & = & -\alpha_1 f_{32} + \alpha_2 f_{33} + v_3 f_{34}, 
 \end{eqnarray*} 
implying that 
 $$f_{22}=f_{44}=f_{42}=f_{33}=0.$$ 
 Equation (\ref{eq-m-2}), entry $(2,4),$ implies $c_{22}=0,$ 
 while the entry $(3,3)$ implies $c_{33}=c_{31}=0.$ 
 From the entry $(1,3)$ we derive $c_{11}=0,$ and from 
the entry $(1,3)$ of (\ref{eq-m-1}) that $f_{11}=0.$ 
 The remaining entries of the two matrix equations force 
 the rest of entries of $\Theta_0,\Theta_1$ to have the 
 form 
 \begin{eqnarray*}
  & & c_{12}=z_0 b_{13}, \hspace{0.2in} 
      c_{32}=z_0 b_{33}, \hspace{0.2in}  \\
  & & c_{41}= - v_3 z_1, \hspace{0.2in} 
      c_{42}= b_{43}z_0 + \alpha_2 z_1, \hspace{0.2in} 
      c_{43}= z_1 \alpha_1, \\
  & & f_{21}= -z_1 a_{24}, \hspace{0.2in} 
      f_{31}=-z_1 a_{34} - xv_2, \hspace{0.2in} 
      f_{41} = - z_1 a_{44}, \\
  & & f_{32}= z_0 v_3, \hspace{0.2in} f_{34}=z_0 \alpha_1. 
 \end{eqnarray*} 
 Therefore, we can write 
 $$ \Theta_0 = Z_1 \Lambda_0 + \Pi_1 Z_0, \hspace{0.2in} 
    \Theta_1 = \Pi_0 Z_1 + Z_0 \Lambda_1,$$ 
for matrices $Z_0,Z_1,$ which have zero entries save for 
$z_0$ in $(3,2)$ in $Z_0$ and $-z_1$   in $(4,1)$ in $Z_1,$ 
with $z_0,z_1\in \Q.$ Thus, $\theta$ is homotopic to $0$ 
 and the lemma follows. $\square$ 

\vspace{0.1in} 

The proof of the lemma mimicked very closely the proof 
 of lemma~\ref{c44not0}. Proofs of the following lemmas 
are omitted, since they are parallel to 
those of lemmas~\ref{f44not0}-\ref{f33is0}. 

\begin{lemma} $f_{33}\not= 0$ if $n=3.$ \end{lemma} 

\begin{lemma} $f_{11}=0$ if $n=3$ and $\theta$ is not an 
isomorphism. \end{lemma} 

\begin{lemma} $c_{22}=0$ if $n=3$ and $\theta$ is not an 
isomorphism. \end{lemma} 

\begin{lemma} $f_{44}=0$ if $n=3$ and $\theta$ is not an 
isomorphism. \end{lemma} 

We can assume that $f_{33}=c_{44}=1.$ 
Entries $(4,3)$ and $(4,4)$ of the equation (\ref{eq-m-2}) 
now become
\begin{eqnarray*} 
  b_{43} & = & \alpha_1 c_{41} + v_3 c_{43} + v_2 , \\
 b_{43} f_{34} & = & \alpha_1 c_{42} - \alpha_2 c_{43}+\alpha_3, 
\end{eqnarray*} 
leading to a contradiction. 
  Lemma~\ref{theta-isom} and 
propositions~\ref{core-isom}, \ref{isom-dir-4} follow 
in the remaining case $n=3.$ 

$\square$


\vspace{0.2in} 

{\bf Shifts.} 
Direct sum decompositions in the above propositions often contain 
terms shifted by $\bracket{1}.$ In general, if there is a direct 
sum decomposition that contains $C(\Gamma_1)$ 
and $C(\Gamma_2)$, for a pair of graphs $\Gamma_1,\Gamma_2,$ 
to determine whether $C(\Gamma_2)$ will be shifted by $\bracket{1}$ 
relative to $C(\Gamma_1)$ modify $\Gamma_1,\Gamma_2$ by substituting 
a pair of arcs for each wide edge of these graphs (right to left 
transformation in figure~\ref{pair2}). The graphs will turn into 
collections of arcs and circles inside a disc. The two collections 
share the set of boundary points. Glue them together along their 
boundaries to get a collection of circles on a 2-sphere. If the number 
of circles plus half the number of boundary points is odd, there is a 
shift by $\bracket{1}.$ Otherwise, $C(\Gamma_1)$ and $C(\Gamma_2)$ are 
unshifted relative to each other. 

\vspace{0.1in} 

{\bf Closed graphs.} We say that $\Gamma$ is closed if it has no 
boundary points. $C(\Gamma)$ is a 2-complex ($d^2=0$) iff 
$\Gamma$ is closed. In this subsection we assume that $\Gamma$ is 
closed. Define the parity of $\Gamma$ (denoted $p(\Gamma)$) 
as the number of circles in the modification of $\Gamma$ 
described above (erase a neighbourhood of each wide edge, 
and add two oriented arcs in place of this edge). 
Let $H(\Gamma)$ be the cohomology groups of $C(\Gamma).$ 
Since $C(\Gamma)$ has cohomology only in degree $p(\Gamma),$ 
we have 
$$H(\Gamma) = H^{p(\Gamma)}(C(\Gamma)).$$ 
Internal grading of $C(\Gamma)$ makes $H(\Gamma)$ into 
a $\Z$-graded $\Q$-vector space. 

\begin{prop} The graded dimension of $H(\Gamma)$ is 
the invariant $P_n(\Gamma),$ 
 $$ \sum_{j\in \Z} \dim H^j(\Gamma)\hspace{0.03in} 
   q^j = P_n(\Gamma).$$ 
\end{prop} 

Proposition follows from the direct sum decompositions 
obtained in this section and the well-known fact that 
skein relations in figure~\ref{graphrel} suffice to evaluate 
$P_n(\Gamma)$ for any graph $\Gamma.$ 

\vspace{0.1in} 

{\bf Examples of homologically regular pairs.}
We continue to assume that $\Gamma$ is closed, and, hence, 
 $C(\Gamma)$ is a 2-complex. It's cohomology $H(\Gamma)$ 
are concentrated 
in one degree only. $C(\Gamma)$ is the cyclic Koszul complex 
$\{\bba,\bbb\}$ for a suitable pair $(\bba,\bbb)$ with $\bba\bbb=0.$ If 
$H^0(\Gamma)\not= 0,$ then $H^1(\Gamma)=0,$ and the pair 
$(\bba,\bbb)$ is homologically regular, in the terminology 
of section~\ref{sec-koszul}. If $H^1(\Gamma)\not=0,$ we can permute 
$a_i$ with $b_i$ in the sequences $\bba,\bbb,$ for some $i,$ 
to produce a homologically regular pair. Thus, any closed 
graph gives rise to a homologically regular pair. 

The simplest example of a homologically regular pair is 
$({\bf 0},\bbb),$ where $\bbb$ is a regular sequence. 

For any pair, $H(\{ \bba,\bbb\})$ is a module over the 
ring $R.$ Note that cohomology $H^0(\{{\bf 0},\bbb\}),$
if $\bbb$ is regular, is a cyclic module over $R,$ isomorphic 
to $R/(\bbb).$ Modifications of this pair $({\bf 0},\bbb)$ 
using symmetries from the group $G$ in section~\ref{sec-koszul} 
do not change the cohomology and its module 
structure. All examples of homologically regular pairs given in 
section~\ref{sec-koszul} have the property that 
$H^0$ is a cyclic $R$-module. 

 \begin{figure} [htb] \drawing{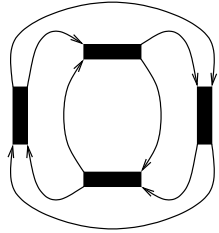}
   \caption{Graph $\Gamma$} \label{hregular}
 \end{figure}

Interestingly, homologically regular pairs $(\bba,\bbb)$ 
assigned to certain graphs $\Gamma$ have the property that 
$H^0$ is not a cyclic $R$-module. The simplest example of 
such $\Gamma$ is shown in figure~\ref{hregular}. We assume 
that $n>2,$ and each oriented edge has one mark (not shown 
on the diagram). The graded dimension of $H^0(C(\Gamma)\bracket{1})$ is 
$$ [n][n-1](2[n-1]+[n-3])$$ 
and has the form $q^{5-3n}(2+x),$ where $x\in q\Z[q].$ 
Since the ring $R$ of polynomials in edge variables 
is non-negatively graded, with $\Q$ in degree $0,$ the 
cohomology of $C(\Gamma)\bracket{1}$ is not a cyclic module 
over $R$ (it's $\Q\oplus \Q$ in the lowest degree).  
Hence, the pair $(\bba,\bbb)$ assigned to this 
graph cannot be obtained from any pair of the form 
$({\bf 0}, \bbb')$ using symmetries from the group $G$ (see the end of 
section~\ref{sec-koszul}.  
Thus, the pair assigned to $\Gamma$ is homologically regular 
and "intrinsically cyclic", unlike pairs $\{ {\bf 0}, \bbb'\}$ 
whose cyclic Koszul complexes are just the Koszul complexes 
of $\bbb'$ with collapsed grading.


 \section{Tangle diagrams and complexes of factorizations} 
 \label{sec-tangle} 

{\bf Complexes of factorizations.} 
To any additive category $\cC$ associate the category  
$K(\cC)$ with objects--bounded complexes of objects of 
$\cC$ and morphisms--morphisms of complexes up to homotopies. 
The category $K(\cC)$ is triangulated. 

Recall that $\hmf_w$ is the category of graded 
factorizations with potential $w$ and finite-dimensional 
cohomology, up to homotopies. 
Let $K_w= K(\hmf_w),$ the homotopy category of 
$\hmf_w.$ An object of $K_w$ is $\Z\oplus \Z \oplus \Z_2$-graded. 

We distinguish the three shift functors in $K_w,$ 
the shifts functors $\langle 1\rangle, \{1\}$ coming from 
$\hmf_w$ and the shift $[1]$ in the category of complexes. 
These three functors pairwise commute. 

\emph{Example:} If $w$ is a potential in the empty set of variables, 
$\hmf_w$ is the homotopy category of 2-periodic complexes of 
graded $\Q$-vector spaces with finite-dimensional cohomology. 
Then, any object of $\hmf_w$ is isomorphic to its cohomology, 
which is a $\Z\oplus \Z_2$-graded vector space, and $\hmf_w$ is 
equivalent to the category of $\Z\oplus \Z_2$-graded 
finite-dimensional vector spaces. $K_w$ is equivalent to the 
homotopy category of this category, and, therefore, to the category 
of $\Z\oplus \Z \oplus \Z_2$-graded finite-dimensional vector spaces.  

 \vspace{0.1in} 

{\bf Tangles in a ball.} By a tangle $L$ we mean a proper embedding of 
an oriented compact $1$-manifold into a ball $B^3.$ We fix 
a great circle on the boundary $2$-sphere of $B^3$ and require 
that the boundary points of the embedded $1$-manifold lie 
on this great circle. A diagram $D$ of $L$ is a generic projection 
of $L$ onto the plane of the great circle. An isotopy of a tangle 
should not move its boundary points. 

\begin{figure} [htb] \drawing{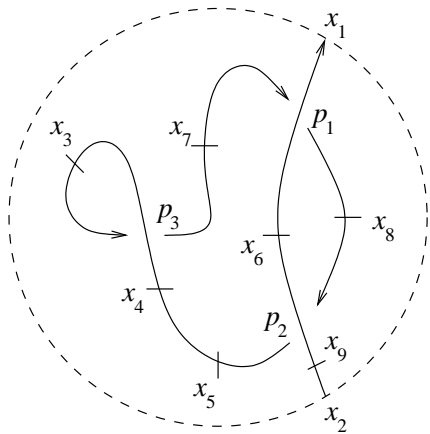}\caption{A marked 
 diagram of a tangle} \label{tangmarks}
 \end{figure}

A marked diagram (also denoted $D$) is a diagram with several 
marks placed on $D$ so that any segment bounded by crossings 
has at least one mark (see an example in figure~\ref{tangmarks}. 
Boundary points also count as marks. 

Let $m(D)$ be the set of marks of $D,$ and $\partial D$ 
the set of boundary points (its a subset of $m(D)$). 
Let $R$ be the ring of polynomials in $x_i,$ over $i\in m(D),$ 
and $R'$ the ring of polynomials in $x_i,$ over $i\in \partial D.$ 
 
\begin{figure} [htb] \drawing{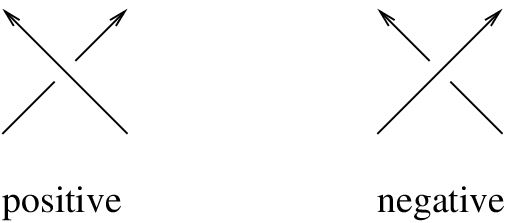}\caption{Positive and 
 negative crossings} \label{posneg}
 \end{figure}

We separate crossings of $D$ into positive and negative as 
explained in figure~\ref{posneg}. 

\begin{figure} [htb] \drawing{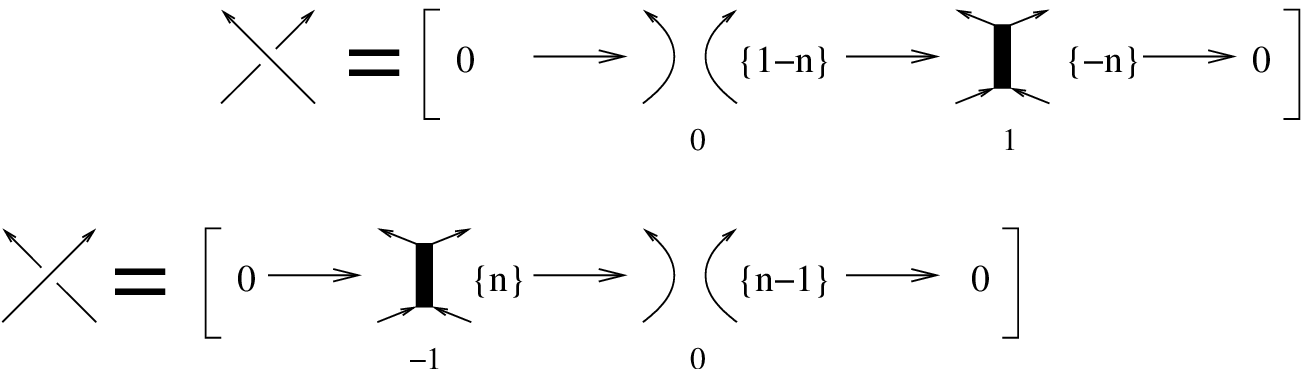}\caption{Complex assigned 
 to a crossing} \label{cones}
 \end{figure}

Given a crossing $p,$ let $\Gamma^0,\Gamma^1$ be its two 
resolutions (see figure~\ref{flats}). 
To a positive crossing assign the complex of factorizations (also see 
 figure~\ref{cones})
 \begin{equation*} 
  0 \lra C(\Gamma^0)\{ 1-n\} \stackrel{\chi_0}{\lra} 
 C(\Gamma^1)\{-n\} \lra 0. 
 \end{equation*} 
To a negative crossing assign the complex 
   \begin{equation*} 
  0 \lra C(\Gamma^1)\{ n\} \stackrel{\chi_1}{\lra} 
 C(\Gamma^0)\{n-1\} \lra 0. 
 \end{equation*}
In both cases we place $C(\Gamma^0)$ in cohomological 
degree $0.$ Denote this complex by $C^p.$ 

To $D$ associate the complex of factorizations $C(D)$ which 
is the tensor product of $C^p,$ over all crossings $p,$ of 
$L^i_j$ over all arcs $j\to i $ and of $A\bracket{1}$ over 
 all crossingless markless circles of $D$ (if such exist). 
The tensoring is done over appropriate polynomial rings so 
that $C(D),$ as an $R$-module, is free of finite rank. 

\vspace{0.1in} 

For instance, to produce $C(D),$ for $D$ in figure~\ref{tangmarks}, 
 we tensor $C^{p_1}$ with $C^{p_2}$ over $\Q[x_6,x_8],$ 
 and tensor the result with $C^{p_3}$ over $\Q[x_7].$ 
 
Proceed by tensoring $C^{p_1}\otimes C^{p_2}\otimes C^{p_3}$ 
with $L^4_5$ over $\Q[x_4,x_5],$ and with $L_2^9$ over 
 $\Q[x_9],$ so that 
 $$C(D) = C^{p_1}\otimes C^{p_2}\otimes C^{p_3} \otimes L^4_5
 \otimes L^9_2.$$ 

\vspace{0.1in} 
 
$C(D),$ for any diagram $D,$ is a complex of graded 
$(R',w)$-factorizations, where 
 $$ w = \sum_{i\in \partial D} \pm x_i^{n+1},$$ 
with signs determined by orientations of $D$ near boundary 
points. Thus, $C(D)$ is an object of the category $K_w.$ 

\begin{prop} Complexes $C(D)$ and $C(D')$ are canonically 
 isomorphic if $D'$ differs from $D$ only by marks. 
 \end{prop} 

This proposition follows at once from proposition~\ref{pr-can-iso}. 
 \hspace{0.06in} $\square$ 

\begin{theorem} \label{main-theorem} 
 Complexes $C(D)$ and $C(D')$ are isomorphic 
 in $\hmf_w$ if $D$ and $D'$ are two diagrams of the same 
 tangle $L.$ 
\end{theorem} 

To prove the theorem, it suffices to check it when $D$ and 
$D'$ are related by a single Reidemeister move. This is 
done in the next section. $\square$ 

\begin{corollary} The isomorphism class of the object $C(D)$ 
in the category $\hmf_w$ is an invariant of tangle $L.$ 
 \end{corollary}

 \vspace{0.1in} 

{\bf Link homology.} When $L$ is a link, the ring $R'=\Q,$ and 
 $\hmf_w$ is isomorphic to the category of finite-dimensional 
 $\Z\oplus \Z\oplus \Z_2$-graded $\Q$-vector spaces. 
 Cohomology groups are nontrivial only in the cyclic degree
 which is the number of components of $L$ modulo $2.$ This 
 reduces the grading of cohomology of $C(D)$ to $\Z\oplus \Z.$ 

 Denote the resulting cohomology groups by 
  $$ H_n(D) = \oplusop{i,j\in \Z} H_n^{i,j}(D).$$ 
It is clear from the construction that the Euler characteristic 
 of $H_n(D)$ is the polynomial $P_n(L),$ 
  $$ P_n(L) = \sum_{i,j\in \Z} (-1)^i q^j \dim_{\Q} H_n^{i,j}(D).$$ 
Isomorphism classes of vector spaces 
$H_n^{i,j}(D)$ depend only on $L.$ 
Denote by $h_n^{i,j}(L)$ the dimension of $H_n^{i,j}(D).$ 
This is an invariant of link $L.$ For each $n,$ we can put them 
together into a 2-variable polynomial invariant of $L,$   
 $$h_n(L) = \sum_{i,j} t^i q^j h_n^{i,j}(L).$$  

\vspace{0.15in} 

{\bf Reduced link homology.} Choose a component of $L,$ and a mark $i$ on $D$ 
that belongs to this component. For each resolution $\Gamma$ of 
$L,$ the vector space $H(\Gamma)$ is a free module over 
the ring $A\cong \Q[x_i]/(x_i^n).$ Let $\widetilde{\Q}$ be 
the one-dimensional graded $A$-module, placed in degree $0,$ 
and 
 $$\widetilde{C}(D) \define C(D) \otimes_A \widetilde{\Q}.$$
Complexes $\widetilde{C}(D)$ and $\widetilde{C}(D')$ are 
quasi-isomorphic if $D$ and $D'$ represent the same link with 
the same preferred component. 
Denote the cohomology groups of $\widetilde{C}(D)$ by 
 $\widetilde{H}_n^{i,j}(D),$ and their 
dimensions by $\widetilde{h}_n^{i,j}(L).$ Then 
 $$ \widetilde{h}_n(L) \define \sum_{i,j} t^i q^j\hspace{0.03in} 
  \widetilde{h}_n^{i,j}(L)$$  
is a two-variable polynomial invariant of $L$ with a preferred 
 component. Its specialization 
 to $t=-1$ is the one-variable polynomial 
 $\frac{P_n(L)}{[n]},$ which is another common normalization 
 for this one-variable specialization of HOMFLY.


 \section{Invariance under Reidemeister moves} 
 \label{sec-invar}  

 \begin{figure} [htb] \drawing{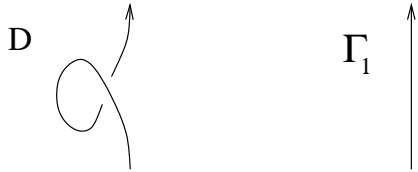}\caption{Type I move} 
 \label{curl1}
 \end{figure}

 {\bf Type I move.} Consider the figure~\ref{curl1} case of the type I 
move. We follow notations from the proof of 
proposition~\ref{prop-digonm1}. 
In particular, $\Gamma,\Gamma_1,$ and $\Gamma_2$ are as in 
figure~\ref{udigon2}. The complex $C(D)$ has 
the form 
$$0\lra C(\Gamma_2)\bracket{1}\{1-n\}\stackrel{\chi_0}{\lra}
 C(\Gamma)\{-n\} \lra 0.$$ 
Let $\widetilde{\alpha}_i,$  for $0\le i\le n-1,$ be the map 
 \begin{eqnarray*} 
  \widetilde{\alpha}_i & : & C(\Gamma_1)\{ 2i+2-2n\} 
    \lra C(\Gamma_2)\bracket{1}\{1-n\},   \\
  \widetilde{\alpha}_i & = & \sum_{j=0}^i m(x_1^j x_2^{i-j})
 \iota'.    
 \end{eqnarray*} 
Let   $Y_1\subset C(\Gamma_2)\bracket{1}\{1-n\}$ be the image of 
 $$\oplusop{0\le i\le n-2} C(\Gamma_1)\{2i+2-2n\}  $$ 
under the map   
$$\widetilde{\alpha} =\sum_{i=0}^{n-2} \widetilde{\alpha}_i,$$ 
and $Y_2\subset C(\Gamma_2)\bracket{1}\{1-n\}$ the image of 
$ C(\Gamma_1)$ under the map $\widetilde{\alpha}_{n-1}.$ 

There is a direct sum decomposition in $\hmf_w$ 
 $$ C(\Gamma_2)\bracket{1}\{1-n\} \cong Y_1 \oplus Y_2.$$
Furthermore, $\chi_0 (Y_2) =0,$ and the restriction of 
$\chi_0$ to $Y_1$ is an isomorphism from $Y_1$ to $C(\Gamma)\{-n\}.$ 
Therefore, in the category $K_w,$ the complex $C(D)$ 
is isomorphic to the direct sum of 
 $$ 0 \lra Y_1 \stackrel{\cong}{\lra} C(\Gamma)\{-n\} \lra 0$$ 
and 
 $$ 0 \lra Y_2 \lra 0,$$ 
with $Y_2$ in cohomological degree $0.$ Since the first summand 
is contractible, $C(D)$ is isomorphic to $Y_2\cong C(\Gamma_1)$ 
in the category $K_w.$ 

The invariance under other cases of the type I move can be 
verified similarly. 

 \vspace{0.2in} 

 \begin{figure} [htb] \drawing{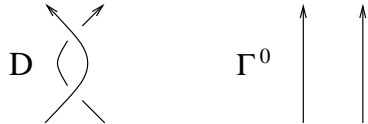}\caption{Type IIa move} 
 \label{ta1}
 \end{figure}

\begin{figure} [htb] \drawing{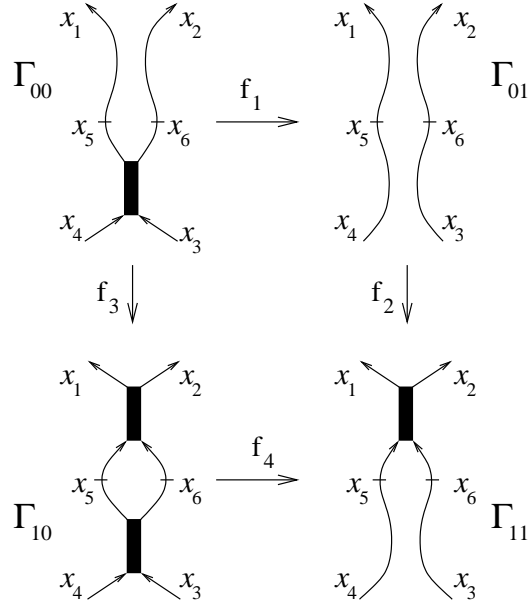}\caption{Four resolutions 
 of $D$ in the type IIa move} 
 \label{reid2a1}
 \end{figure}

 {\bf Type IIa move.} Consider diagrams $D$ and $\Gamma^0$ in 
figure~\ref{ta1}. The complex $C(D),$   
$$ 0 \lra C^{-1}(D) \stackrel{\partial^{-1}}{\lra} 
  C^0(D) \stackrel{\partial^0}{\lra} C^1(D) \lra 0 $$
has the form (see figure~\ref{reid2a1})  
 $$ 0 \lra C(\Gamma_{00})\{ 1 \} \xrightarrow{( f_1, f_3)^t} 
   \begin{array}{c} C(\Gamma_{01}) \\   \oplus   \\
  C(\Gamma_{10}) \end{array} \xrightarrow{(f_2, -f_4)} 
   C(\Gamma_{11})\{-1\} \lra 0 $$ 
where $f_1,f_4$ are given by the map $\chi_1$ (corresponding to the 
topology change in the lower halves of diagrams 
$\Gamma_{00}, \Gamma_{10}$), and $f_2,f_3$ are given by $\chi_0$ 
 (topology change in the upper halves
 of diagrams $\Gamma_{00},\Gamma_{01}$). 
  
We know that 
\begin{eqnarray} 
\label{iso-101} 
C(\Gamma_{10})& \cong & C(\Gamma^1)\{ 1\} \oplus C(\Gamma^1)\{-1\}, \\ 
 \label{iso-010} C(\Gamma_{01})& \cong & C(\Gamma^0), \\
 \label{iso-0011} 
  C(\Gamma_{00})& \cong & C(\Gamma_{11}) \cong C(\Gamma^1).
\end{eqnarray}  
where $\Gamma^0,\Gamma^1$ are diagrams in figure~\ref{pair2}, and 
here and further on the isomorphisms are in the category $\hmf_w,$ 
for the potential $w=x_1^{n+1}+x_2^{n+1}-x_3^{n+1}-x_4^{n+1}.$ 

Fix isomorphisms (\ref{iso-101}), (\ref{iso-010}), 
(\ref{iso-0011}). Under these isomorphisms, $f_3,f_4$ become 
two-component maps 
$$ f_3=(f_{03},f_{13})^t, \hspace{0.2in} 
    f_4= (f_{04}, f_{14}),$$
where, for instance $f_{03}$ (respectively, $f_{13}$) 
 is a degree $0$ (respectively, degree $2$) endomorphism of 
$C(\Gamma^1).$ 

\begin{lemma} \label{lemma-3term} 
\begin{enumerate}\item $C(\Gamma^1)$ has no negative degree 
endomorphisms in $\hmf_w.$ The only degree $0$ endomorphisms 
are multiples of the identity endomorphism. 
\item If $n>2,$ the space of degree $2$ endomorphisms of 
$C(\Gamma^1)$ is 3-dimensional. Multiplications by $x_1,x_2,x_3,x_4$
span this space, with the only relation 
 $m(x_1+x_2-x_3-x_4)=0.$ 
\item If $n=2,$ the space of degree $2$ endomorphisms of 
$C(\Gamma^1)$ is 2-dimensional. Multiplications by $x_1,x_2,x_3,x_4$ 
span this space, with relations $m(x_1+x_2)=0$ and $m(x_3+x_4)=0.$ 
\end{enumerate} 
\end{lemma} 

\emph{Proof:} since $C(\Gamma^1)\{1\}$ is the Koszul factorization 
for the pair 
$$((u_1,u_2),(x_1+x_2-x_3-x_4, x_1x_2-x_3x_4)),$$ 
the complex $\Hom_R(C(\Gamma^1),C(\Gamma^1))$ is isomorphic to the 
Koszul complex of the sequence
\begin{equation}\label{eq-fourt} 
(x_1+x_2-x_3-x_4,x_1x_2-x_3x_4,  u_1, u_2),
\end{equation} 
with the grading collapsed from $\Z$ to $\Z_2$ (for the definition 
of $u_1,u_2$ see Section~\ref{sec-planar}). This sequence 
is regular, so that the cohomology of this 2-complex is 
$$\Q[x_1,x_2,x_3,x_4]/(x_1+x_2-x_3-x_4,x_1x_2-x_3x_4,  u_1, u_2).$$
This cohomology is also isomorphic, as a graded $\Q$-algebra, 
to $\mathrm{End}_{\HMF}(C(\Gamma^1)),$ which implies part $1$ 
of the lemma.   
If $n>2,$ all terms in (\ref{eq-fourt}) but the first are homogeneous 
polynomials in $x$'s that are at least quadratic. Part $2$ follows. 
If $n=2,$ the term $u_2$  is linear, and we get 
additional relation $m(x_3+x_4)=0.$  \hspace{0.15in} 
$\square$ 

\emph{Remark:} If $n=1,$ there is nothing to prove, since $C(\Gamma)=0$ 
whenever $\Gamma$ has a wide edge, and the whole story becomes 
trivial. 

From the lemma we deduce that $f_{03}$ and $f_{14}$ are rational 
multiples of the identity endomorphism, while $f_{13}, f_{04}$ are 
multiplications by linear combinations 
of $x_1,x_2,x_3,x_4.$ 

$C(D)$ is a complex, so that $\partial^{0}\partial^{-1}=0.$ We compute 
 $$0 = \partial^0\partial^{-1}= f_2f_1 - f_4f_3= 
  m(x_1-x_3) - f_{14}f_{13}-f_{04}f_{03}.$$ 
Therefore, either $f_{14}\not= 0,$ or $f_{03}\not= 0,$ since the multiplication by 
$x_1-x_3$ is a nontrivial endomorphism of $C(\Gamma^1).$ 

Assume that $f_{14}\not=0.$ Then $f_{14}$ is a nonzero multiple of the identity.  
By rescaling, we normalize $f_{14}$ to be the identity.

\begin{lemma} \label{lemma-notz} $f_{03}\not= 0.$ \end{lemma} 

\emph{Proof:} Assume otherwise, $f_{03}=0.$ Then  $f_{13}= m(x_1-x_3).$  
The  composition of $f_3$ with the map $\chi_1$ in the 
opposite direction, where the topology change takes place 
around the upper wide edge of $\Gamma_{10},$ is the multiplication 
by $x_1-x_6,$ which is the same as multiplication by $x_1-x_2$ 
(since multiplications by $x_2$ and $x_6$ are homotopic 
endomorphisms of $C(\Gamma_{00})$). The map $\chi_1,$ restricted 
to the summand $C(\Gamma^1)\{-1\}$ of $C(\Gamma_{10}),$ is a 
multiplication by some rational number $z.$ Under the assumption 
$f_{03}=0,$ we have 
$$ m(x_1-x_2) = \chi_1 f_3 = \chi_1 f_{13}=z\cdot  m(x_1-x_3).$$ 
Therefore, 
$$ m((1-z)x_1-x_2+z x_3) =0.$$ 
This is impossible, by lemma~\ref{lemma-3term}. Lemma~\ref{lemma-notz}
follows. $\square$  

\vspace{0.1in} 

Likewise, if $f_{03}\not=0,$ we can show that $f_{14}\not=0.$ Thus, 
both $f_{03}$ and $f_{14}$ are nonzero multiples of the identity morphism. 
We normalize so that each of them is the identity map. The differential 
$$\partial^{-1}\hspace{0.1in} : \hspace{0.1in} C^{-1}(D) \lra C^0(D) $$ 
is split injective. Since $f_{14}$ is the identity, the 
differential
$$\partial^0 \hspace{0.1in} : \hspace{0.1in} C^0(D) \lra C^1(D) $$
is split surjective. We can  twist the direct sum decomposition 
$$C^0(D) \cong C(\Gamma^0) \oplus C(\Gamma^1)\{1\} 
\oplus C(\Gamma^1)\{-1\}$$ 
so that $f_{03},f_{14}$ become the only nonzero entries in the 
matrices describing differentials $\partial^{-1},\partial^0,$ and  
 $C(D)$ breaks into the direct sum of three complexes 
\begin{eqnarray*} 
  0 \lra & C(\Gamma^0) & \lra 0,    \\
 0 \lra C(\Gamma^1)\{1\} \stackrel{\cong}{\lra} & C(\Gamma^1)\{1\} & 
  \lra 0,   \\
 0 \lra & C(\Gamma^1)\{-1\}&  \stackrel{\cong}{\lra} C(\Gamma^1)\{-1\}  
  \lra 0.  
\end{eqnarray*} 
The last two complexes are contractible. Therefore, $C(D)$ and 
$C(\Gamma^0)$ are isomorphic in the category $K_w.$ 

\vspace{0.2in}

 \begin{figure} [htb] \drawing{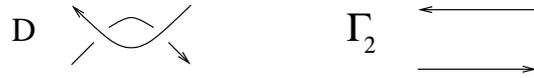}\caption{Type IIb move} 
 \label{tb1}
 \end{figure}

\begin{figure} [htb] \drawing{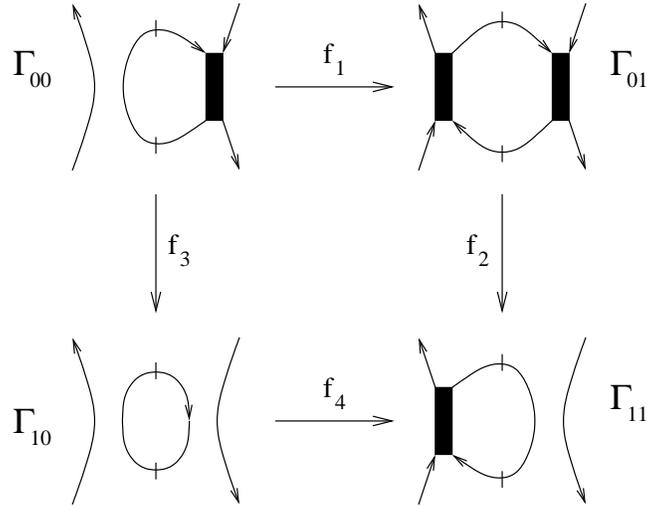}\caption{Commutative square 
of resolutions of $D,$ type IIb move} 
 \label{reid2b1}
 \end{figure}

 {\bf Type IIb move.} Consider diagrams $D$ and $\Gamma_2$ in figure~\ref{tb1}.
The complex $C(D)$ has the form (see figure~\ref{reid2b1}) 
 $$ 0 \lra C(\Gamma_{00})\{ 1 \} \xrightarrow{( f_1, f_3)^t} 
   \begin{array}{c} C(\Gamma_{01}) \\   \oplus   \\
  C(\Gamma_{10}) \end{array} \xrightarrow{(f_2, -f_4)} 
   C(\Gamma_{11})\{-1\} \lra 0. $$ 

There are isomorphisms in $\hmf_w$ 
\begin{eqnarray} 
  C(\Gamma_{00})\{1\} & \cong &  \label{iso-h1} 
   {\mathop{\oplus}\limits_{i=0}^{n-2}} C(\Gamma_1)\bracket{1}\{3-n+2i\}   \\
  C(\Gamma_{11})\{-1\} & \cong & {\mathop{\oplus}\limits_{i=0}^{n-2}}
 C(\Gamma_1)\bracket{1}\{1-n+2i\}   \\
  C(\Gamma_{10}) & \cong & {\mathop{\oplus}\limits_{i=0}^{n-1}}
  C(\Gamma_1)\bracket{1}\{1-n+2i\}  \\     \label{iso-h4} 
  C(\Gamma_{01}) & \cong & \biggl( {\mathop{\oplus}\limits_{i=0}^{n-3}}
 C(\Gamma_1)\langle 1 \rangle \{ 3-n+2i\} \biggr)\oplus  C(\Gamma_2)  
\end{eqnarray} 
where $\Gamma_1$ and $\Gamma_2$ are as in figure~\ref{square1}. 

The proof of the invariance under type I move implies that $f_4$ 
is split surjective. Likewise, $f_3$ is split injective. Since 
the category $\hmf_w$ has splitting idempotents, we can decompose 
$C^0(D)$ into the direct sum 
$$C^0(D) \cong  \mathrm{Im}(\partial^{-1})  \oplus Y_1  \oplus Y_2  $$
such that $\partial^0$ restricts to an isomorphism from $Y_1$ to 
 $C(\Gamma_{11})\{-1\}$ and $\partial^0(Y_2)=0.$ Therefore, 
$C(D)$ is isomorphic to the direct sum of complexes  
\begin{eqnarray*} 
    0 \lra & Y_2  & \lra 0,    \\
   0 \lra C(\Gamma_{00})\{1\} \stackrel{\cong}{\lra} & 
 \mathrm{Im}(\partial^{-1}) & \lra 0,   \\
   0 \lra & Y_1 & \stackrel{\cong}{\lra} C(\Gamma_{11})\{-1\} \lra 0 . 
\end{eqnarray*} 
From formulas (\ref{iso-h1})-(\ref{iso-h4}) we obtain  
$$C^0(D) \cong C(\Gamma_{01})\oplus C(\Gamma_{10}) \cong 
  C(\Gamma_{00})\{1\}\oplus C(\Gamma_{11})\{-1\} \oplus C(\Gamma_2).$$
Category $\hmf_w$ is Krull-Schmidt (objects have unique direct sum 
decompositions into indecomposables).  
Therefore, $Y_2\cong C(\Gamma_2)$ and complexes $C(D)$ and 
 $ 0 \lra C(\Gamma_2) \lra 0$ are isomorphic. This concludes our 
 proof of the invariance under the type IIb move. 

 \vspace{0.2in}

\begin{figure} [htb] \drawing{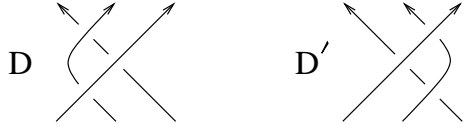}\caption{Type III move} 
 \label{braid1}
 \end{figure}

 {\bf Type III move.} We need to show that $C(D)$ and $C(D')$ are isomorphic  for 
$D,D'$ in figure~\ref{braid1}. 

\begin{figure} [htb] \drawing{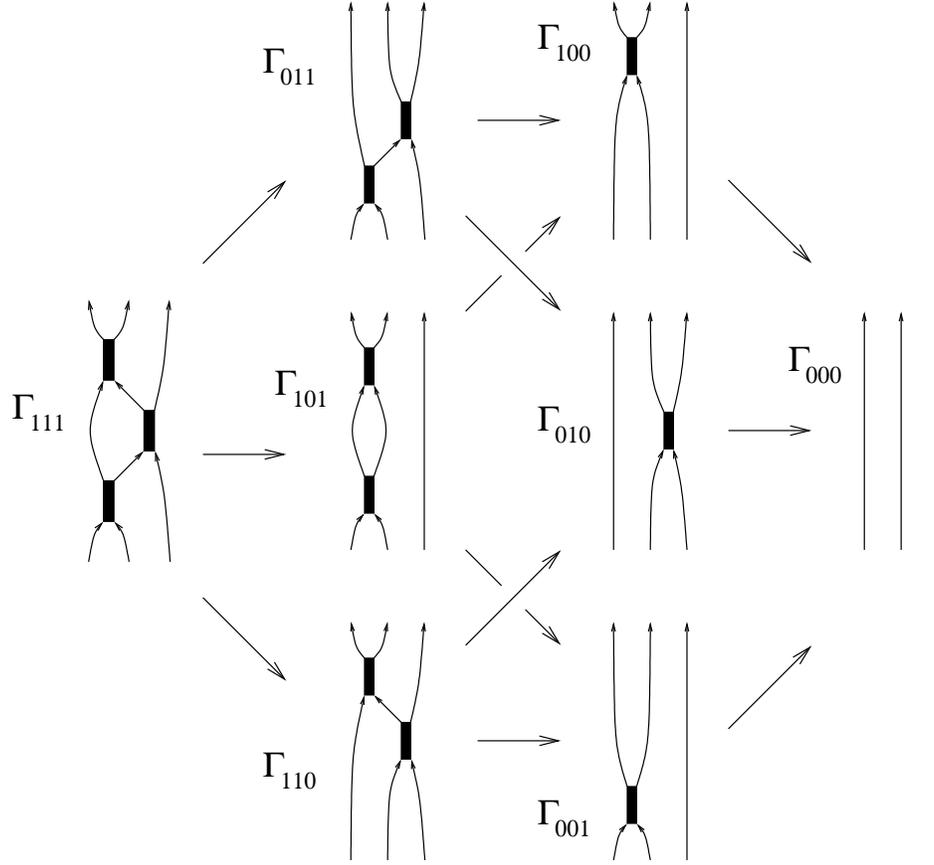}\caption{Resolution cube 
 of $D$} 
 \label{rthree1}
 \end{figure}

The diagram $D$ has $8$ resolutions, denoted 
$\Gamma_{ijk},$ for $i,j,k\in \{0,1\},$ see figure~\ref{rthree1}. The complex 
$C(D)\{-3n\}$ has the form 
\begin{eqnarray*} 
 0 \lra C(\Gamma_{111}) \stackrel{d^{-3}}{\lra}  & \left( \begin{array}{c} C(\Gamma_{011})\{-1\} \\ 
       C(\Gamma_{010}) \{-1\} \\   C(\Gamma_{110})\{-1\} \end{array} \right) &  
     \stackrel{d^{-2}}{\lra}  \\
    \stackrel{d^{-2}}{\lra}  &    \left( \begin{array}{c} C(\Gamma_{100})\{-2\} \\ 
       C(\Gamma_{010}) \{-2\} \\   C(\Gamma_{001})\{-2\} \end{array}\right)  &  
    \stackrel{d^{-1}}{\lra} C(\Gamma_{000})\{-3\} \lra 0 
\end{eqnarray*} 
with $C(\Gamma_{000})\{-3\}$ in cohomological degree $0.$ We depicted this complex 
 in figure~\ref{rthree1}. 

Diagrams $\Gamma_{100}$ and $\Gamma_{001}$ are isotopic, so that $C(\Gamma_{100})$ 
and $C(\Gamma_{001})$ are isomorphic factorizations. Moreover, 
 \begin{equation} \label{a-dir-sum}   
 C(\Gamma_{101}) \cong  C(\Gamma_{100})\{1\} \oplus C(\Gamma_{100})\{-1\},
 \end{equation} 
 and from proposition~\ref{core-isom} we know that 
 $$C(\Gamma_{111}) \cong C(\Gamma_{100}) \oplus \Upsilon.$$ 
Differential $d^{-3}$ is injective on $C(\Gamma_{100}) \subset C(\Gamma_{111}).$ 
In fact, its middle component (the map to $C(\Gamma_{101})\{-1\}$) is injective, which follows 
from our construction of the inclusion  $C(\Gamma_{100}) \subset C(\Gamma_{111})$
and the proof of the invariance under type IIa move. 

The factorization $d^{-3}(C(\Gamma_{100}))$ is a direct summand of $C^{-2}(D)\{-3n\}.$ 
Thus,  $C(D)\{-3n\}$ contains a contractible summand 
 \begin{equation}\label{dir-s-1}
  0 \lra C(\Gamma_{100}) \stackrel{d^{-3}}{\lra} C(\Gamma_{100}) \lra 0  .
 \end{equation} 

 Direct sum decomposition (\ref{a-dir-sum}) can be selected so that 
 $$ C(\Gamma_{101})\{-1\} \cong p_{101}d^{-3} C(\Gamma_{100}) \oplus C(\Gamma_{100})\{-2\}$$ 
where $p_{101}$ is the projection onto the middle summand of $C^{-2}(D)\{-3n\}.$ 

\begin{figure} [htb] \drawing{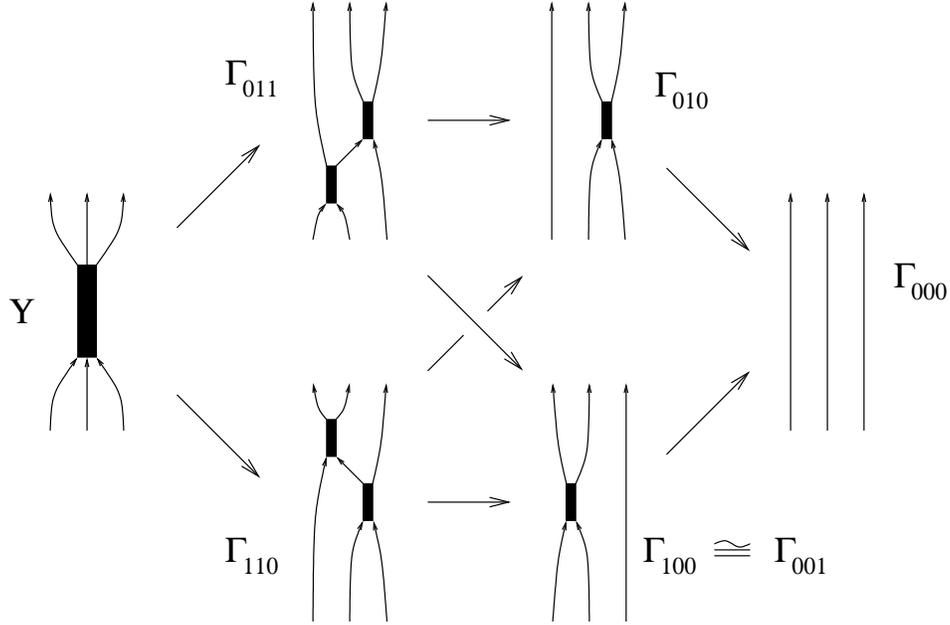}\caption{Complex $C$} 
 \label{rthree2}
 \end{figure}

The differential $d^{-2}$ is injective on $C(\Gamma_{100})\{-2\}\subset C(\Gamma_{101})\{-1\},$ 
since its middle component is a nonzero multiple of the 
identity. Furthermore, the image of $C(\Gamma_{100})\{-2\}\subset C(\Gamma_{101})\{-1\}$
under $d^{-2}$ is a direct summand of $C^{-1}(D).$ Hence, the complex $C(D)\{-3n\}$ 
contains a contractible direct summand isomorphic to  
  \begin{equation}\label{dir-s-2} 
  0 \lra C(\Gamma_{100})\{-2\}  \stackrel{d^{-2}}{\lra} 
                 C(\Gamma_{100})\{-2\} \lra 0.
  \end{equation}  
After splitting off contractible direct summands (\ref{dir-s-1}) and (\ref{dir-s-2}), the 
complex $C(D)\{-3n\}$ reduces to the complex $C$ of the form 
  \begin{eqnarray*} 
 0 \lra \Upsilon\stackrel{d^{-3}}{\lra}  & \left( \begin{array}{c} C(\Gamma_{011})\{-1\} \\ 
       C(\Gamma_{110})\{-1\} \end{array} \right) &  
     \stackrel{d^{-2}}{\lra}  \\
    \stackrel{d^{-2}}{\lra}  &    \left( \begin{array}{c} C(\Gamma_{010})\{-2\} \\ 
       C(\Gamma_{100})\{-2\} \end{array}\right)  &  
    \stackrel{d^{-1}}{\lra} C(\Gamma_{000})\{-3\} \lra 0, 
\end{eqnarray*} 
see figure~\ref{rthree2}. To diagram $Y$ on the left of the figure we assign 
factorization $\Upsilon.$ 

In the rest of the section we assume that $n>2$ (proofs for $n=2$ 
are easier since then $\Upsilon=0$). 

\begin{lemma} For each arrow in figure~\ref{rthree2} 
with some diagram $Z_1$ as the source and $Z_2$ as the target, the 
space of grading-preserving morphisms 
   $$ \Hom_{\hmf}(C(Z_1), C(Z_2)\{-1\})$$ 
is one-dimensional. 
\end{lemma} 

The lemma admits a straightforward proof similar to the proof of corollary~\ref{space-q}. 
Details are omitted. $\square$ 

\vspace{0.1in} 

\begin{figure} [htb] \drawing{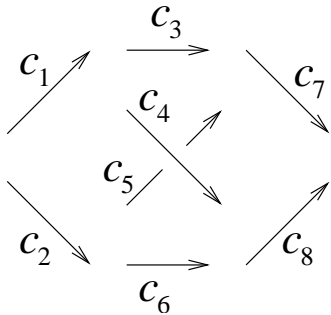}\caption{Indices and arrows} 
 \label{rthree3}
 \end{figure}

Denote by $c_1, \dots, c_8$ nontrivial morphisms from the above lemma. Each morphism 
is determined up to multiplication by a nonzero rational number. Figure~\ref{rthree3} 
shows how the indices $i$ of $c_i$'s match the eight arrows of figure~\ref{rthree2}. 
We choose $c_i$'s so that each square in figure~\ref{rthree3} commutes. 

The differential in $C$ is grading-preserving. Therefore, it has the form 
$$ d^{-3}= \left( \begin{array}{c} \lambda_1 c_1 \\  \lambda_2 c_2 \end{array} 
   \right), \hspace{0.1in} 
    d^{-2} = \left( \begin{array}{cc} \lambda_3 c_3 & \lambda_5 c_5  \\
                                                           \lambda_4 c_4 & \lambda_6 c_6 
                                   \end{array} \right), \hspace{0.1in} 
    d^{-1} = \left( \begin{array}{c} \lambda_7 c_7 \\  \lambda_8 c_8  \end{array} 
                 \right)
$$ 
for some rational numbers $\lambda_1, \dots, \lambda_8.$ 

\begin{lemma} For any two composable arrows $Z_1 \lra Z_2 \lra Z_3$ in 
 figure~\ref{rthree2}, the corresponding homomorphism of factorizations 
   $$  c_j  c_i \hspace{0.1in} : \hspace{0.1in} C(Z_1) \lra C(Z_3)\{-2\} $$ 
 is nonzero in $\hmf_w.$ 
\end{lemma} 

\begin{figure} [htb] \drawing{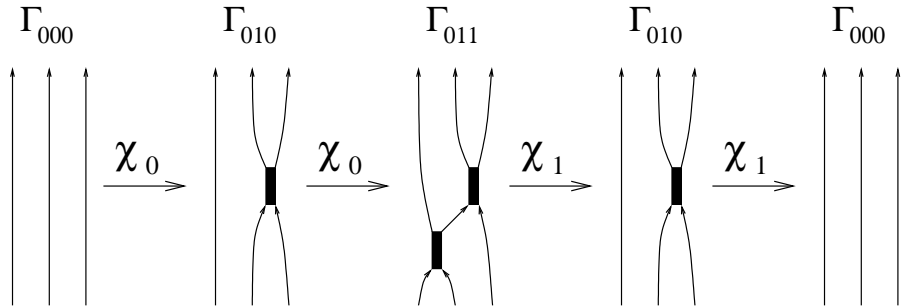}\caption{Composing $c_7c_3$ with its dual} 
 \label{rthree4}
 \end{figure}

 \emph{Proof:} For instance, to show that $c_7 c_3\not=0,$ 
we note that this composition is the product of two $\chi_1$ maps (up to rescaling by 
a nonzero rational number), one for each wide 
edge of $\Gamma_{011}.$ Compose it with the "dual" product of two $\chi_0$ maps, 
as in figure~\ref{rthree4}. Assign variables $x_1, \dots, x_6$ to the endpoints of 
our diagrams in the same way as in figure~\ref{mapgone}. Then the figure~\ref{rthree4} 
map is $m(x_1-x_5)m(x_2-x_6).$ This is a nontrivial endomorphism of $C(\Gamma_{000}),$ 
so that the composition $c_7 c_3$ is nonzero. 

Identical arguments 
take care of the other three compositions that end in $\Gamma_{000}.$ Similar 
but longer computations can be used to check nontriviality of each of the four 
compositions $c_j c_i$ that start at $Y.$ \hspace{0.06in} $\square$ 

The category $K_w$ has an automorphism (tensoring with the complex assigned to the inverse 
braid of $D$) that takes $C(D)$ to $C(D'')$ where 
$D''$ is the diagram made of three disjoint oriented arcs. Since $C(D'')$ is indecomposable 
in $K_w,$ the same is true of $C(D)$ and $C\cong C(D).$ 
Indecomposability of $C$ together with 
the above lemma implies that none of the coefficients $\lambda_1, \dots, \lambda_8$ is $0.$  
By rescaling, we can set 
  $$ \lambda_1 = \lambda_2 =\lambda_4=\lambda_5=\lambda_7=\lambda_8=1, 
  \hspace{0.1in} \lambda_3=\lambda_6=-1. $$ 

The complex $C$ is therefore invariant under the "flip" which transposes $x_1$ with $x_3$ 
and $x_4$ with $x_6.$ This flip takes $C(D)$ to $C(D'),$ thus, $C$ is isomorphic in $K_w$ 
to $C(D').$ We have $C(D) \cong C \cong C(D').$ The invariance under the Reidemeister 
move III in figure~\ref{braid1} follows.


 \section{Factorizations and 2-dimensional TQFTs with corners} 
 \label{sec-corners}

Let $I$ be a finite set and $s:I\to \{ 1,-1\}$ a map, which we call 
the orientation map. We say that $s$ is \emph{balanced} if  the sets $s^{-1}(1)$ and 
$s^{-1}(-1)$ have the same cardinality (i.e.,  $s$ takes as many elements of $I$ to $1$ 
as it does to $-1$). From now on we assume that $s$ is balanced. 

To $(I,s)$ we assign the category 
$\mathrm{Cob}_{I,s}$ whose objects are oriented one-dimensional 
manifolds $N$ with boundary $I,$ such that the orientation of $N$ induces 
orientation $s$ on $I= \partial N.$  The 
morphisms from $N_0$ to $N_1$ are oriented 2-dimensional surfaces $S$ 
with boundary $N_0 \cup -N_1 \cup I\times [0,1]$ and corners 
 $I\times \{0\}\sqcup I \times \{1\}.$    

Let $R_I$ be the ring of polynomials in $x_i,i\in I,$ and 
$$w(I,s) = \sum_{i\in I} s(i) x_i^{n+1} \hspace{0.03in} \in R_I.$$

\begin{figure} [htb] \drawing{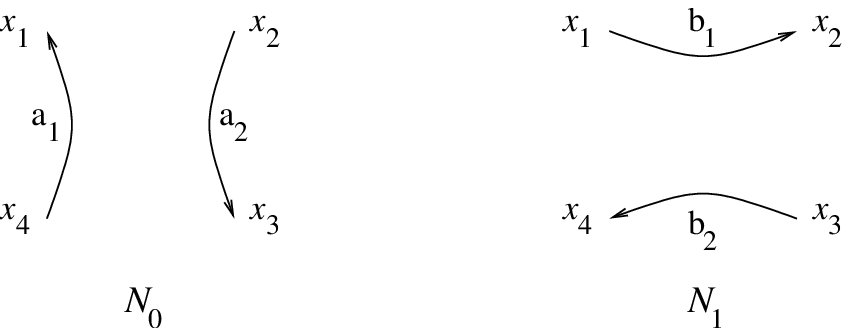}
 \caption{} \label{saddle}
 \end{figure}

To an object $N\in \mathrm{Cob}_{I,s}$ we can assign a  factorization 
 $C(N)$ with potential $w(I,s).$ The factorization 
is the tensor product (over $\Q$) of $L^i_j,$ over 
all arcs $j\to i$ in $N,$ and of $A\bracket{1},$ one for each circle in $N.$ 
We would like to extend this assignment to a functor from 
$\mathrm{Cob}_{I,s}$ to the category of factorizations. 
To do this, pick a morphism $S$ in $\mathrm{Cob}_{I,s}$ and 
write it as a composition of simple cobordisms (cobordisms 
with only one critical point). To the cobordisms of creation and 
annihilation of a circle we assign maps $\iota, \varepsilon,$
 defined in Section~\ref{sec-functors}. To the saddle point cobordism between 
one-manifolds $N_0,N_1$ in figure~\ref{saddle} we want to associate  
a map 
 $$ \eta: C(N_0) \lra C(N_1)\bracket{1}\{ 1-n\}.$$ 
To define $\eta,$ assign labels $a_1,a_2,b_1,b_2$ to the 
components of $N_0,N_1$ as indicated in figure~\ref{saddle}.

Factorization $C(N_0)$ has the form 
\begin{equation*}
  \left( \begin{array}{c}R(\emptyset) \\ R(a_1a_2)\{2-2n\}\end{array}
   \right)
 \xrightarrow{P_0}
  \left( \begin{array}{c} R(a_1)\{1-n\} \\ R(a_2)\{1-n\}  \end{array}
 \right)
 \xrightarrow{P_0'}
 \left( \begin{array}{c}R(\emptyset) \\ R(a_1a_2)\{2-2n\}\end{array}
   \right), 
\end{equation*}
where
\begin{equation*}
P_0=\left(\begin{array}{cc} \pi_{14}  & x_3-x_2 \\
      \pi_{32}  & x_4-x_1 \end{array} \right), \hspace{0.1in}
P_0'=\left(\begin{array}{cc} x_1-x_4  & x_3-x_2 \\
     \pi_{32}  & - \pi_{14} \end{array} \right) , 
\end{equation*} 
and $R=\Q[x_1,x_2,x_3,x_4].$

Factorization $C(N_1)\bracket{1}$ is given by 
\begin{equation*}
  \left( \begin{array}{c}R(b_1) \{1-n\}\\ R(b_2)\{1-n\}\end{array}
   \right)
 \xrightarrow{P_1}
  \left( \begin{array}{c} R(\emptyset) \\ R(b_1b_2)\{2-2n\}  \end{array}
 \right)
 \xrightarrow{P_1'}
 \left( \begin{array}{c}R(b_1) \{1-n\}\\ R(b_2)\{1-n\}\end{array}
   \right), 
\end{equation*}
where
\begin{equation*}
P_1=\left(\begin{array}{cc} x_2-x_1  & x_4-x_3 \\
      -\pi_{34}  & \pi_{12} \end{array} \right), \hspace{0.1in}
P_1'=\left(\begin{array}{cc} -t_{12}  & x_4-x_3 \\
     -\pi_{34}  & x_1-x_2 \end{array} \right) , 
\end{equation*} 
and $R=\Q[x_1,x_2,x_3,x_4].$ 

$\eta$ is given by the pair of matrices 
\begin{eqnarray*}
  & &   \left( \begin{array}{cc}  e_{123}+e_{124}+(x_4-x_3)r &  1 \\
                -e_{134}-e_{234} + (x_1-x_2) r  & 1 \end{array}\right) ,     \\
  & &   \left( \begin{array}{cc}  -1  &  1 \\
               -e_{123}-e_{234} + (x_1-x_4) r & -e_{134}-e_{124}+(x_3-x_2)r  \end{array}\right) ,
\end{eqnarray*} 
with 
 $$ e_{ijk} = \sum_{a+b+c=n-1} x_i^a x_j^b x_k^c$$ 
and $r$ being an arbitrary polynomial of degree $n-2$ in $x_1, \dots, x_4.$ Up to chain 
homotopy, $\eta$ does not depend on $r.$ 

 $\eta$  generates the $R$-module of homomorphisms 
 $\Hom_{\HMF}(C(N_0), C(N_1)\bracket{1}),$ and we should assign $\eta$ to the saddle 
point cobordism.  There is a problem, though. The diagrams 
in figure~\ref{saddle} are invariant under the rotation by $180^{\circ},$ but $\eta$ acquires a 
sign after the rotation. Namely, if we transpose $x_1$ with $x_3,$ $x_2$ with $x_4,$ 
$a_1$ with $a_2$ and $b_1$ with $b_2$ in the formula for $\eta,$ the resulting 
map is $-\eta.$ Thus, $\eta$ can be canonically defined only up to a sign. 

If the saddle point cobordisms takes place not between arcs, but between components of $N_0,N_1$ 
some of which are circles,  we  add marks to each component, 
select a suitable pair of arcs bounded by marks (and, possibly, by boundary points),  
apply $\eta$ to this pair, and, finally, erase marks. It is easy to see that, up to a sign,
 the resulting map does not depend on intermediate choices. 

Let $\mathrm{H}'_{I,s}$ be the category with the same objects as 
$\hmf_{w(I,s)}$ but morphisms being Ext groups $\mathrm{Ext}_{\HMF}(M_0,M_1),$  with 
$f$ and $-f$ identified for all $f\in \mathrm{Ext}_{\HMF}(M_0,M_1).$ We use Ext groups 
rather than just Hom's since the shift $\bracket{1}$ is built into $\iota, \varepsilon, $ and $\eta.$

Given a surface $S$ which is an object of $\mathrm{Cob}_{I,s},$ write it as a product 
of cobordisms with only one critical point and define $C(S)$ as the corresponding product of 
 $\iota, \varepsilon,$ and $\eta$'s. Recall that our definition of $\varepsilon$ contained 
a parameter $\zeta\in \Q^{\ast}.$ To make $\varepsilon$ compatible with $\eta$ we 
must set $\zeta$ to either $\frac{1}{n+1}$ or $-\frac{1}{n+1}$ (so that  
 $\varepsilon \eta= \pm\mathrm{Id}$  if the saddle point cobordism 
goes from an arc to the union of a circle and an arc).   

\begin{prop} $\pm C(S)$ does not depend on the presentation of $S$ as a product 
of elementary cobordisms. \end{prop} 

The proof is left to the reader. $\square$ 

Thus, we obtain a functor from the cobordism category $\mathrm{Cob}_{I,s}$ to 
the category $\mathrm{H}'_{I,s}.$ The shifts $\upsilon_1(S),\upsilon_2(S)$ in 
 $$ \pm C(S)  \hspace{0.1in} : \hspace{0.1in} C(N_0) \lra C(N_1) \bracket{v_1(S)}\{ v_2(S)\}$$ 
are as follows. Glue $N_0$ and $N_1$ along the common boundary $I$, and 
count the number $\upsilon$ of components in the closed 1-manifold that results. $\upsilon_1(S)$ is the 
parity of $\upsilon+\frac{|I|}{2},$  while 
$$\upsilon_2(S)= (n-1)(\chi(S) - \frac{|I|}{2}) $$ 
where $\chi(S)$ is the Euler characteristic of $S.$ 

\vspace{0.15in} 

{\bf 2-functor.}  In the above construction we fixed the boundary of one-manifolds. If we 
put together all categories $\mathrm{Cob}_{I,s},$ over various $I$ and $s$ (and consider 
decompositions of $I$ into pairs of disjoint sets, to view $N$ with $\partial N = I$ as 
a morphism) we get a 2-category of oriented surfaces with corners. 
Functors 
 $$ C: \mathrm{Cob}_{I,s} \lra \mathrm{H}'_{I,s}$$ 
extend to a two-functor from this 2-category of cobordisms to the 2-category of 
functorizations with potentials as objects, factorizations as morphisms, and 
elements of Ext groups between factorizations (with $f$ and $-f$ identified) 
as 2-morphisms. We leave details to the reader.

\section{Projective invariance for tangle cobordisms}

As before, we consider oriented tangles in a ball $B^3.$  Fix a great circle on the ball's boundary, 
choose a finite subset $I$ of this circle and a
balanced "orientation" function $s: I \to \{ 1,-1\}.$  
Let $\mathrm{TC}_{I,s}$ be the category of tangle cobordisms with objects--oriented tangles $L$ 
in $B^3$ with oriented boundary $(I,s)$ and morphisms from $L_0$ to $L_1$--oriented 
surfaces $S$ embedded in  $B^3\times [0,1]$ with boundary 
 $$ \partial S = L_0 \times \{ 0\} \cup L_1 \times \{ 1\} \cup I\times [0,1]$$ 
and corners $I\times \{0\} \cup I \times \{ 1\},$ up to isotopies that fix the boundary.  

 \begin{figure} [htb] \drawing{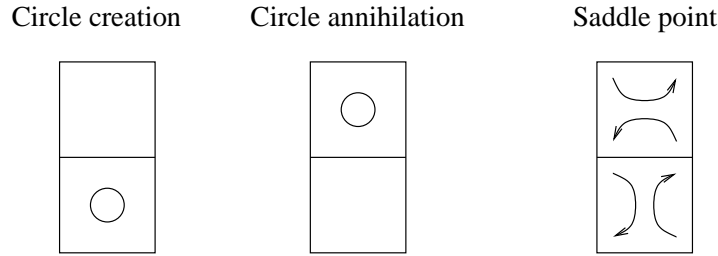}
 \caption{Movie moves of simple cobordisms} \label{morse}
 \end{figure}

\begin{figure} [htb] \drawing{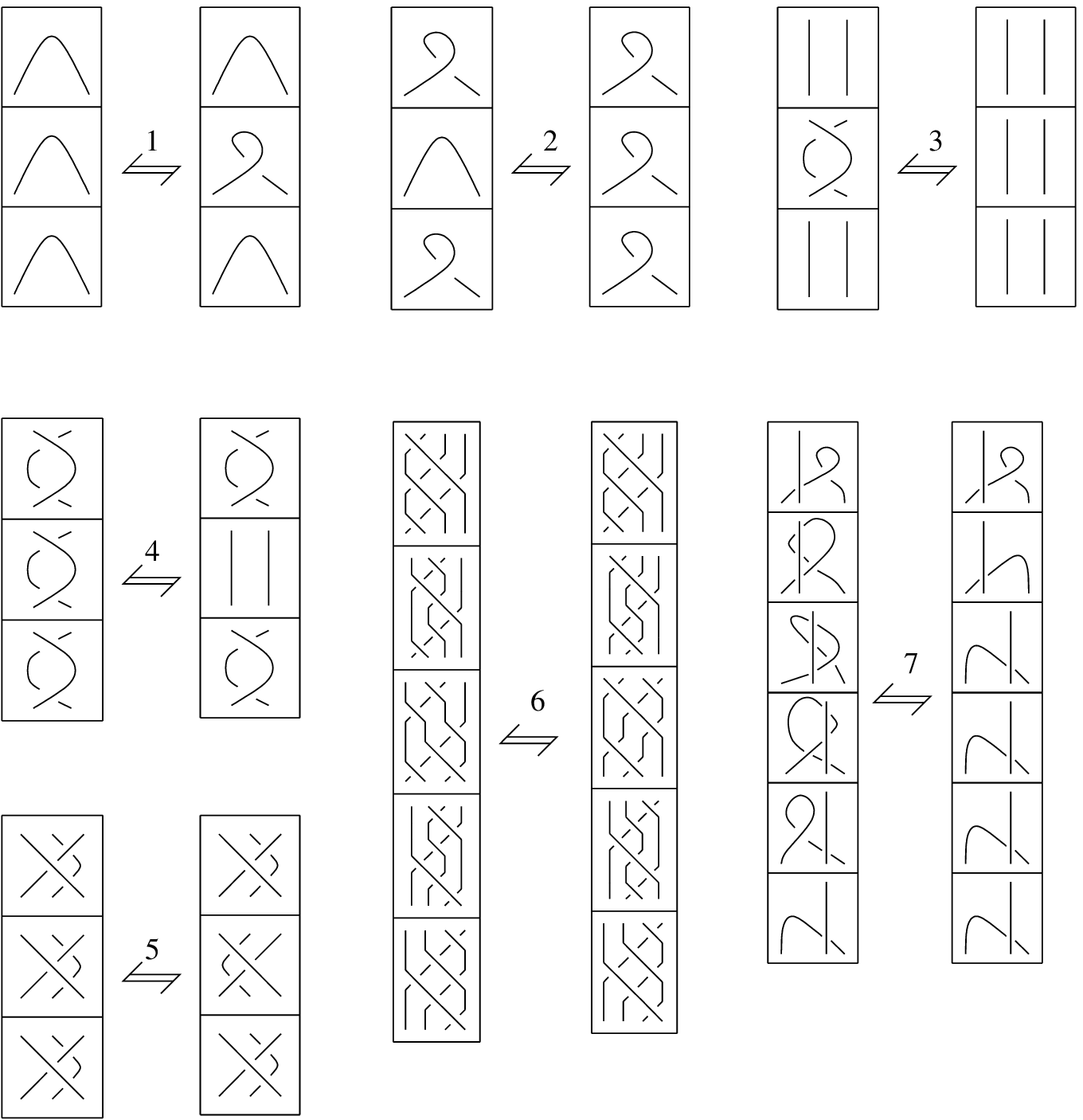}
 \caption{Movie moves 1-7} \label{mm1}
 \end{figure}

\begin{figure} [htb] \drawing{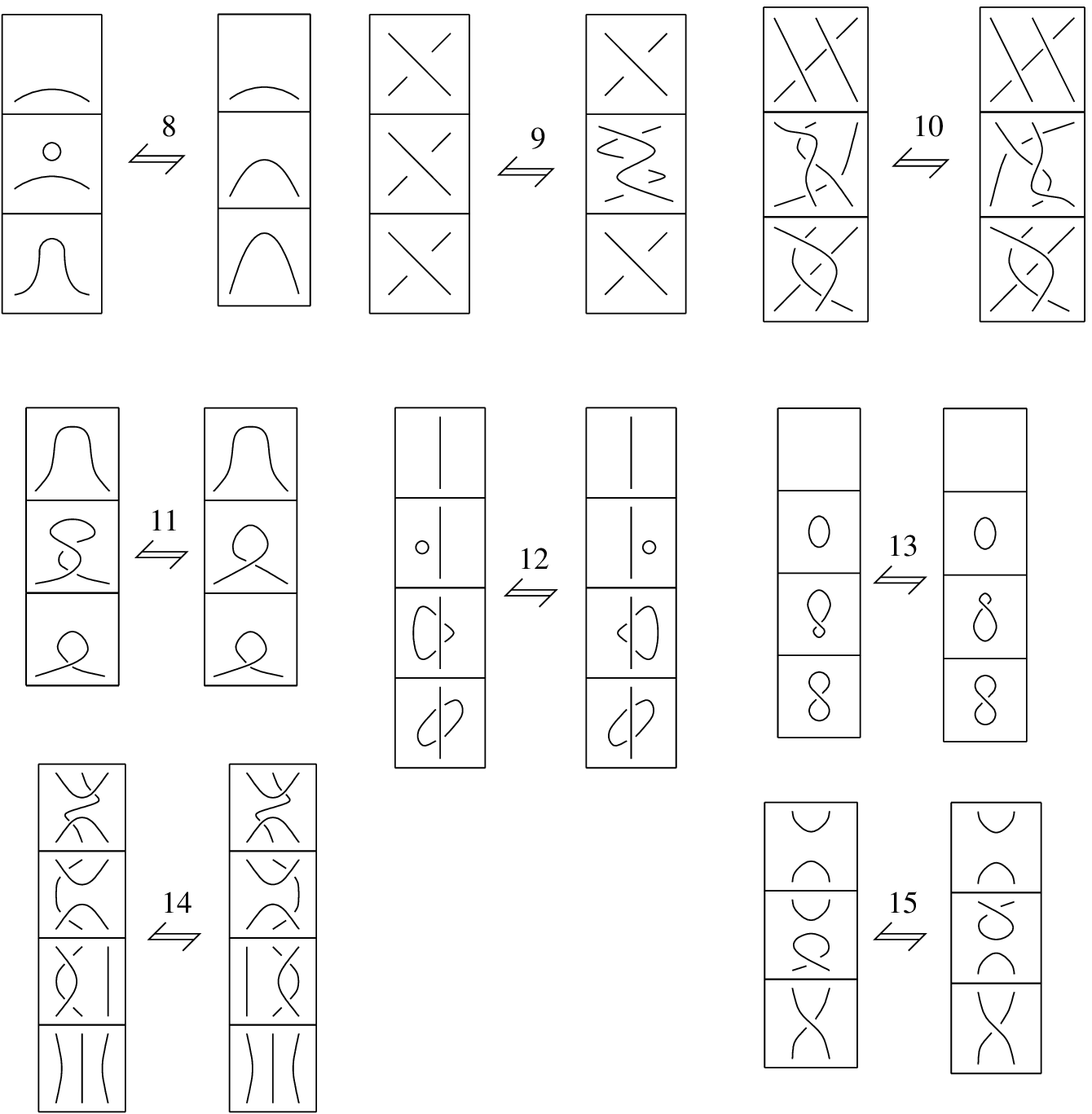}
 \caption{Movie moves 8-15} \label{mm2}
 \end{figure}

A cobordism $S$ admits a combinatorial description via a sequence of plane diagrams of 
its cross-sections with $B^3\times k,$ for various $k\in [0,1].$ Each consecutive pair 
of diagrams differ either by a Reidemeister move, or a Morse move, the latter describing 
a simple cobordism with one critical point (figure~\ref{morse}). Such sequences are 
referred to as \emph{movies}. 
Two sequences describe the same cobordism if the can be connected through a finite 
sequence of \emph{movie moves}, shown in figures \ref{mm1}, \ref{mm2} at the end of the paper. 
We assume that the reader is familiar with this theory, and refer to [CS1], [CS2] and references 
therein for details.  

Let $w=w(I,s)$ be the potential defined in the preceeding section. 
Earlier, we associated an isomorphism $C(D_1) \cong C(D_2)$ in the category 
$K_w$ to a Reidemeister move between tangle diagrams $D_1$ and $D_2.$  Aping the 
last section, assign the maps $\iota, \varepsilon, \eta$ to the circle creation, circle annihilation, and 
saddle point moves (see figure~\ref{morse}), respectively. 

Given a movie $z=(z_0,\dots, z_m)$ representing surface $S,$ assign to 
$z$ the morphism in $K_w$: 
   $$ C(z)\hspace{0.1in} : \hspace{0.1in} C(z_0) \lra C(z_m)\bracket{\upsilon_1(S)}\{ \upsilon_2(S)\} $$
 by composing the maps associated to each move $z_i\to z_{i+1}$ in $z.$ The quantities 
$\upsilon_1(S),\upsilon_2(S)$ were defined at the end of the previous section.

\begin{prop} \label{cob-inv-prop} 
Up to overall multiplication by nonzero rational numbers, the map $C(z)$ (viewed 
as a morphism in the category $K_w$) does not depend on the movie presentation $z$ 
of $S.$ 
\end{prop} 

\emph{Proof:} For a given movie move in figures~\ref{mm1}, \ref{mm2}, denote 
the top frame by $b_1,$ the bottom frame by $b_2,$ the left movie by $S_l$ and 
the right movie by $S_r.$ We need to show that morphisms $C(S_l)$ and $C(S_r)$ 
are proportional, $C(S_l)=\lambda C(S_r)$ for some $\lambda\in \Q^{\ast}.$ 

Movies $S_l$ and $S_r$ in movie move $6$ are compositions of Reidemeister moves. 
Therefore, 
 $$C(S_l), C(S_r)\hspace{0.1in} : \hspace{0.1in} C(b_1)  \lra C(b_2)$$ 
are isomorphisms. Let $b$ be the tangle diagram which consists of four parallel 
 disjoint segments (the crossingless diagram of the trivial braid). 
There is an automorphism of $K_w$ (tensoring with $ C(b_1'),$ 
 where $b_1'$ is the "inverse" of the braid diagram $b_1$) which 
takes $C(b_1)$ to $C(b).$ Therefore, 
  $$\Hom_{K_w}(C(b_1), C(b_1)) \cong \Hom_{K_w}(C(b),C(b)).$$ 
The vector space on the right hand side is isomorphic to $\Q$ (the only degree $0$ 
endomorphisms of $C(b)$ are multiples of the identity). Hence 
  $$\Hom_{K_w}(C(b_1),C(b_2)) \cong \Hom_{K_w}(C(b_1),C(b_1)) \cong \Q $$
(the isomorphisms are not canonical, though), and $C(S_l), C(S_r)$ are proportional. 

This argument applies to moves $1,$   $2,$  $3,$   $4,$  $5,$   $7,$  $9,$  $10,$   $11$ as well, 
and to all versions
of these moves (various orientations,  overcrossing/undercrossing variations, etc.)  

The move $8$ follows from the compatibility of $\iota, \varepsilon,$ and $\eta,$ see 
section~\ref{sec-corners}. 

In move $12,$ the maps 
 $$C(S_l), C(S_r) \hspace{0.1in} : \hspace{0.1in} C(b_1) \lra C(b_2)\bracket{1}\{n-1\}$$ 
are non-trivial and lie in the one-dimensional $\Q$-vector space 
$$\Hom_{K_w}(C(b_1), C(b_2)\bracket{1}\{ n-1\}). $$ 
Therefore, the maps are proportional. Same argument works for other versions of this 
move, and for all versions of moves $13,14,15$ (for moves $14, 15$ change $\{n-1\}$
to $\{1-n\}$)  Proposition~\ref{cob-inv-prop} follows. $\square$ 

\vspace{0.1in} 

The above proof is based on the observation that the space of homs  
between $C(b_1)$ and $C(b_2)$ (with a suitable shift) is one-dimensional. The same 
approach was used in [Kh4] to show functoriality of the homology theory $\mc{H}$ from 
[Kh1] (see Jacobsson [J] for a different proof), and by Dror Bar-Natan [BN2] to prove 
functoriality of his refinement of $\mc{H}.$ 

We denote by $C(S)$ the set $\{ \lambda C(z)| \lambda \in \Q^{\ast} \}.$ This is an 
invariant of the cobordism $S.$ 

In particular, given a diagram $D$ of a tangle $L,$ the object $C(D)$ of $K_w$ is 
canonically (up to rescaling) associated to $L.$ We denote this object by $C(L)$ (recall 
that it's a complex of graded matrix factorizations, treated as an object of $K_w$). 

Given a diagram $D$ of an oriented link $L,$ homology groups $H_n(D)$ are $\Q$-vector 
spaces canonically
assigned to $L$ (up to overall rescaling by nonzero rational numbers). We denote these 
groups by $H_n(L)$ and their graded summands by $H_n^{i,j}(L).$  

Thus, an oriented link cobordism $S$ between $L_0$ and $L_1$ induces a homomorphism 
  $$ C(S) \hspace{0.1in} : \hspace{0.1in} H_n(L_0) \lra H_n(L_1) , $$ 
well-defined up to rescalings by nonzero rationals, and for each $i,j$ restricts to 
homomorphisms 
   $$ H_n^{i,j}(L_0) \lra H_n^{i,j+(1-n)\chi(S)}(L_1),$$
where $\chi(S)$ is the Euler characteristic of $S.$  
The collection of maps $C(S),$ over all link cobordisms $S,$ is a functor from the 
category of link cobordisms to the category of bigraded $\Q$-vector spaces (with morphisms 
being graded linear maps with the equivalence relation $f\sim\lambda f,$ for $\lambda\in \Q^{\ast}$). 
The Euler characteristic of $H_n(L)$ is the polynomial $P_n(L)$: 
  $$P_n(L) = \sum_{i,j} (-1)^i q^j \dim_{\Q} H_n^{i,j}(L).$$ 

\vspace{0.1in}

{\bf 2-functor.} So far we considered tangles with a fixed oriented 
boundary $(I,s).$ By switching from tangles in a ball to tangles in $\R^2\times [0,1]$ 
and varying possible boundaries one can form the 2-category of tangle cobordisms 
(see [F], [BL], and references therein).  Our construction can be  extended to a 2-functor 
from the 2-category of oriented tangle cobordisms to a 2-category with potentials $w(I,s)$ as 
objects, complexes of matrix factorizations as 1-morphisms, and homomorphisms 
between (suitably shifted) complexes as 2-morphisms (of course, we'll have to 
quotient by null-homotopic morphisms, and identify morphisms that are multiples 
of each other, $f\cong \lambda f,$ for $\lambda\in \Q^{\ast}$).  This 2-functor is 
braided monoidal.

\section{A generalization} 
  
Each complex simple Lie algebra $\mf{g}$ gives rise to a polynomial 
invariant of links whose components are decorated by finite-dimensional 
irreducible representations of $\mf{g},$ see [RT]. The polynomial $P_n$ results if 
$\mf{g}=\mf{sl}_n$ and 
every component is assigned the fundamental $n$-dimensional representation $V$
of $\mf{sl}_n.$ Murakami, Ohtsuki, and Yamada [MOY] develop a calculus 
of trivalent graphs that helps in understanding the polynomial invariant of links with 
components colored by arbitrary exterior powers of $V$  (the $i$-th exterior power 
$\Lambda^iV$ is often called the $i$-th  fundamental representation of $\mf{sl}_n$). 
Although their construction also extends the invariant to
spacial trivalent graphs, here we will only look at planar graphs.   
Each edge of a graph is oriented and labelled by a number from $1$ to $n-1.$ 
Every vertex is trivalent and the sum of labels at edges entering
the vertex minus the sum of labels at edges leaving the vertex is a multiple of $n.$ 

\begin{figure} [htb] \drawing{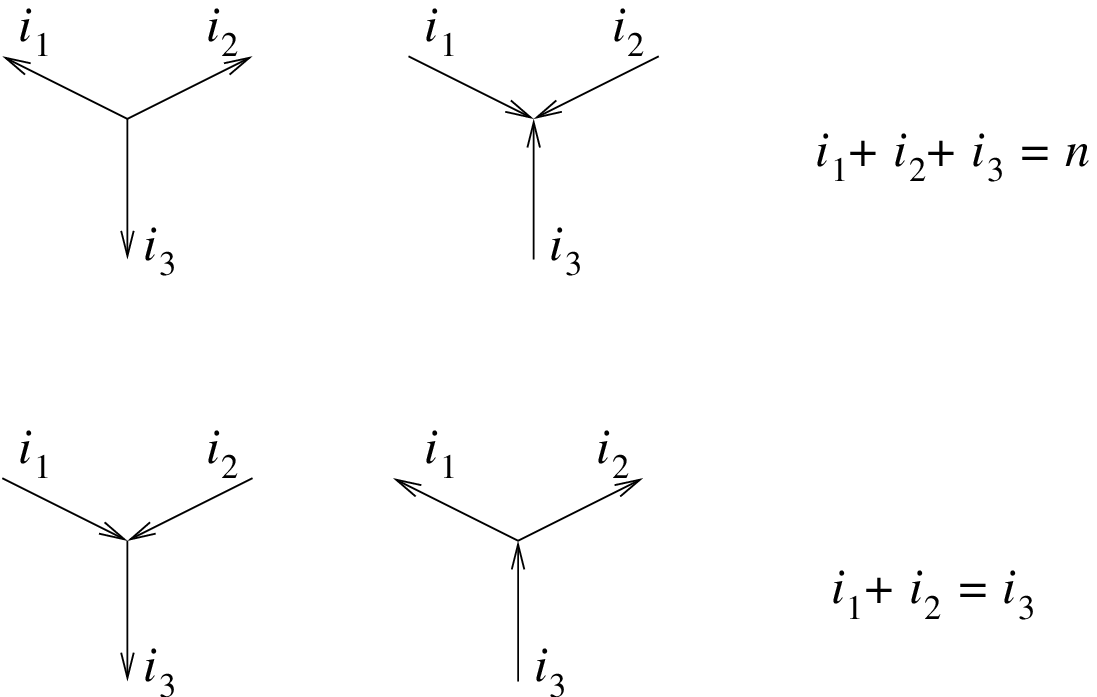}
 \caption{The four types of vertices} \label{orient4}
 \end{figure}

To every such graph $\Gamma$ an invariant $\bracket{\Gamma}$  is assigned, 
taking values in $\Z[q,q^{-1}].$ The invariant is unchanged if the orientation of 
an edge is reversed simultaneously with changing the label from $i$ to $n-i.$ 
We use this transformation to reduce our consideration to graphs with all labels at most 
$\frac{n}{2}.$ We split the vertices into four types by the number (zero to three) 
of edges oriented into the vertex, see figure~\ref{orient4} in the back of the paper.  

The construction of Section 6  generalizes to a homology theory 
$H(\Gamma)$ for graphs $\Gamma$ as above. We will now sketch this generalization 
and its conjectural extension to link homology. For convenience, assume that $n$ is even. 

Denote by $i(e)$ the number assigned to an edge $e,$ and select a set 
$s(e)$ of cardinality $i(e)$ such that the sets assigned to different 
edges are disjoint.  Suppose that a vertex $v$ bounds the edges $e_1,e_2,e_3.$ 
To $v$ we assign the potential 
 $$w_v = \sum_{j\in s(e_1) \sqcup s(e_2) \sqcup s(e_3)} \pm x_j^{n+1} $$ 
where the sign is $+$ if  edge $e_k$ leaves $v$ and $j\in s(e_k),$ and $-$ otherwise
(note that the potential of an edge $e$ is $\pm \sum_{j\in s(e)} x_j^{n+1}$).  

Let  $\Q[s(e)]$ be the ring  of polynomials in variables $x_j,j\in s(e),$ and 
denote by $S(e)$ its subring of symmetric polynomials. Let  
$$R_v \define S(e_1) \otimes_{\Q} S(e_2) \otimes_{\Q} S(e_3)$$ 
be the tensor product of the three rings. Then $w_v\in R_v.$ 

Consider first the case when the edges $e_1,e_2,e_3$ are oriented away from $v.$ 
Then $i(e_1) + i(e_2) + i(e_3) =n.$ 
Let $\sigma_k(v)$ be the $k$-th elementary symmetric function 
in $x_j$'s,  for $j$ in the set $s(e_1) \sqcup s(e_2) \sqcup s(e_3).$ Write 
 $$ w_v =\sum_j x_j^{n+1}= \sum_{k=1}^n    \sigma_k(v) g_k(v) $$ 
for some $g_k(v)$'s (which are not uniquely defined, of course).  

To $v$ we assign the factorization $C_v$ which is the tensor product of 
 $$ R_v \xrightarrow{g_k(v)} R_v \xrightarrow{\sigma_k(v)} R_v$$ 
over $1\le k\le n.$ 

If the edges $e_1,e_2,e_3$ are all oriented towards $v,$ then 
 $$  w_v =-\sum_j x_j^{n+1}=- \sum_{k=1}^n    \sigma_k(v) g_k(v), $$ 
and to $v$ we assign the factorization $C_v$ which is the tensor product 
of 
 $$ R_v \xrightarrow{g_k(v)} R_v \xrightarrow{-\sigma_k(v)} R_v$$ 
over $1\le k\le n.$ 

Suppose now that, say, $e_1,e_2$ are oriented away from $v$ and 
$e_3$ towards $v.$ Then $i(e_3)=i(e_1)+i(e_2)$ and 
 $$  w_v =\sum_{j\in s(e_1)\sqcup s(e_2)} x_j^{n+1}  -  
   \sum_{j\in s(e_3)} x_j^{n+1}.$$ 
Let $\sigma_k'$ be the $k$-th elementary symmetric function 
in variables $x_j,$ over $j\in s(e_1) \sqcup s(e_2),$ and 
$\sigma_k''$ be the $k$-th elementary symmetric function 
in variables $x_j,$ over $j\in s(e_3).$ Since $w_v$ belongs to 
the ideal of $R_v$ generated by the differences $\sigma_k'-\sigma_k'',$  
we can write 
 $$ w_v = \sum_{k=1}^{i(e_3)} (\sigma_k'-\sigma_k'')g_k(v)$$ 
 for some $g_k(v) \in R_v,$ and $1\le k \le i(e_3).$ 
Now assign factorization $C_v,$ the tensor product of 
$$ R_v \xrightarrow{g_k(v)} R_v \xrightarrow{\sigma_k'-\sigma_k''} R_v,$$ 
for $1\le k \le i(e_3),$ to the vertex $v.$ 

Likewise, assign the tensor product of 
$$ R_v \xrightarrow{g_k(v)} R_v \xrightarrow{\sigma_k''-\sigma_k'} R_v$$ 
to $v$ if two edges $e_1,e_2$ are oriented in, and $e_3$  out. 

If $\Gamma$ does not contain loops, define $H(\Gamma)$ as the 
homology of the 2-complex $C(\Gamma)$ which is the tensor product 
of $C_v,$ over all vertices $v$ of $\Gamma.$ The tensor product is 
taken over intermediate rings $S(e),$ for various edges $e,$ so that 
$C(\Gamma)$ is a finite-rank module over the ring $\otimesop{e} S(e).$ 

When loops are present in $\Gamma,$  a loop labelled $i$ will ``contribute'' to the 
tensor product $C_v$ 
the cohomology of the Grassmannian $\mathrm{Gr}(i,n)$ of $i$-dimensional 
subspaces in $\C^n.$

We conjecture that $H(\Gamma)$ has cohomology in only one out of its
two degrees. An additional $\Z$-grading on $H(\Gamma)$ comes from grading 
on the rings $S(e),$ with  $\deg(x_j)=2.$  We conjecture that the graded dimension 
of $H(\Gamma)$ is the invariant $\bracket{\Gamma}$ (up to obvious 
normalizations such as multiplication by a power of $q,$ etc.).  

The above construction should work when $n$ is even. If $n$ is odd, add a 
variable $x_v$ to each vertex with all edges oriented in, and to 
each vertex with all edges out. Then add $x_v^2$ to the potential $w_v,$ enlarge 
the ring $R_v$ by adjoining variable $x_v,$ and modify the factorization $C_v$ 
correspondingly. We omit the details.  

Given a diagram $D$ of an oriented framed link $L$ colored by numbers from $1$ to $n-1,$ 
each crossing can be resolved in a number of ways into planar graphs, and 
the invariant $P_n(L)$ is a linear combination of $\bracket{\Gamma}$ for 
various resolutions $\Gamma$ of $D,$ with coefficients which are plus or minus 
powers of $q,$ see Section 5 of [MOY]. We conjecture that, likewise, 
$\Q$-vector spaces $H(\Gamma)$ can be naturally strung together into a complex $C(D)$ 
whose (bigraded) cohomology groups are invariants of $L$ and have graded 
Euler characteristic $P_n(L).$ For each $n,$ this homology theory of colored 
oriented framed links in $\R^3$ should be functorial under (oriented, framed, and
suitably decorated) link cobordisms, should extend to a homology theory of 
(decorated) spacial trivalent graphs, and be functorial under (carefully defined) 
graph cobordisms in $\R^3 \times [0,1].$


\section{References} 

[AGV] V.~I.~Arnold, S.~M.~Gusein-Zade and A.~N.~Varchenko, 
Singularities of differentiable maps, vol.I, Monographs in 
Mathematics {\bf 82}, Burkh\"auser, Boston, 1985. 

[BL] J.~Baez and L.~Langford, Higher-Dimensional Algebra IV: 2-Tangles, 
 arXiv math.QA/9811139. 

[BN1] D.~Bar-Natan,   On Khovanov's categorification of the Jones polynomial,
 \emph{Algebraic and Geometric Topology,}  2 (2002) 337-370, math.QA/0201043.

[BN2] D.~Bar-Natan, Khovanov's homology for tangles and cobordisms, in 
preparation. 

[Be] D.~J.~Benson, Representations and cohomology I. 
Basic representation theory of finite groups and associative 
algebras, {\it Cambridge studies in advanced mathematics} {\bf 30}, Cambridge U 
Press, 1995. 

[BV] B.~Blok and A.~N.~Varchenko, Topological conformal field 
theories and the flat coordinates, {\it Int. J. Mod. Phys.} {\bf A7} (1992),
no. 7, 1467--1490. 
 
[B] R.-O.~Buchweitz, Maximal Cohen-Macaulay modules and 
Tate cohomology over Gorenstein rings, preprint circa 1986. 
 
[BEH] R.-O.~Buchweitz, D.~Eisenbud, and J.~Herzog, Cohen-Macaulay 
modules on quadrics (with an appendix by R.-O.~Buchweitz), in 
\emph{Lecture Notes in Mathematics} {\bf 1273}, 58--116, 1987. 

[BGS] R.-O.~Buchweitz, G.-M.~Greuel and F.-O.Schreyer, 
Cohen-Macaulay modules on hypersurface singularities II, 
\emph{Invent. Math.} 88, 165--182, 1987. 
 
[CS1] J.~S.~Carter and  M.~ Saito, Reidemeister moves for surface isotopies and their
 interpretation as moves to movies, \emph{Journal  of Knot Theory and its Ramifications,} {\bf 2},
 3, 251--284, 1993. 

[CS2] J.~S.~Carter and M.~Saito, Knotted surfaces and their diagrams, 
 \emph{Mathematical surveys and monographs, vol.55,} AMS, 1998. 

[Di] A.~Dimca, Topics on real and complex singularities, 
 {\it Vieweg Advanced Lectures in Mathematics,} 1987. 

[D] B.~Dubrovin, Geometry and analytic theory of Frobenius 
manifolds, Proceedings of the International Congress of 
Mathematicians, Vol. II (Berlin, 1998).  Doc. Math.  1998,  
Extra Vol. II, 315--326, arxiv  math.AG/9807034.  

[E1] D.~Eisenbud, Homological algebra on a complete intersection, 
with an application to group representations, {\it Trans. Amer. Math. Soc.} 
{\bf 260} (1980), 35--64.

[E2] D.~Eisenbud, Commutative algebra, with a view toward algebraic 
geometry,  Graduate Texts in Mathematics, 150. Springer-Verlag, 
New York, 1995. 

[EP] V.~Ene and D.~Popescu, Rank one maximal Cohen-Macaulay 
modules over singularities of type $Y_1^3+Y_2^3+Y_3^3+Y_4^3,$ 
arXiv math.AC/0303151. 
 
[F] J.~E.~Fischer, Jr., 2-categories and 2-knots, {\it Duke Math. J.} 75 (1994), 
 no. 2, 493--526. 

[FKS] I.~B.~Frenkel, M.~Khovanov, and O.~Schiffmann, Homological 
realization of Nakajima varieties and Weyl group actions, 
arXiv math.QA/0311485. 

[GH] P.~Griffiths and J.~Harris, Principles of algebraic 
geometry, 1978. 

[HOMFLY] P.~Freyd, D.~Yetter, J.~Hoste, W.~B.~R.~Lickorish, 
 K.~Millett and A.~Ocneanu, A new polynomial invariant of knots 
and links, {\it Bull. AMS (N.S.)} {\bf 12}, 2, 239--246, 1985. 

[HP] J.~Herzog and D.Popescu, Thom-Sebastiani problems for maximal 
Cohen-Macaulay modules, {\it Math. Ann.} {\bf 309} (1997), no.4, 677--700. 

[J] M.~Jacobsson, An invariant of link cobordisms from Khovanov's homology 
theory, arXiv math.GT/0206303. 

[KL1] A.~Kapustin and Y.~Li, D-Branes in Landau-Ginzburg models and 
algebraic geometry, arXiv hep-th/0210296.  

[KL2] A.~Kapustin and Y.~Li, Topological correlators in 
Landau-Ginzburg models with boundaries, arXiv hep-th/0305136. 

[KL3] A.~Kapustin and Y.~Li,  D-branes in topological minimal 
models: the Landau-Ginzburg approach, arXiv hep-th/0306001. 

[KS] L.H.~Kauffman and H.~Saleur, Free fermions and the 
Alexander-Conway polynomial,  Comm. Math. Phys.  141  (1991),  
no. 2, 293--327. 

[Kh1] M.~Khovanov, A categorification of the Jones polynomial, 
{\it Duke Math J.} {\bf 101}, 3, 359--426, 1999, arXiv math.QA/9908171. 

[Kh2] M.~Khovanov, A functor-valued invariant of tangles, 
{\it Algebraic and Geometric Topology} {\bf 2} (2002), 665--741, 
 arXiv math.QA/0103190.  

[Kh3] M.~Khovanov, sl(3) link homology I, math.QA/0304375. 

[Kh4] M.~Khovanov, An invariant of tangle cobordisms, math.QA/0207264. 

[Kn] H.~Kn\"orrer, Cohen-Macaulay modules on hypersurface 
singularities. I, {\it Invent. Math.}  {\bf 88}  (1987),  no. 1, 153--164.

[Ku] V.~Kulikov, Mixed Hodge structures and singularities, 
Cambridge U Press, 1998. 

[M] E.~Martinec, Algebraic geometry and effective Lagrangians, 
 {\it Phys. Lett.} B217 (1989) 431--437. 

[MC] D.~McDuff, D.~Salamon, J-holomorphic curves and quantum 
cohomology, University Lecture Series {\bf 6}, 1994, AMS. 

[MOY] H.~Murakami, T.~Ohtsuki and S.~Yamada, HOMFLY polynomial 
via an invariant of colored plane graphs,  Enseign. Math. (2)  44  
(1998),  no. 3-4, 325--360.
 
[NN] T.~Nakayama, C.~Nesbitt, Note on symmetric algebras, 
{\it Annals of Math.} {\bf 39}, no.3, (1938) 659--668. 

[O] D.~Orlov, Triangulated categories of singularities and D-branes 
in Landau-Ginzburg models, arXiv math.AG/0302304. 

[OS] P.~Ozsv\'ath and Z.~Szab\'o, Holomorphic discs and 
knot invariants, math.GT/0209056. 

[P] D.~Popescu, Cohen-Macaulay representation, in 
{\it Algebra---representation theory (Constanta, 2000),} 249--256, NATO Sci. Ser. II Math. 
Phys. Chem., 28, Kluwer Acad. Publ., Dordrecht, 2001. 

[Ra] J.~Rasmussen, Floer homology and knot complements,
 math.GT/0306378. 

[RT] N.~Reshetikhin and V.~Turaev,  Ribbon graphs and their 
invariants derived from quantum groups, Comm. Math. Phys. 127 
(1990), no. 1, 1--26.

[R] L.~Rozansky, Topological A-models on seamed Riemann surfaces, 
 arXiv hep-th/0305205. 

[Sa] K.~Saito, Period mapping associated to a primitive form, 
 {\it Publ RIMS} {\bf 19} (1983), 1231--1264. 

[Sh] A.~Shumakovitch,  http://www.geometrie.ch/KhoHo/ 

[S] F.-O.~Schreyer, Finite and countable CM-representation type, 
in {\it Singularities, representation of algebras, and vector bundles (Lambrecht, 1985),}  9--34, 
Lecture Notes in Math. {\bf 1273},  Springer, Berlin, 1987.

[VW] C.Vafa and N.Warner, Catastrophes and the classification of 
conformal theories, {\it Phys. Lett.} B218 (1989) 51--58. 

[Y1] Y.~Yoshino, Cohen-Macaulay modules over Cohen-Macaulay rings, 
London Math. Soc. Lect. Note Ser. 146, Cambridge University Press, 
 1990. 

[Y2] Y.~Yoshino,  Tensor products of matrix factorizations,   
 {\it Nagoya Math. J.}  152,   (1998), 39--56.

\vspace{0.1in} 

Mikhail Khovanov, Department of Mathematics, University of California, 
Davis, CA 95616, mikhail@math.ucdavis.edu

\vspace{0.07in} 

Lev Rozansky, Department of Mathematics, University of North Carolina
Chapel Hill, NC 27599, rozansky@math.unc.edu

\end{document}